\documentclass[poms,nonblindrev]{poms1} 

\OneAndAHalfSpacedXI 

\usepackage{graphicx}
\usepackage{epsfig}
\usepackage{natbib}
\def\bibfont{\small}%
\usepackage[english]{babel}
\usepackage[ansinew]{inputenc}
\usepackage{sansmath}
\usepackage{amsmath}
\usepackage{enumitem}
\usepackage{array}[r]
\usepackage{eqparbox}[r]
\usepackage[noend]{algcompatible}

\usepackage{etoolbox}
\usepackage{cases}
\usepackage{amsfonts}
\usepackage{longtable}
\usepackage{graphicx}
\usepackage[caption=false]{subfig}
\usepackage[dvipsnames, svgnames, x11names]{xcolor}
\usepackage{breqn}
\usepackage{bm}
\usepackage{pdflscape}
\usepackage{setspace} 
\usepackage{amsfonts,amssymb} 
\usepackage{multirow}
\usepackage{booktabs}
\usepackage[ruled,noend,norelsize,linesnumbered]{algorithm2e}
\usepackage{amssymb}
\usepackage{etoolbox}
\newtheorem{modl}{Model}

\newenvironment{model}{\begin{samepage}\begin{modl}}{\end{modl}\end{samepage}}
\makeatletter
\AtBeginDocument{\def\@citecolor{blue}}
\AtBeginDocument{\def\@linkcolor{blue}}
\AtBeginDocument{\def\@urlcolor{blue}}
\makeatother
\usepackage[colorlinks=true,citecolor=azure,linkcolor=red]{hyperref}
\usepackage[noend]{algcompatible}
\newcommand{\Fig}{\textcolor{black}{Figure}\space}
\newcommand{\T}{\textcolor{black}{Table}\space}
\newcommand{\eq}{\textcolor{black}{Eq.}}

\newcommand{\Alg}{\textcolor{black}{Algorithm}\space}
\newcommand{\M}{\textcolor{black}{Model}\space}
\newcommand{\s}{\textcolor{black}{Section}\space}

\newcommand{\SDL}{\textcolor{blue}{SD-L}}
\newcommand{\SDQ}{\textcolor{blue}{SD-Q}}
\newcommand{\MPECQ}{\textcolor{blue}{MPEC-Q}}
\newcommand{\MPECL}{\textcolor{blue}{MPEC-L}}
\newcommand{\SLMFGQ}{\textcolor{blue}{SLMFG-Q}}
\newcommand{\SLMFGL}{\textcolor{blue}{SLMFG-L}}

\newcommand{\V}{\textcolor{blue}{3.1}}
\newcommand{\VI}{\textcolor{blue}{3.2}}
\newcommand{\III}{\textcolor{blue}{2.1}}
\newcommand{\IV}{\textcolor{blue}{2.2}}
\newcommand{\I}{\textcolor{blue}{1.1}}
\newcommand{\II}{\textcolor{blue}{1.2}}
\newcommand{\Bard}{$\mathsf{Bard \& Moore}$}
\newcommand{\BP}{$\mathsf{SD}$-$\mathsf{BP}$}
\newcommand{\Diff}{$\mathsf{SD}$-$\mathsf{Diffob}$}
\newcommand{\Weight}{$\mathsf{SD}$-$\mathsf{Wi}$}
\newcommand{\SDBB}{$\mathsf{SD}$-$\mathsf{B\&B}$}

\allowdisplaybreaks
\OneAndAHalfSpacedXI 
\usepackage{natbib}
 \bibpunct[, ]{(}{)}{,}{a}{}{,}%
 \def\bibfont{\small}%
\newcommand{\DA}[1]{\textcolor{black}{#1}}
\newcommand{\DR}[1]{\textcolor{black}{#1}}
\newcommand{\HX}[1]{\textcolor{black}{#1}}
\TheoremsNumberedThrough     
\ECRepeatTheorems
\newtheorem{observation}{Observation}

\begin{document}


\RUNAUTHOR{Xi et al.}

\RUNTITLE{Single-leader multi-follower games for the regulation of two-sided Mobility-as-a-Service markets}

\TITLE{Single-leader multi-follower games for the regulation of two-sided Mobility-as-a-Service markets}

\ARTICLEAUTHORS{%
\AUTHOR{Haoning Xi}
\AFF{School of Civil and Environmental Engineering, UNSW Sydney, NSW, 2052, Australia.\\
Data61, CSIRO, Canberra ACT 2601, Australia.\\} 
\AUTHOR{Didier Aussel}
\AFF{PROMES UPR CNRS 8521, Universit\'{e} de Perpignan Via Domitia, Tecnosud, 66100 Perpignan, France.\\}
\AUTHOR{Wei Liu}
\AFF{School of Civil and Environmental Engineering, UNSW Sydney, NSW, 2052, Australia.\\School of Computer Science and Engineering, University of New South Wales, Sydney, NSW 2052, Australia.\\\URL{}}
\AUTHOR{S.Travis Waller}
\AFF{School of Civil and Environmental Engineering, UNSW Sydney, NSW, 2052, Australia.\\ \URL{}}
\AUTHOR{David Rey}
\AFF{SKEMA Business School, Universit\'e C\^ote d'Azur, Sophia Antipolis Campus, France.\\
School of Civil and Environmental Engineering, UNSW Sydney, NSW, 2052, Australia.\\ \EMAIL{david.rey@skema.edu} \URL{}}
} 

\ABSTRACT{%
Mobility-as-a-Service (MaaS) is an emerging business model driven by the concept of ``Everything-as-a-Service'' and enabled through mobile internet technologies. In the context of economic deregulation, a MaaS system consists of a typical two-sided market, where travelers and transportation service providers (TSPs) are two groups of agents interacting with each other through a MaaS platform. In this study, we propose a modeling and optimization framework for the regulation of two-sided MaaS markets. We consider a name-your-own-price (NYOP)-auction mechanism where travelers submit purchase-bids to accommodate their travel demand via MaaS platform, and TSPs submit sell-bids to supply mobility resources for the MaaS platform in exchange for payments. We cast this problem as a single-leader multi-follower game (SLMFG) where the leader is the MaaS regulator and two groups of follower problems represent the travelers and the TSPs. The MaaS regulator aims to maximize its profits by optimizing operations. In response to the MaaS regulator's decisions, travelers (resp. TSPs) adjust their participation level in the MaaS platform to minimize their travel costs (resp. maximize their profits). We analyze cross-group network effects in the MaaS market, and formulate SLMFGs without and with network effects leading to mixed-integer linear bilevel programming and mixed-integer quadratic bilevel programming problems, respectively. We propose customized branch-and-bound algorithms based on strong duality reformulations to solve these SLMFGs. Extensive numerical experiments conducted on large scale simulation instances generated from realistic mobility data highlight that the performance of the proposed algorithms is significantly superior to a benchmarking approach, and provide meaningful managerial insights for the regulation of two-sided MaaS markets in practice.
}%


\maketitle

%


\section{Introduction}
 
The accelerated evolution of the digital and shared economies in recent years has brought profound shifts in how various services are delivered. In this context, the transportation industry is undergoing a massive revolution from infrastructure-focused  towards  service-focused business models brought by the ``Mobility-as-a-Service (MaaS)'' concept. According to the MaaS Alliance \citep{MaaSalliance}, MaaS integrates various forms of transportation services into a single mobility service accessible on demand. To meet a traveler request, the MaaS operator can provide a diverse combination of mobility options across multiple travel modes, including public transport, ride-, car- or bike-sharing, taxi; hereby referred to as a \textit{MaaS bundle}. \HX{MaaS is a user-centric framework where customized services and mode-agnostic mobility resources are priced in a unified framework. Yet, in the vast majority of studies on MaaS systems, mobility resource pricing is based on segmented travel modes.  \cite{xi2020incentive} made an initial attempt to address this research gap by introducing innovative MaaS mechanisms where users can bid for any quantity of mobility resources in a mode-agnostic fashion with preference requirements, and the MaaS regulator allocates the mobility resources integrated from various transportation service providers (TSPs) to meet heterogeneous traveler requests. In this MaaS framework, the mobility service provided by different travel modes are unified as \textit{mobility resources} with continuous provision by capturing the travel distance and the average service speed.} Building on this flexible MaaS framework, we consider a MaaS system under economic deregulation \citep{wong2020mobility}, which admits a natural two-sided market representation wherein a MaaS regulator aims to promote the participation of both travelers and TSPs. The operation of a two-sided MaaS market is challenging due to the interactions among different stakeholders, i.e. the regulator, travelers and TSPs which may have divergent and conflicting objectives \citep{meurs2017mobility}. To capture the interactions across stakeholders of a two-sided MaaS market, we propose a single-leader multi-follower game (SLMFG) where the MaaS regulator is the leader, travelers and TSPs are two groups of followers. In the proposed SLMFG, the MaaS regulator makes operational decisions, i.e., price, \textit{MaaS bundles}, to allocate mobility resources to travelers, maximizing its profits by anticipating the participation of followers. Each traveler aims to minimize her travel costs by deciding the participation level defined as the proportion of mobility demand fulfilled through the MaaS platform in comparison with outside options. Each TSP aims to maximize its profits by deciding the participation level defined as the proportion proportion of mobility resources supplied to the MaaS platform in comparison with its reserve options. Two-sided markets represent a refinement of the concept of cross-group network effects, which is captured by the supply-demand gap in this study. We formulate  SLMFGs  without  and  with  network  effects  which  lead  to  mixed-integer linear bilevel programming (MILBP) and mixed-integer quadratic bilevel programming (MIQBP) problems with multiple constraints representing followers' strategies, respectively. Then we propose exact solution methods to find optimal solutions for both MILBP and MIQBP problems.

We next review the literature on two-sided markets in flexible mobility systems (Section \ref{Two-sided}), single-leader multi-follower games (Section \ref{Single-leader}), solution approaches of bilevel programming problems (Section \ref{Solution approaches}) and outline the main contributions of this paper (Section \ref{contribution}).

\subsection{Two-sided markets in flexible mobility systems} 
\label{Two-sided}

A two-sided market is one in which two groups of agents interact through a platform, and the decisions of each group of agents affect the outcomes of the other group of agents, typically through an externality \citep{rysman2009economics}. The goal of the platform regulator is to incentivize the participation of both sides of the market \citep{rochet2003platform}. Two-sided platforms have become omnipresent in literature on flexible mobility systems, especially for pricing in ridesourcing markets. \cite{wang2016pricing} introduced matching mechanisms to model the cross-group externalities between customers and taxi drivers on an e-hailing platform and reveal the impacts of the pricing strategy of an e-hailing platform in taxi service. \cite{djavadian2017agent} explored MaaS in the framework of two-sided markets based on agent-based stochastic, day-to-day adjustment processes using Ramsey's pricing criterion for social optimum.  \cite{bai2019coordinating} examined how various factors affect the optimal price, wage, and commission by considering heterogeneous earning-sensitive drivers and price-sensitive passengers. \cite{nourinejad2020ride} introduced a dynamic model in a two-sided ridesourcing market such that supply and demand are endogenously dependent on the fare charged from riders, the wage paid to drivers, riders' waiting time and drivers' cruising time. \cite{he2021beyond} proposed an integrated model of two-sided platforms
by addressing the joint design of incentives and spatial capacity allocations in shared-mobility systems. 
\subsection{Single-leader multi-follower games}
\label{Single-leader}

Single-leader multi-follower games (SLMFGs) can be viewed as a particular class of hierarchical decision problems belonging to the class of bilevel optimization \citep{dempe2002foundations}. SLMFGs have a various of practical applications in \HX{operations and management}. \HX{\cite{wu2016bi} developed a leader-follower game model to evaluate the potential gains from the merger of different organizations with constrained resources and discussed its applications in banking operations. \cite{zha2018surge} investigated the performance of surge pricing in ridesharing system where the leader's problem determines the prices and the followers' problem represents drivers' equilibrium in work hour choices. \cite{rey2019branch} considered the network maintenance problem, where the leader problem describes the manager's policy problem and each lower level problem represents a user equilibrium problem per time period. \cite{xiongrobust} investigated a resource recovery planning problem where the local authority (leader) solves a joint waste and cost sharing problem to minimize its expected payout, and each private operator (follower) solves a facility location and resource recovery operations problem to maximize its expected total
profits.}
\subsection{Solution approaches for bilevel optimization problems}
\label{Solution approaches}

Solving bilevel optimization problems is a challenging task since even the linear single-leader single-follower problems are NP-hard \citep{bard1991some}. The most common solution approach for linear bilevel programming (LBP) problems is to replace the lower-level problem with its Karush-Kuhn-Tucker (KKT) conditions, leading to a Mathematical Program with Equilibrium Constraints (MPEC) \citep{luo1996mathematical}. \cite{fortuny1981representation} proposed to reformulate the complementary slackness condition by introducing a large positive constant $M$ and a binary variable, thus converting LBPs into large mixed-integer programming problems. This approach has been studied in several papers, e.g., \cite{luathep2011global,dempe2012karush}. However, since the MPEC reformulation of LBPs require binary variables for each complementary slackness condition, this may affect computational scalability. To address this drawback, an alternative method is to employ the strong-duality (SD) property of linear programming to reformulate a linear lower-level using a primal-dual inequality \citep{pandzic2012epec,beheshti2016exact}. \cite{zare2019note} compared the performance of SD-based reformulations against KKT-based ones and showed that the SD-based approaches could reduce computational runtime by several orders of magnitude for certain classes of instances. In contrast to the MPEC reformulation, the number of additional binary variables in SD-based reformulations does not depend on the size of the follower's problem. Further, if the leader variables interacting with the follower problem are integers, the resulting nonlinear constraints induced by the SD reformulation can be linearized by applying classical linearization approaches for integer-linear bilinear terms such as so-called big-$M$ approaches \citep{glover1974converting,liberti2009reformulations}. However, \cite{pineda2019solving} raised concerns about the widespread  use of big-$M$ approaches to solve LBP and showed that the heuristic procedure employed to select the big-$M$ in many published works may actually fail and provide suboptimal solutions. \cite{kleinert2020there} considered the hardness to find a correct big-$M$ for LBP and showed that identifying such value is NP-hard.

To obviate the aforementioned issues induced by follower-problem mixed-integer reformulations, exact algorithms based on branch-and-bound (B\&B) have been developed \citep{falk1969algorithm,bard1982explicit}. \cite{bard1990branch} developed a B\&B algorithm based on the MPEC reformulation, referred to as \Bard\space in this paper, which works by relaxing nonlinear constraints corresponding to complementary slackness conditions and branching on these constraints. However, this approach may require a substantial amount of branching on complementary conditions which may be impractical on large problems. Recently, \cite{kleinert2020closing} derived valid inequalities for LBP by exploiting the SD conditions of the follower-level problem, which proves to be very effective in closing the gap for some instances.

The above review of the literature highlights that solving bilevel optimization problem is highly non-trivial, even in the single-follower case. Using big-$M$ approaches may fail and provide suboptimal solutions, thus it is necessary to develop exact solution approaches.
\subsection{Our contributions}
\label{contribution}
In this study, we present a novel approach for the regulation of two-sided MaaS markets. We propose a SLMFG where the MaaS regulator is the leader while travelers and TSPs are represented as two groups of multiple followers. The proposed SLMFG aims to capture the  interaction of all three types of stakeholders (MaaS regulator, travelers and TSPs). The main contributions of this study are summarized as follows:
\begin{itemize}
\item We develop a name-your-own-price (NYOP) auction-based mechanism for two-sided MaaS markets which allows travelers and TSPs to submit heterogeneous MaaS requests, including their preference and willingness to pay/sell in terms of purchase/sell-bids simultaneously. 
\item We formulate the SLMFG revealing the hierarchical interactions between the MaaS regulator (the leader) and the travelers/TSPs (the followers) through the pricing strategies and MaaS bundles, as well as the cross network effects between travelers and TSP by tracking the supply and demand in the two-sided MaaS markets. Specifically, we formulate SLMFGs without and with network effects leading to MILBP  and MIQBP problems, respectively.  
\item   We provide constraint qualifications for MPEC reformulations and theoretically prove the equivalence between proposed SLMFGs and their MPEC reformulations. 
\item We propose a new exact solution approach to solve the proposed SLMFGs that combines SD reformulation with a customized branch-and-bound algorithm, referred to as \SDBB. We conduct extensive numerical experiments which highlight that the proposed algorithm (\SDBB) considerably outperforms the classical B\&B algorithm (\Bard) and draw managerial insights that can be used as operational management guidelines in practice.
\end{itemize}
To the best of  our knowledge, this is the first tractable framework where single-leader multi-follower game approach is used for modelling a two-sided MaaS market \DA{and the first SD-based branch and bound algorithm with customized branching rules for solving both MILBP and MIQBP problems}.

The rest of this paper is organized as follows: Section \ref{PS} gives the problem description; Section \ref{SLMFGs} introduces SLMFGs in the two-sided MaaS market; Section  \ref{Solution methods} introduces the exact solution methods for the SLMFGs; Section \ref{numerical experiments}  presents the computational study; Section \ref{Conclusion} concludes the study and discusses future research directions.  

\section{Problem description}
\label{PS}
This section introduces preliminaries and notations (Section \ref{Preliminaries}), presents the NYOP-auction mechanism (Section \ref{mechanism}) and provides an overview of the SLMFG (Section \ref{overview}).

\subsection{Preliminaries and Notations}
\label{Preliminaries}
\HX{In this study, we take the perspective of a MaaS regulator who aims to make the optimal operational strategies, i.e., pricing policy, and allocate \textit{mobility resources} integrated from TSPs to travelers in terms of multi-modal \textit{MaaS bundles} in a two-sided MaaS system.} In this MaaS framework, the average speed of the travel modes is assumed known and representative of the  \HX{``Commercial speed'' \footnote{\HX{Commercial speed refers to the average speed over all operational stops and is a key factor in the operation of public transport systems since it represents a direct measure of the quality of service and also considerably affects system costs \citep{cortes2011commercial}.}}}, \HX{which incorporates service delays, e.g. waiting time, transfer time, operational stop time as well as the average driving speed.} In addition, to quantify the ``user experience'' of a multi-modal trip, each travel mode is assigned an inconvenience cost per unit of time which aims to capture discomfort in shared mobility modes and is dependent on the vehicle occupancy (the number of shared riders on a vehicle), e.g., inconvenience cost per unit time of taxi is smaller than that of public transit. Travelers can request \HX{total service time} and the corresponding willingness to pay for each trip, they are also required to submit the personal information while registering for the MaaS platform, such as inconvenience tolerance and travel delay budget representing the maximum acceptable inconvenience cost and delay for a mobility service. 

\HX{To quantify the mobility services consisting of different travel modes, we introduce the quantity of \textit{mobility resource} defined as speed-weighted travel distance. Let $D_{i}$ and $T_{i}$ denote Traveler $i$'s requested travel distance and total service time, the corresponding mobility resources are defined as $Q_{i} \triangleq \frac{D_i}{T_i}D_i$, where $\frac{D_i}{T_i}$ denotes Traveler $i$'s requested average service speed. Using this concept, users' MaaS requests are converted into mobility resources which can be fulfilled through a multi-modal MaaS bundle with the required commercial speed (service time).} 

\HX{ Since TSPs are limited in mobility resource capacity, we assume that mobility capacity of the network is known and can be expressed in terms of mobility resources. We assume that heterogeneous users have preferences towards mobility services such as travelling distance, total service time, inconvenience cost, travel delay and willingness to pay (WTP), and heterogeneous TSPs have preferences towards the capacity of mobility resources and willingness to sell (WTS). }

\subsection{NYOP-auction mechanism}
\label{mechanism}
\HX{MaaS aims to deliver a portfolio of multi-modal mobility services that places user experience at the centre of the offer. Thus we propose a name-your-own-price (NYOP) auction-based mechanism, where the user determines the price by submitting a bid and the auctioneer can either accept or reject the user by comparing with an internal threshold price which is invisible for users \citep{spann2004measuring}. 
 \cite{terwiesch2005online} indicated that the NYOP auction provides a niche market where consumers are sensitive to
price or psychologically prefer this type of auction.} \HX{\cite{cai2009optimal} evaluated the optimal reserve prices in the NYOP channel with bidding options in various situations. However, since consumers do not have complete information of the products, the existing NYOP auctions are opaque, and consumers have to trade off their convenience with the low price. In order to solve this issue, travelers and TSPs are allowed to submit their heterogeneous preference and WTP/WTS in terms of purchase/sell-bids in the proposed user-centric MaaS system.} 

Let $\mathcal{I}$ denote the set of travelers, $\mathcal{M}$ denote the set of travel modes and $\mathcal{N}_{m}$ denotes the set of TSPs of mode $m$, $\forall m\in \mathcal{M}$. For each traveler $i \in \mathcal{I}$, the purchase-bid $\mathcal{B}_{i}=(D_{i},T_{i}, b_{i}, B_{i}, R_{i}, \Gamma_{i})$ is defined as a tuple consisting of the following data: $D_{i}$ is Traveler $i$'s requested travel distance, $T_{i}$ is Traveler $i$'s requested total service time, $b_{i}$ is Traveler $i$'s unit purchase-bidding price (WTP), $B_{i}$ is Traveler $i$'s travel expenditure budget, $R_{i}$ is Traveler $i$'s travel delay budget and $\Gamma_{i}$ is Traveler $i$'s inconvenience tolerance. For each TSP $mn$, $m\in \mathcal{M}, n\in\mathcal{N}_m$, its sell-bid $\mathcal{B}_{mn}=(C_{mn},\beta_{mn},\bar{B}_{mn})$ is defined as a tuple, where $C_{mn}$ is the capacity of TSP $mn$ in terms of quantity of mobility resources, $\beta_{mn}$ is TSP $mn$'s unit sell-bidding price (WTS) and $\bar{B}_{mn}$ is TSP $mn$'s operating cost budget.
\begin{figure}[!t]
		\centering
	\caption{NYOP auction-based mechanism in the two-sided MaaS market}
	\vspace{0.1in}
 \includegraphics[width=1\textwidth]{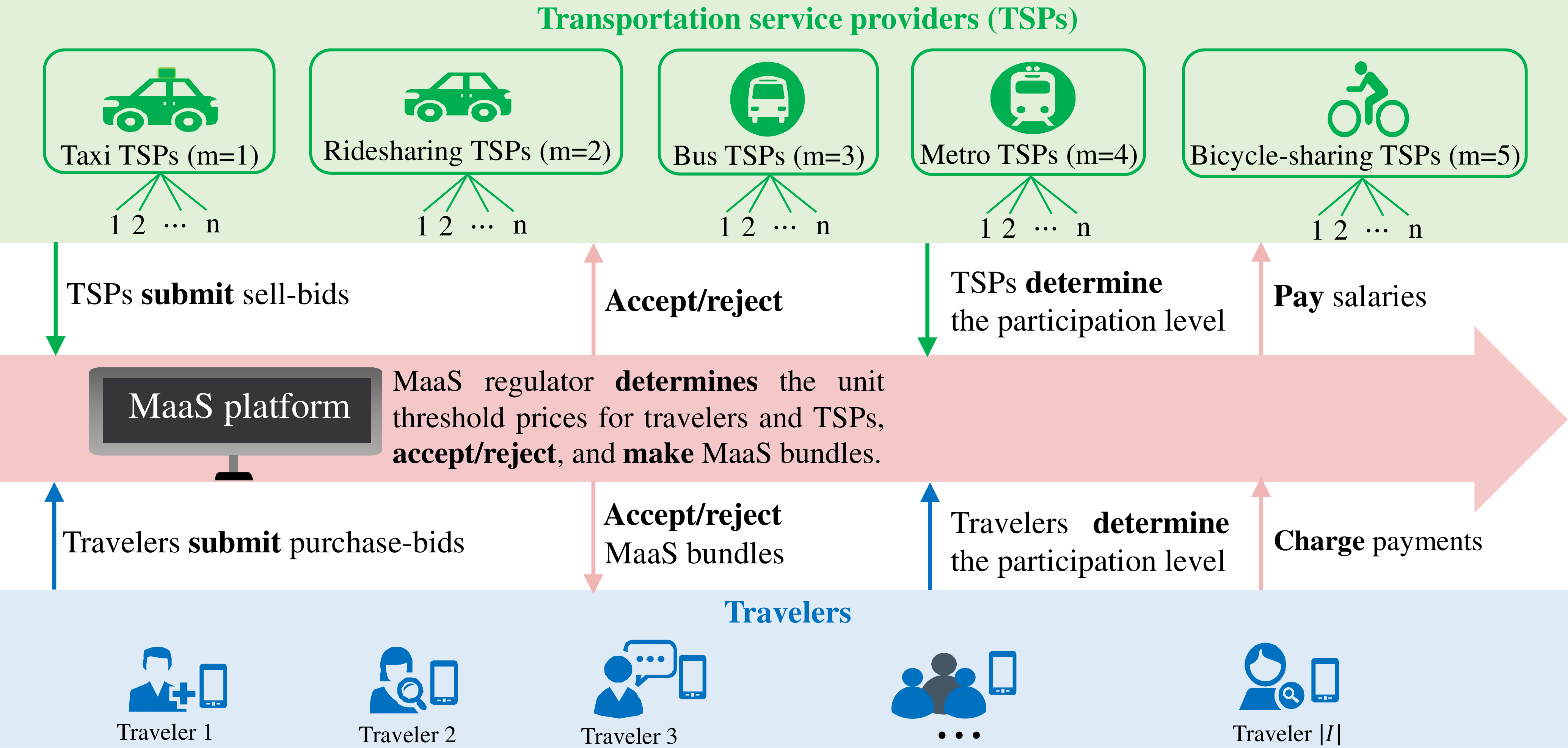}		\label{F1}
\end{figure}

 \HX{Given the purchase/sell-bids of travelers/TSPs, \HX{ The MaaS regulator solves a SLMFG to determine the unit threshold price for travelers and TSPs, respectively, allocate mobility resources to users in terms of MaaS bundles, so as to maximize its profits while coordinating travelers' mobility demand and the mobility resources supplied by TSPs. On the one hand, if a traveler's WTP on per unit mobility is higher than or equal to the unit threshold price for travelers, she will be accepted; otherwise, she will be rejected. On the other hand, if a TSP's WTS on per unit mobility is lower than or equal to the unit threshold price for TSPs, it will be accepted, otherwise, it will be rejected by the MaaS regulator. More details on the SLFMG are provided in Section \ref{overview}. Considering five travel modes, the proposed NYOP auction process is also illustrated  in \Fig\ref{F1}}.}
\subsection{SLMFG model overview}
\label{overview}
We propose a SLMFG where the MaaS regulator is the leader, while travelers and TSPs are two groups of followers. We next introduce the main decision variables of the SLMFG before discussing cross network effects. A summary of the  parameters, variables and sets used in the mathematical models is provided in \T\ref{A} of Online Appendix \ref{EC1}.

\subsubsection{Decision variables.} 
In the proposed NYOP mechanism, we denote $p$ and $q$ the non-negative real variable representing \HX{the unit threshold price} for travelers and TSPs, respectively. Then the MaaS regulator accepts or rejects participants by comparing their bidding prices with the internal threshold unit price. We denote $u_{i} \in \{0,1\}$ the binary variable indicating whether the MaaS regulator accepts traveler $i$, $\forall i\in\mathcal{I}$, to join the MaaS platform (1) or not (0). Analogously, we denote $w_{mn} \in \{0,1\}$ the binary variable indicating whether the MaaS regulator accepts TSP $mn$, $\forall m\in\mathcal{M},n\in\mathcal{M}_{n}$, to join the MaaS platform (1) or not (0).  \HX{In the proposed SLMFG, the travelers/TSPs decide whether to join the MaaS platform or not, and to which extent they wish to use/supply the MaaS platform.} Note that $u_i = 1$ (resp. $w_{mn}=1$) only means that traveler $i$ (resp. TSP $mn$) has the option to join the MaaS platform, i.e., if a traveler is accepted by the MaaS regulator, she may join the MaaS platform with a certain participation level or not join the MaaS platform. However, $u_i=0$ (resp. $w_{mn}=0$) means that traveler $i$ (resp. TSP $mn$) cannot join the MaaS platform.

The allocation of mobility resources supplied by TSPs to travelers is coordinated via real variables $l_{i}^{m} \geq 0$, which represent the \HX{service time} of travel mode $m \in \mathcal{M}$ in the \textit{MaaS bundle} of Traveler $i \in \mathcal{I}$. Cross network effects in the two-sided market are captured by the real variable $\Delta$ which represents the supply-demand gap based on the decisions of followers (more details are provided in Section \ref{cross}). The decision space of the leader is thus the tuple  $(p,q,\bm{u},\bm{w},\bm{l},\Delta)$ where bold-face symbols represent vectors of variables of appropriate dimensions.
 
The decision variables of followers are divided into two groups of variables $\bm{x}=[x_i]_{i \in \mathcal{I}}$ for travelers and $\bm{y} = [y_{mn}]_{m\in\mathcal{M},n\in\mathcal{N}_m}$ for TSPs. For each traveler $i \in \mathcal{I}$, $x_i \in [0,1]$ is a real variable representing the proportion of mobility demand that $i$ decides to use the MaaS platform \HX{in comparison with her other options, e.g., private vehicles and public transit}. Analogously, for each TSP $mn$, $m \in \mathcal{M}, n \in \mathcal{N}_m$, $y_{mn} \in [0,1]$ is a real variable representing the proportion of mobility resources that TSP $mn$ decides to supply the MaaS platform \HX{as opposed to its reserve options}. Hereby $x_i$ and $y_{mn}$ are  referred to as Traveler $i$'s and TSP $mn$'s participation level on MaaS platform, respectively.

\subsubsection{Cross network effects.} \label{cross} Cross network effects aim to capture the interaction between both sides of the MaaS market. E.g., an increase in the mobility resources supplied by TSPs will reduce travelers' average waiting time and may thus increase their travel demand. Reciprocally, a decrease of travelers' mobility demand for MaaS may increase TSPs' average idle time and thus disincentivize participation. We propose to capture these cross network effects through the supply-demand gap in the two-sided MaaS market.

Recall that Traveler $i$'s total demand is $Q_{i}, \forall i\in\mathcal{I}$, and TSP $mn$'s maximum supply is $C_{mn}, \forall m\in\mathcal{M},n\in\mathcal{N}_{m}$. Travelers' total demand on MaaS is written as $\sum_{i \in \mathcal{I}}Q_{i}x_{i}$ and TSPs' total supply is written as $\sum_{m \in  \mathcal{M}}\sum_{n \in \mathcal{N}_{m}} C_{mn}y_{mn}$. The supply-demand gap in the MaaS market is given in \eq\eqref{gap1}:
\begin{equation}
\Delta= \sum_{m \in  \mathcal{M}}\sum_{n \in \mathcal{N}_{m}} C_{mn} y_{mn} - \sum_{i \in \mathcal{I} } Q_{i} x_{i}.\label{gap1}
\end{equation}

The range of $\Delta$ is $[\underline{C},\overline{C}]$, where $\underline{C}$ is the reserved capacity of the MaaS regulator and $\overline{C}$ is the capacity of  mobility resources provided by all TSPs, namely, $\overline{C}=\sum_{m \in  \mathcal{M}}\sum_{n \in \mathcal{N}_{m}} C_{mn}$.

\HX{Let $\Delta_{-i}$ and $\Delta_{-mn}$ denote the supply-demand gap perceived by Traveler $i$ and TSP $mn$, which are defined in \eqref{Deltai} and \eqref{Deltamn}, respectively, }
\begin{equation}
\label{Deltai}
\Delta_{-i}=\sum_{m \in  \mathcal{M}}\sum_{n \in \mathcal{N}_{m}} C_{mn} y_{mn} - (\sum_{i \in \mathcal{I} } Q_{i} x_{i}- Q_{i} x_{i})=\Delta+Q_{i}x_{i},
\end{equation}
\begin{equation}
\label{Deltamn}
\Delta_{-mn}=(\sum_{m \in  \mathcal{M}}\sum_{n \in \mathcal{N}_{m}} C_{mn} y_{mn}-C_{mn} y_{mn})- \sum_{i \in \mathcal{I} } Q_{i} x_{i}=\Delta-C_{mn}y_{mn},
\end{equation}

\HX{It is widely accepted that traveler's waiting time and TSP's idle time are
highly dependent on the supply of participating TSPs and demand of participating travelers, e.g., \cite{bai2019coordinating}. We assume that the MaaS regulator estimates travelers' average waiting time and TSPs' average idle time based on the supply-demand gap $\Delta$. These estimated average waiting/idle time are  reported to travelers/TSPs. We assume that, in the long run, Traveler $i$ (resp. TSP $mn$) observes $\Delta$ based on the estimated waiting (resp. idle) time information provided by the MaaS regulator, and perceives average waiting (resp. idle) time as functions of $\Delta_{-i}$ (resp. $\Delta_{-mn}$), and then use this perceived time information to make their own decisions.}

\HX{As Traveler $i$'s perceived supply-demand gap ($\Delta_{-i}$) increases, her perceived average waiting time decreases. Similarly, as TSP $mn$'s perceived supply-demand gap ($\Delta_{-mn}$) increases, her perceived average idle time increases. Let $\Psi_{i}(x_{i},\Delta)$ denote the perceived waiting time function of traveler $i \in \mathcal{I}$ and $\Phi_{mn}(y_{mn},\Delta)$ denote the perceived idle time of TSP $mn$, $m\in \mathcal{M}, n \in \mathcal{N}_m$. We consider two types of perceived waiting/idle time functions for travelers/TSPs: linear and quadratic. }

Linear waiting/idle functions of Traveler $i$/TSP $mn$ are given in \eq\eqref{L1}/\eq\eqref{L2}.
\begin{subequations}
\begin{align}
&\text{Traveler $i$'s perceived waiting time}:\quad\Psi_{i}^{l}(x_{i},\Delta) =-\frac{\zeta}{\overline{C}}(\Delta+Q_{i}x_{i})+\vartheta,\label{L1}\\
&\text{TSP $mn$'s perceived idle time}:\quad\Phi_{mn}^{l}(y_{mn},\Delta)= \frac{\kappa}{\overline{C}}(\Delta-C_{mn}y_{mn})+\xi ,\label{L2}
\end{align}
\label{eq:waiting}
\end{subequations}
Quadratic waiting/idle functions of Traveler $i$/TSP $mn$ are given in \eq\eqref{Q1}/\eq\eqref{Q2}.
\begin{subequations}
\begin{align}
&\text{Traveler $i$'s perceived waiting time}:\Psi_{i}^{q}(x_{i},\Delta)=\frac{\zeta}{\overline{C}^{2}}(\Delta+Q_{i}x_{i})^{2}-\frac{2\zeta}{\overline{C}}(\Delta+Q_{i}x_{i})+\vartheta,\label{Q1}\\
&\text{TSP $mn$'s perceived idle time}:\quad\Phi_{mn}^{q}(y_{mn},\Delta)= \frac{\kappa}{\overline{C}^{2}}(\Delta-C_{mn}y_{mn})^{2}+\xi.\label{Q2}
\end{align}
\label{eq:idle}
\end{subequations}
\indent In \eq\eqref{eq:waiting}-\eqref{eq:idle}, $\zeta$, $\kappa$, $\vartheta$ and $\xi$ are positive parameters which ensure that perceived waiting/idle time functions are nonnegative. Observe that if $\Delta+Q_{i}x_{i}\in [0,\overline{C}]$, then $\Psi_{i}^{l}(x_{i},\Delta)$ and $\Psi_{i}^{q}(x_{i},\Delta)$ are strictly monotone decreasing with regards to $\Delta$ and satisfy $\Psi_{i}^{l}(x_{i},\Delta), \Psi_{i}^{q}(x_{i},\Delta)\in [\vartheta-\zeta,\vartheta]$. Analogously, if $\Delta-C_{mn}y_{mn}\in [0,\overline{C}]$, then $\Phi_{mn}^{l}(y_{mn},\Delta)$ and $\Phi_{mn}^{q}(y_{mn},\Delta)$ are strictly monotone increasing with regards to $\Delta$ and satisfy $\Phi_{mn}^{l}(y_{mn},\Delta), \Phi_{mn}^{q}(y_{mn},\Delta)\in [\xi,\kappa+\xi]$.

\begin{figure}[!t]
		\centering
	\caption{Overview of the SLMFG model in the two-sided MaaS markets}
	\vspace{0.1in}
 \includegraphics[width=1\textwidth]{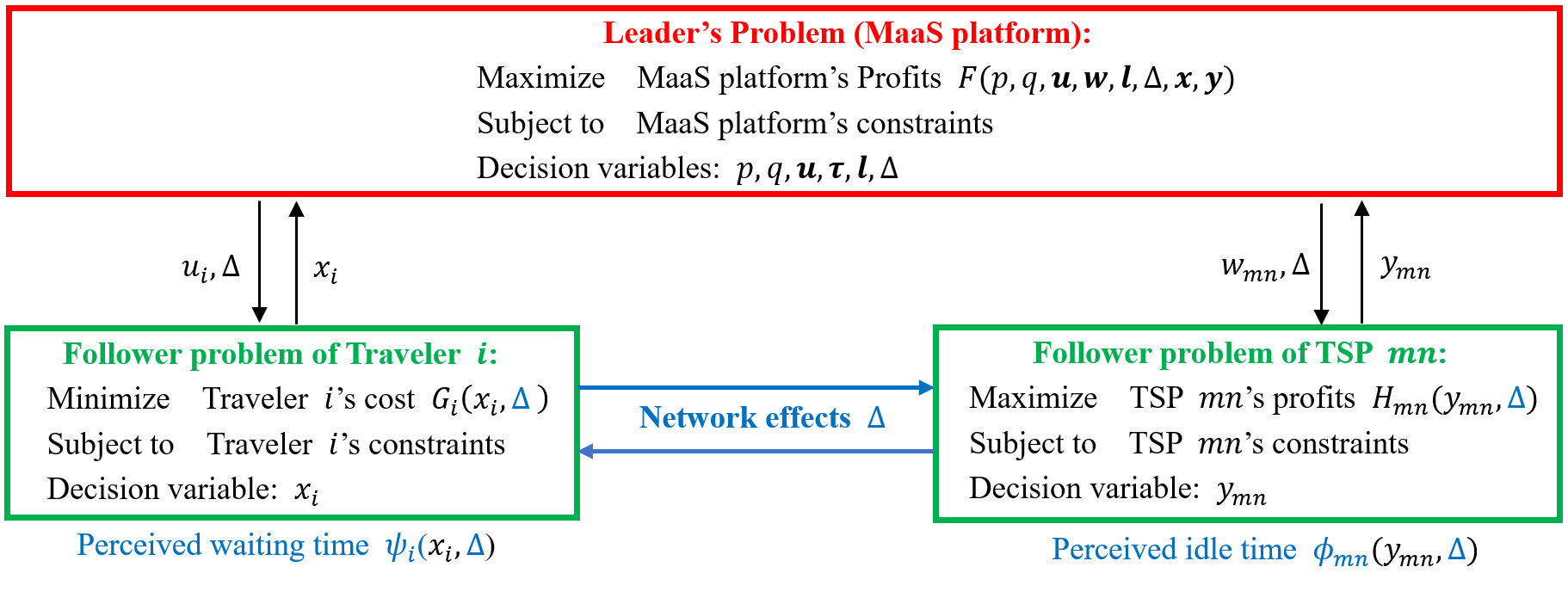}		\label{net}
\end{figure}

The proposed SLMFG model for two-sided MaaS markets is illustrated in \Fig\ref{net}. The MaaS regulator aims to maximize profits by adjusting the \HX{unit threshold price} for travelers ($p$) and the \HX{unit threshold price} for TSPs ($q$), as well as making \textit{MaaS bundles} ($\bm{l}$), while anticipating the choices of followers on participation levels of the MaaS platform. On the one hand, each traveler aims to minimize her travel costs by choosing the proportion of their mobility demand ($x_i$) fulfilled via the MaaS platform. On the other hand, TSPs aim to maximize its profits by choosing the proportion of their mobility resources ($y_{mn}$) supplied to the MaaS platform. \HX{The ``supply'' of participating TSPs and the ``demand'' of travelers' requests
are endogenously dependent on the \HX{unit threshold price} set by the leader and other followers' decisions.} The network effects information is captured by travelers' perceived waiting time and TSPs' perceived idle time, both of which are functions of the supply-demand gap $\Delta$. As $\Delta$ increases, travelers' perceived waiting time will decrease and TSPs' perceived idle time will increase. Thus, in the proposed two-sided MaaS market, the participation levels of travelers (resp. TSPs) is affected by the participation levels of TSPs (resp. travelers).

\section{SLMFG formulations for two-sided MaaS markets}
\label{SLMFGs}
In this section, we introduce SLMFG formulations in \s \ref{formulation}, the SLMFG base formulation in \s \ref{basicmodel}, and dedicate SLMFG formulations with and without network effects in \s \ref{Sec:Network}.

\subsection{SLMFG formulations \label{formulation}} We first introduce two follower problems before presenting the leader problem.
\subsubsection{Follower problems.\label{sectionfollower}} We consider two groups of follower problems, one per side of the two-sided MaaS market, i.e., travelers and TSPs. 

\indent\textbf{Travelers' problem}: The objective of travelers is to minimize their total travel costs. For each traveler $i \in \mathcal{I}$, let $r_{i}$ denote Traveler $i$'s unit reserve price and $\alpha_{i}$ denotes Traveler $i$'s unit waiting time cost. The travel cost of Traveler $i$ include travel cost for using the MaaS platform, travel cost for her reserve travel options, and the perceived waiting time cost $\alpha_{i}\Psi_{i}(x_{i},\Delta)$. The objective function of traveler's problem is denoted by $G_{i}(x_{i},\Delta)$ and given in \eq\eqref{BL1a}.
\begin{subequations}
\allowdisplaybreaks
\begin{align}
&\min_{x_{i}} \quad \underbrace{G_{i}(x_{i},\Delta)}_{\text{Total travel cost}}= \underbrace{b_{i}Q_{i}x_{i}}_{\text{MaaS travel cost}}+\underbrace{r_{i}Q_{i}(1-x_{i}) }_{\text{reserve travel cost}}+\underbrace{\alpha_{i}\Psi_{i}(x_{i},\Delta)}_{\text{waiting time cost}}.\tag{T.1} \label{BL1a}
\end{align}

Recall that $B_i$ is Traveler $i$'s travel expenditure budget for using MaaS, hence we require:
\begin{align}
&b_{i}Q_{i}x_{i}\leq B_{i}.\tag{T.2}\label{BL1b} 
\end{align}

To capture accept/reject decisions of the MaaS regulator in the strategy of Traveler $i$, leader variable $u_i$ is introduced as an upper bound for $x_i$, i.e.,
 \begin{align}
&x_{i} \leq u_{i},\tag{T.3}\label{BL1c}\\
&x_{i}\geq 0.\tag{T.4}\label{BL1d}     
 \end{align}
\end{subequations}

In sum, the follower problem of Traveler $i \in \mathcal{I}$ consists of Eqs.\eqref{BL1a}-\eqref{BL1d}.\\

\indent\textbf{TSPs' problem:} The objective of TSPs is to maximize their profits. For any $m\in\mathcal{M}$ and $n\in\mathcal{N}_{m}$, let $\gamma_{mn}$ denote TSP $mn$'s unit operating cost, $\rho_{mn}$ denote TSP $mn$'s unit reserve price, and $\eta_{mn}$ denote TSP $mn$'s unit idle time cost. Recall that $\beta_{mn}$ is the unit sell-bidding price of TSP $mn$. The objective function of TSP $mn$ is denoted $H_{mn}(y_{mn},\Delta)$ and given in \eq\eqref{BL2a}:
\begin{subequations}
\begin{align}
&\max_{y_{mn}} \hspace*{0.35cm} \underbrace{H_{mn}(y_{mn},\Delta)}_{\text{Total profits}}=\underbrace{\left (\beta_{mn}-\gamma_{mn}\right)C_{mn}y_{mn}}_{\text{revenue from MaaS}}+\underbrace{(\rho_{mn}-\gamma_{mn})C_{mn}(1-y_{mn})}_{\text{revenue from reserve options}}-\underbrace{\eta_{mn}\Phi_{mn}(y_{mn},\Delta)}_{\text{idle time cost }}.\tag{P.1}\label{BL2a} 
\end{align}

Recall that $\bar{B}_{mn}$ is TSP $mn$'s operating cost budget for supplying the MaaS platform, we require: 
\begin{align}
 &C_{mn}\gamma_{mn}y_{mn}\leq \bar{B}_{mn}.\tag{P.2}\label{BL2b}
\end{align}

To capture accept/reject decisions of the MaaS platform in the strategy of TSP $mn$, leader variable $w_{mn}$ is introduced as an upper bound for $y_{mn}$, i.e.:
\begin{align}
&y_{mn}\leq w_{mn},\tag{P.3} \label{BL2c}\\
&y_{mn}\geq 0.\tag{P.4}
\label{BL2d}
\end{align}
\end{subequations}

In sum, the follower problem of TSP $mn$, $m \in \mathcal{M}, n \in \mathcal{N}_m$ consists of Eqs.\eqref{BL2a}-\eqref{BL2d}.

\subsubsection{Leader's problem. \label{sectionleader}} The objective of the MaaS regulator is to maximize its profits, which is the difference between the total revenue obtained from travelers and the total payments made to TSPs. The  objective function of the leader is denoted by $F(\cdot)$ and given in \eq\eqref{B1a}.
\begin{subequations}
\begin{equation}
\max_{p,q,\bm{l},\bm{u},\bm{w},\Delta}F(p,q,\bm{l},\bm{u},\bm{w},\Delta,\bm{x},\bm{y})= \left\{
\underset{\bm{x},\bm{y}}{\text{``max''}} \left\{\sum_{i \in \mathcal{I}} b_{i}Q_{i}x_{i}-\sum_{m \in \mathcal{M}}\sum_{n \in \mathcal{N}_{m}} \beta_{mn}C_{mn} y_{mn}\right\}
\right\},\tag{M.1}\label{B1a}\\
\end{equation}

The inner ``max'' in \eq\eqref{B1a} is equivalent to an optimistic SLMFG formulation in case the response of followers to a leader decision is not unique. This optimistic approach to bilevel optimization attempts to capture the best equilibrium reaction of the followers with regard to the leader's objective; see \cite{AusselSvenssonJNCA,aussel2019towards,AusselSvenssonChapter} for other options.

The constraints of the leader problem link follower variables ($\bm{x}$ and $\bm{y}$) with the \textit{MaaS bundle} allocation variables ($\bm{l}$), \HX{unit threshold prices} for travelers and TSPs ($p$ and $q$), as well as accept/reject decisions for travelers ($\bm{u}$) and TSPs ($\bm{w}$). Recall that variable $l_{i}^{m}$ represents the service time allocated to mode $m \in \mathcal{M}$ in Traveler $i$'s \textit{MaaS bundle}. Let $v_m$ be the commercial speed of mode $m \in \mathcal{M}$. 
\HX{Let $l_{i}^{m}$ denote the service time of travel mode $m$ in the \textit{MaaS bundle} allocated to Traveler $i$.} Traveler $i$'s decision $x_{i}$ is linked to variable $l_{i}^m$ via constraint \eqref{B1b} which requires that the proportion of travel distance $D_i$ fulfilled via MaaS is distributed across multi-travel modes.
\begin{align}
&\sum_{m \in \mathcal{M}} v_{m} l_{i}^{m} = D_{i} x_{i},&& \forall i \in \mathcal{I}.\tag{M.2}\label{B1b}      
\end{align}

Recall that $R_i$ denote Traveler $i$'s travel delay budget, combining this upper bound on travel delay along with the requested total service time $T_{i}$ yields:
\begin{align}
& 0\leq\sum_{m \in \mathcal{M}}l_{i}^{m}-T_{i}x_{i}\leq {R}_{i},&& \forall i \in \mathcal{I}.\tag{M.3}\label{B1c}
\end{align}

Further, we assume that travelers perceive inconvenience cost of using different travel modes differently, which is a common assumption in the literature, e.g., \cite{Bian2019}. Let $\sigma_m$  denote the inconvenience cost of mode $m \in \mathcal{M}$ per unit of time, and let $\Gamma_i$ denote the inconvenience tolerance of user $i$. The total travel inconvenience of Traveler $i$ must not exceed $\Gamma_i$:
\begin{align}
&\sum_{m \in  \mathcal{M}}\sigma_{m} l_{i}^{m} \leq  \Gamma_{i},&& \forall i \in \mathcal{I}. \tag{M.4}\label{B1d}
\end{align}

Supply-side constraints require that, for each travel mode $m \in \mathcal{M}$, the quantity of mobility resources allocated to travelers does not exceed the quantity of mobility resources supplied by TSPs. This links variables $y_{mn}$ to variables $l_{i}^m$ as follows:
\begin{align}
&\sum_{i \in \mathcal{I}} v_{m}^{2}l_{i}^{m}\leq \sum_{n \in \mathcal{N}_{m}} C_{mn}y_{mn},&& \forall m \in  \mathcal{M}.\tag{M.5}\label{B1e}
\end{align}

Let $p_{min}$ and $p_{max}$ be lower and upper bounds on $p$. As introduced in NYOP mechanism, constraints \eqref{B1f}-\eqref{B1g} indicate that if Traveler $i$'s unit bidding price $b_{i}$ is not smaller than the \HX{unit threshold price} $p$, then Traveler $i$ may use the MaaS platform (accept); otherwise, Traveler $i$ is rejected, i.e. if $b_{i}\geq p$, then $u_{i}=1$; otherwise, $u_{i}=0$. 
\begin{align}
&b_{i}-p\geq(1-u_{i})(b_{i}-p_{max}),&& \forall i \in \mathcal{I},\tag{M.6}\label{B1f}\\
&b_{i}-p\leq u_{i}(b_{i}-p_{min}),&& \forall i \in \mathcal{I}.\tag{M.7}\label{B1g}
\end{align}

Let $q_{min}$ and $q_{max}$ be lower and upper bounds on $q$, respectively. Analogously, constraints \eqref{B1h}-\eqref{B1i} indicate that if the unit bidding price of TSP $mn$ is not greater than the \HX{unit threshold price} $q$, then TSP $mn$ may join the MaaS platform (accept); otherwise, TSP $mn$ is rejected, i.e. if $q\geq \beta_{mn}$, then $w_{mn}=1$; otherwise, $w_{mn}=0$. 
\begin{align}
&q-\beta_{mn}\geq(1-w_{mn})(q_{min}-\beta_{mn}),&& \forall m \in \mathcal{M},n \in \mathcal{N}_{m},\tag{M.8}\label{B1h}\\
&q-\beta_{mn}\leq w_{mn}(q_{max}-\beta_{mn}),&& \forall m \in \mathcal{M},n \in \mathcal{N}_{m}.\tag{M.9}\label{B1i}
\end{align}

Constraint \eqref{B1j} links the supply-demand gap variable $\Delta$ to followers' decisions $\bm{x}$ and $\bm{y}$. We assume that the leader aims to maintain the supply-demand gap in the range $[\underline{C},\overline{C}]$, \HX{where $\underline{C}$ is the reserved capacity of the MaaS regulator and $\overline{C}$ is the capacity of  mobility resources provided by all TSPs, namely, $\overline{C}=\sum_{m\in\mathcal{M}}\sum_{n\in\mathcal{N}_{m}}C_{mn}$.} \HX{From a practical standpoint, we assume that the MaaS regulator observes followers' decisions $\bm{x}$ and $\bm{y}$, and adjusts $\Delta$ accordingly. This is reflected in terms of the estimated waiting/idle time information provided to the followers.} 
\begin{equation}
\Delta= \sum_{m \in \mathcal{M}}\sum_{n \in \mathcal{N}_{m}}C_{mn} y_{mn} - \sum_{i \in \mathcal{I}} Q_{i} x_{i}.\tag{M.10}\label{B1j}  
\end{equation}

This is summarized, along with other bound constraints on leader variables in \eq\eqref{B1k}-\eqref{B1p}.
\begin{align}
&u_{i}\in \left\{0,1\right\},&& \forall i \in \mathcal{I},\tag{M.11}\label{B1k} \\
&w_{mn}\in \left\{0,1\right\},&& \forall m \in \mathcal{M},n \in \mathcal{N}_{m},\tag{M.12}\label{B1l}\\
&l_{i}^{m}\geq 0, && \forall i \in \mathcal{I},m \in \mathcal{M},\tag{M.13}\label{B1m}\\
&p_{min}\leq p\leq p_{max},\tag{M.14}\label{B1n}\\
&q_{min}\leq q\leq q_{max},\tag{M.15}\label{B1o}\\
&\underline{C}\leq\Delta\leq \overline{C}.\tag{M.16}\label{B1p}
\end{align}
\end{subequations}

\subsection {\DA{SLMFG base model}} \label{basicmodel}

Our modelling of the two-sided MaaS market is based on the SLMFG approach. We first present the base model of the proposed SLMFG. Let us define \DR{the feasible set $\mathcal{S}$ of the leader problem:}
\begin{equation}
\begin{array}{c}
\DA{
\mathcal{S}=\left\{(p,q,\bm{l},\bm{u},\bm{w},\Delta,\bm{x},\bm{y}): \hspace*{0.1cm}\text{Eqs.}\hspace*{0.1cm} \eqref{B1b}-\eqref{B1p}\right\}.
}
\end{array}
\end{equation}

The feasible sets of travelers and TSPs problems are parameterized by accept/reject binary variables $\bm{u}$ and $\bm{w}$ and denoted $\mathcal{S}_{i}(u_{i})$ and $\mathcal{S}_{mn}(w_{mn})$, respectively:
\begin{equation}
\mathcal{S}_{i}(u_{i})=\left\{x_{i}\in \mathbb{R}:\hspace*{0.1cm}\text{Eqs.}\hspace*{0.1cm}\eqref{BL1b}-\eqref{BL1d}\right\}, \quad \forall i \in\mathcal{I}.
\end{equation}
\begin{equation}
\mathcal{S}_{mn}(w_{mn})=\left\{y_{mn}\in \mathbb{R}:\hspace*{0.1cm}\text{Eqs.}\hspace*{0.1cm}\eqref{BL2b}-\eqref{BL2d}\right\}, \quad \forall m\in \mathcal{M}, n\in \mathcal{N}_{m}.
\end{equation}

The SLMFG base model for two-sided MaaS market is summarized in Model \ref{SLMFG}.
\begin{model}[SLMFG]
\label{SLMFG}
\begin{subequations}
\allowdisplaybreaks
\begin{align}
&\max_{p,q,\bm{l},\bm{u},\bm{w},\Delta,\bm{x},\bm{y}}\quad\text{\emph{MaaS regulator's profits}} \quad\eqref{B1a},\nonumber\\
&\text{\emph{subject to:}}
\quad  && \nonumber \\
&\text{\emph{MaaS regulator's MaaS bundles constraints \quad  \eqref{B1b}-\eqref{B1e}}}\nonumber\\
&\text{\emph{MaaS regulator's acceptance/rejection constraints\quad  \eqref{B1f}-\eqref{B1i}}}\nonumber\\
&\text{\emph{Supply-demand gap\quad \eqref{B1j}}}\nonumber\\
&\text{\emph{MaaS regulator's variable bounds\quad  \eqref{B1k}-\eqref{B1p}}}\nonumber\\
&\qquad \emph{Traveler $i$'s follower problem} \quad x_{i}\in \arg \min _{\hat {x}_{i}} G_{i}(\hat{x}_{i},\Delta),&& \forall i \in \mathcal{I},\\
&\qquad \qquad\text{\emph{subject to:}} \quad \hat{x}_{i}\in \mathcal{S}_{i}(u_{i}),&& \forall i \in \mathcal{I},\\
&\qquad\emph{TSP $mn$'s follower problem:} \quad y_{mn}\in \arg \max_{\hat{y}_{mn}} H_{mn}(\hat{y}_{mn},\Delta),&& \forall m \in \mathcal{M},n \in \mathcal{N}_{m},\\
&\qquad\qquad \text{\emph{subject to:}} \quad \hat{y}_{mn}\in \mathcal{S}_{mn}(w_{mn}),&& \forall m \in \mathcal{M},n \in \mathcal{N}_{m}.
\end{align}
\end{subequations}
\end{model}

\DA{This base model can be rewritten in the following more compact form:
\begin{equation}
\max\left\{F(p,q,\bm{l},\bm{u},\bm{w},\Delta,\bm{x},\bm{y}):(p,q,\bm{l},\bm{u},\bm{w},\Delta,\bm{x},\bm{y}) \in \text{IR} \right\}   
\end{equation}
where IR denotes the {\em Inducible Region} and represents the feasible set of the SLMFG:
\begin{small}
\begin{equation}
\mathrm{IR} \triangleq \left\{\left(p, q, \bm{l}, \bm{u}, \bm{w}, \Delta, \bm{x}, \bm{y}\right) \in \mathcal{S}: x_{i} \in \mathcal{P}_{i}\left(u_{i}, \Delta\right), \forall i \in \mathcal{I}, y_{mn} \in \mathcal{P}_{m n}\left(w_{mn}, \Delta\right), \forall m \in \mathcal{M}, n \in \mathcal{N}_{m}\right\}
\end{equation}
\end{small}
with $\mathcal{P}_{i}(u_{i},\Delta)$ and $\mathcal{P}_{mn}(w_{mn},\Delta)$ being respectively the {\em sets of rational reactions} of the travelers and of the TSPs when the MaaS regulator decision is $(\bm{u},\bm{w},\Delta)$:
\begin{equation}
\mathcal{P}_{i}(u_{i},\Delta) \triangleq \left\{x_{i} \in \arg \min \left\{G_{i}({\hat{x}_{i},\Delta)}:\hat{x}_{i} \in \mathcal{S}_{i}(u_{i})\right\}\right\},\quad \forall i \in\mathcal{I}.
\end{equation}
\begin{equation}
\mathcal{P}_{mn}(w_{mn},\Delta) \triangleq \left\{y_{mn} \in \arg \max \left\{H_{mn}(\hat{y}_{mn},\Delta):\hat{y}_{mn} \in \mathcal{S}_{mn}(w_{mn}) \right\}\right\},\quad \forall m\in \mathcal{M}, n\in \mathcal{N}_{m}.
\end{equation}
}
\indent One can observe that the feasible set $\mathcal{S}_{i}(u_{i})$ and the objective function $G_{i}(\hat{x}_{i},\Delta)$ of each traveler only depends on the leader variables, i.e., $u_i$ and $\Delta$, and do not depends on the other followers' variables. Since this observation also holds for the TSP problems, the follower problems are decoupled and thus the lower level Nash game turns out to be a ``concatenation'' of independent optimization problems. Thus the travelers' problems can be grouped in one where the objective function is the sum of the objective functions of all travelers and the feasible set is the product of the feasible sets of all travelers. \DR{Let $\mathcal{MN}=\left\{(m,n):\forall m \in \mathcal{M}, \forall n \in \mathcal{N}_{m}\right\}$ be the index set of TSPs,} the same grouping process can be also applied for the TSPs, then we can state Proposition \ref{P1}.

\begin{proposition}[Decoupled SLMFG]
\label{P1}
\DA{The SLMFG for two-sided MaaS markets summarized in \M\ref{SLMFG} is equivalent to the following Single-Leader Decoupled-Follower Game (SLDFG):
\begin{align}
&\max_{p,q,\bm{l},\bm{u},\bm{w},\Delta,\bm{x},\bm{y}}\quad \text{\emph{MaaS regulator's profits}} \quad \eqref{B1a},\nonumber\\
&\text{\emph{subject to:}}\nonumber \\
& \text{\emph{MaaS regulator's constraints}}\quad\eqref{B1b}-\eqref{B1p},\nonumber\\
&\qquad \bm{x}\in\arg\min_{\bm{\hat{x}}}\quad G(\bm{\hat{x}},\Delta)= \sum_{i\in \mathcal{I}}\left[b_{i}Q_{i}\hat{x}_{i}+(1-\hat{x}_{i}) r_{i}Q_{i}+\alpha_{i}\Psi_{i}(\hat{x}_{i},\Delta)\right], \\
&\qquad\text{\emph{subject to }} \hat{\bm{x}}=(\hat{\bm{x}}_1,\dots,\hat{\bm{x}}_{|\mathcal{I}|})\in \prod_{i\in \mathcal{I}}\mathcal{S}_i(u_i)\\
&\qquad \bm{y}\in\arg\max_{\bm{\hat{y}}} \quad H(\bm{\hat{y}},\Delta)=\sum_{m\in\mathcal{M}}\sum_{n\in\mathcal{N}_m}[\left (\beta_{mn}-\rho_{mn}\right)C_{mn}\hat{y}_{mn}+(\rho_{mn}-\gamma_{mn})C_{mn}-\eta_{mn}\Phi_{mn}(\hat{y}_{mn},\Delta)],\\
&\qquad\text{\emph{subject to }} \hat{\bm{y}}=(\hat{\bm{y}}_1,\dots,\hat{\bm{y}}_{|\mathcal{M}\mathcal{N}|})\in \prod_{m\in\mathcal{M}}\prod_{n\in\mathcal{N}_m}\mathcal{S}_{mn}(w_{mn})
\end{align}
in the sense that their solution sets coincide.}
\end{proposition}

Note that the Decoupled SLMFG model will not be adapted to the numerical treatment, in the sequel, we will base our development on the Basic SLMFG (Model \ref{SLMFG}) and variants of it. 


\subsection{SLMFG with and without network effects}\label{Sec:Network}
As announced in Section \ref{PS}, we consider two classes of perceived time functions for the followers: linear and quadratic functions. In the forthcoming subsections \ref{subsec:withoutnet} and \ref{subsec:withnet}, we show that linear perceived time functions lead to model without network effects while with quadratic perceived time functions, network effects are intrinsically associated to the model. 
\DA{
\subsubsection{SLMFG without network effects.\label{subsec:withoutnet}}
When the perceived waiting/idle time functions of the followers are linear, e.g., Eq.\eqref{L1}/\eqref{L2}, the resulting mixed-integer linear bilevel programming (MILBP) formulation is \DR{denoted \SLMFGL\ and summarized in Model \I}.\\}
\noindent\textbf{Model 1.1} (SLMFG-L)
\label{SLMFG-Linear}
\begin{subequations}
\allowdisplaybreaks
\begin{align}
&\max_{p,q,\bm{l},\bm{u},\bm{w},\Delta,\bm{x},\bm{y}} \quad\quad\text{MaaS regulator's profits}\quad\eqref{B1a} \nonumber \\
&\text{subject to:} \nonumber\\
&\text{MaaS regulator's constraints}\quad\eqref{B1b}-\eqref{B1p},\nonumber\\
&\quad x_i\in\arg\min_{\hat{x}_{i}}\quad G_{i}(\hat{x}_{i},\Delta)=b_{i}Q_{i}\hat{x}_{i}+r_{i}Q_{i}(1-\hat{x}_{i}) +\alpha_{i}\Psi_{i}^{l}(\hat{x}_{i},\Delta),\forall i\in\mathcal{I},\label{L1a}\\
&\qquad\qquad\text{subject to:} \quad \text{Traveler $i$'s constraints} \quad\eqref{BL1b}-\eqref{BL1d},\forall i\in\mathcal{I},\\
&\quad y_{mn}\in\arg\max_{\hat{y}_{mn}}~ H_{mn}(\hat{y}_{mn},\Delta)=\left (\beta_{mn}-\rho_{mn}\right)C_{mn}\hat{y}_{mn}+(\rho_{mn}-\gamma_{mn})C_{mn}(1-\hat{y}_{mn})\nonumber\\
&\quad \hspace*{5cm}-\eta_{mn}\Phi_{mn}^{l}(\hat{y}_{m,n},\Delta),\forall\,m\in\mathcal{M},n\in\mathcal{N}_{m},\\
&\qquad\qquad\text{subject to:}  \quad \text{TSP $mn$'s constraints}\quad \eqref{BL2b}-\eqref{BL2d},\forall\,m\in\mathcal{M},n\in\mathcal{N}_{m}.
\end{align}
\end{subequations}
\HX{\indent Then we propose an equivalent formulation for \SLMFGL\ in Proposition \ref{prop:SLMFG_lin}, where the solution set of the follower's problem do not depend on $\Delta$, and thus defined as SLMFG without network effects. The proof of Proposition \ref{prop:SLMFG_lin} is provided in Online Appendix \ref{ECProp2}}.
\begin{proposition}
\label{prop:SLMFG_lin}
If the perceived waiting/idle time functions of the followers are linear, then \SLMFGL\space admits the same set of solutions
with the following SLMFG without network effects:

\begin{subequations}
\allowdisplaybreaks
\begin{align}
&\max_{p,q,\bm{l},\bm{u},\bm{w},\Delta,\bm{x},\bm{y}} \quad\text{MaaS regulator's profits}\quad\eqref{B1a} \nonumber \\
&\text{subject to:}\nonumber\\
&\text{MaaS regulator's constraints}\quad\eqref{B1b}-\eqref{B1p},\nonumber\\
&\quad x_i\in\arg\min_{x_{i}}\quad G_{i}^1(\hat{x}_{i})=b_{i}Q_{i}\hat{x}_{i}+r_{i}Q_{i}(1-\hat{x}_{i}) +\frac{\alpha_{i}\zeta Q_{i}}{\bar{C}}\hat{x}_{i}, \forall i\in\mathcal{I}, \\
&\qquad\quad\text{subject to:} \quad \text{Traveler $i$'s constraints} \quad\eqref{BL1b}-\eqref{BL1d},\forall i\in\mathcal{I},\\
&\quad y_{mn}\in\arg\max_{\hat{y}_{mn}}~ H_{mn}^1(\hat{y}_{mn})=\left (\beta_{mn}-\rho_{mn}\right)C_{mn}\hat{y}_{mn}+(\rho_{mn}-\gamma_{mn})C_{mn}(1-\hat{y}_{mn})\nonumber\\
&\quad\hspace*{4.5cm}-\frac{\eta_{mn}\kappa C_{mn}}{\bar{C}}\hat{y}_{mn},\forall m\in\mathcal{M},n\in\mathcal{N}_{m}, \label{L2a}\\
&\qquad\quad\text{subject to:} \quad \text{TSP $mn$'s constraints}\quad\eqref{BL2b}-\eqref{BL2d},\forall m\in\mathcal{M},n\in\mathcal{N}_{m}.
\end{align}
\end{subequations}
\end{proposition}
\subsubsection{SLMFG with network effects.\label{subsec:withnet}}
When the perceived waiting/idle time functions of the followers are quadratic, i.e., Eq.\eqref{Q1}/\eqref{Q2}, then the resulting mixed-integer quadratic bilevel programming (MIQBP) formulations \DR{is denoted \SLMFGQ\ and summarized in Model \II.\\
}
\noindent\textbf{Model 1.2} (SLMFG-Q)
\label{SLMFG-NE}
\begin{subequations}
\allowdisplaybreaks
\begin{align}
&\max_{p,q,\bm{l},\bm{u},\bm{w},\Delta,\bm{x},\bm{y}}\quad \text{MaaS regulator's profits} \quad \eqref{B1a},\nonumber\\
&\text{subject to:} \nonumber\\
& \text{MaaS regulator's constraints}\quad\eqref{B1b}-\eqref{B1p},\nonumber \\
&\qquad x_i\in\arg\min_{\hat{x}_{i}}~ G_{i}(\hat{x}_{i},\Delta)=b_{i}Q_{i}\hat{x}_{i}+r_{i}Q_{i}(1-\hat{x}_{i}) +\alpha_{i}\Psi_{i}^{q}(\hat{x}_{i},\Delta),\forall i\in\mathcal{I},\label{Q1a}\\
&\qquad\qquad\qquad\text{subject to: Traveler $i$'s constraints} \quad \eqref{BL1b}-\eqref{BL1d},\forall i\in\mathcal{I},\\
&\qquad y_{mn}\in\arg\max_{\hat{y}_{mn}}\quad  H_{mn}(\hat{y}_{mn},\Delta)=\left (\beta_{mn}-\rho_{mn}\right)C_{mn}\hat{y}_{mn}+(\rho_{mn}-\gamma_{mn})C_{mn}(1-\hat{y}_{mn}),\nonumber\\
&\qquad\hspace*{4.5cm}-\eta_{mn}\Phi_{mn}^{q}(\hat{y}_{m,n},\Delta),\forall m\in\mathcal{M}, n\in\mathcal{N}_{m}, \label{Q2a}\\
&\qquad\qquad\qquad\text{subject to: TSP $mn$'s constraints}  \quad \eqref{BL2b}-\eqref{BL2d},\forall m\in\mathcal{M}, n\in\mathcal{N}_{m}.
\end{align}
\end{subequations}

Further, taking advantage of the quadratic characteristic of the perceived waiting/idle functions, we can simplify \SLMFGQ\ in Proposition \ref{prop:SLMFG_quadra}, and the proof is provided in Online Appendix \ref{ECProp3}.
\begin{proposition}
\label{prop:SLMFG_quadra}
If  the  perceived  waiting/idle  time  functions  of  the  followers  are  quadratic,  then:
\begin{itemize}
    \item \DR{for any pair of leader variables $(u_i,\Delta)$, each of the traveler problems admits a unique solution, denoted by $x_i(u_i,\Delta)$;
    \item for any pair of leader variables $(w_{mn},\Delta)$, each of the TSP's problem admits a unique solution, denoted by $y_{mn}(w_{mn},\Delta)$;}
    \item \SLMFGQ\ admits  the  same  set  of  solutions with the following SLMFG with network effects:
    \begin{subequations}
    \allowdisplaybreaks
    \begin{align}
    &\max_{p,q,\bm{l},\bm{u},\bm{w},\Delta}\quad \sum_{i \in \mathcal{I}} b_{i}Q_{i}\DR{x_{i}(u_i,\Delta)}-\sum_{m \in \mathcal{M}}\sum_{n \in \mathcal{N}_{m}} \beta_{mn}C_{mn} \DR{y_{mn}(w_{mn},\Delta)},\nonumber\\
    &\text{subject to:} \nonumber\\
    & \text{MaaS regulator's constraints}\quad\eqref{B1b}-\eqref{B1p},\nonumber \\
    &\quad \DR{x_i(u_i,\Delta)}\text{is the unique solution of }\nonumber\\
    &\quad\min_{\hat{x}_{i}}\quad G_{i}(\hat{x}_{i},\Delta)=b_{i}Q_{i}\hat{x}_{i}+r_{i}Q_{i}(1-\hat{x}_{i}) +\alpha_{i}\Psi_{i}^{q}(\hat{x}_{i},\Delta),\forall i\in\mathcal{I},\label{Q1a_S}\\
    &\qquad\qquad\text{subject to: Traveler $i$'s constraints} \quad \eqref{BL1b}-\eqref{BL1d},\forall i\in\mathcal{I},\\
    &\quad \DR{y_{mn}(w_{mn},\Delta)}\text{ is the unique solution of }\nonumber\\
    &\quad\max_{\hat{y}_{mn}}\quad  H_{mn}(\hat{y}_{mn},\Delta)=\left (\beta_{mn}-\rho_{mn}\right)C_{mn}\hat{y}_{mn}+(\rho_{mn}-\gamma_{mn})C_{mn}(1-\hat{y}_{mn})\nonumber\\
    &\hspace*{4cm}-\eta_{mn}\Phi_{mn}^{q}(\hat{y}_{m,n},\Delta),\forall m\in\mathcal{M},\forall n\in\mathcal{N}_{m},\label{Q2a_S} \\
    &\qquad\qquad\text{subject to: TSP $mn$'s constraints}  \quad \eqref{BL2b}-\eqref{BL2d},\forall m\in\mathcal{M},\forall n\in\mathcal{N}_{m}.
    \end{align}
    \end{subequations}
\end{itemize}
\end{proposition}

\indent The fact that the unique solutions \DR{$x_i(u_{i},\Delta)$ and $y_{mn}(w_{mn},\Delta)$ depend on the leader's accept/reject variables $u_{i}$ and $w_{mn}$, respectively,} emphasizes that, in the case of quadratic perceived waiting/idle time functions, network effects are intrinsically associated to \SLMFGQ problem.
\subsection{Illustration of the two-sided MaaS market}
\label{example}

In this section, we give an example to illustrate the behavior of the proposed SLMFG formulations, \SLMFGL\space and \SLMFGQ,  for two-sided MaaS markets.

\begin{figure}[ht!]
\caption{The relationship between $\Delta$ and followers' participation level/objective value in \SLMFGQ.}
		\vspace{-0.1in}
\subfloat[$\Delta-\bm{x}$ ]{\includegraphics[width=0.33\textwidth]{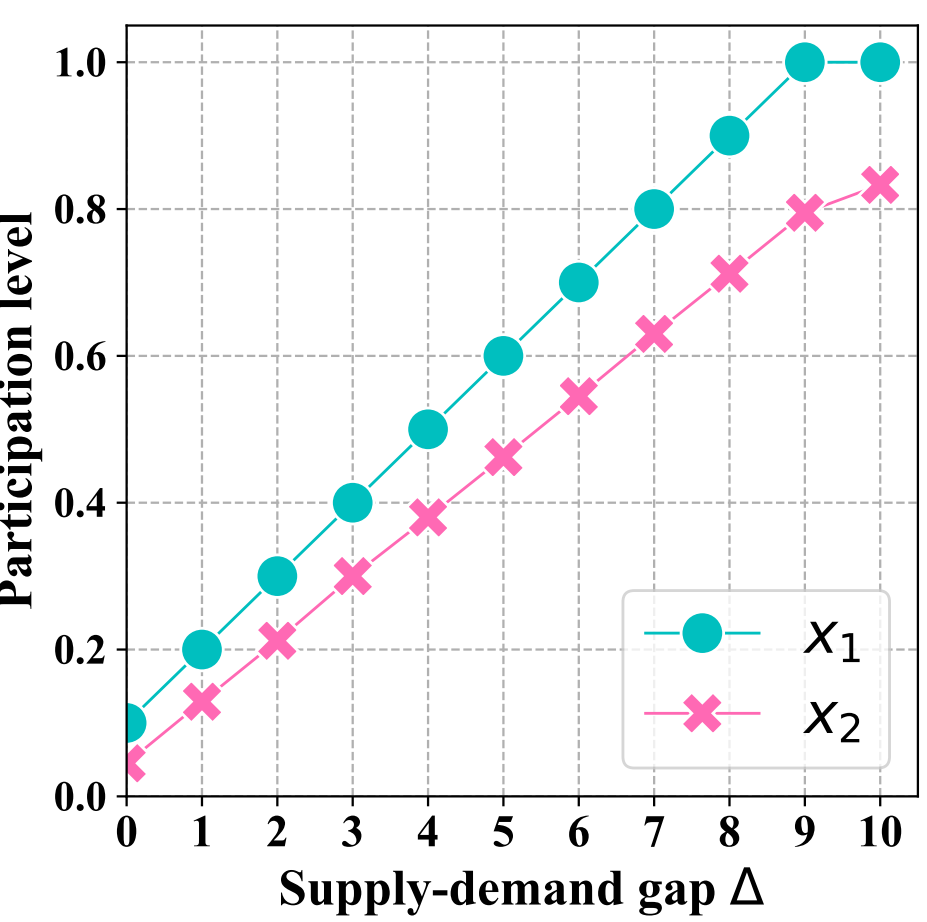}\label{N2a}}
\subfloat[$\Delta-\bm{y}$ ]{\includegraphics[width=0.33\textwidth]{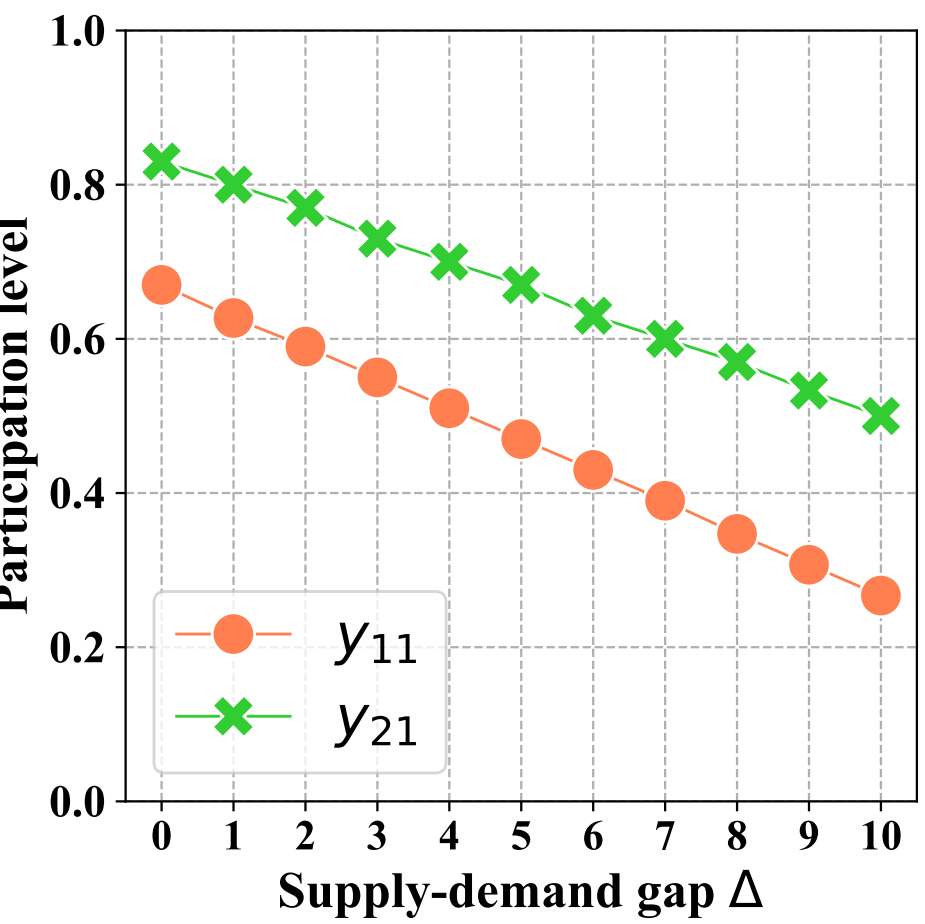}\label{N2b}}
\subfloat[$\Delta-G_{1},\Delta-H_{1}$ ]{\includegraphics[width=0.33\textwidth]{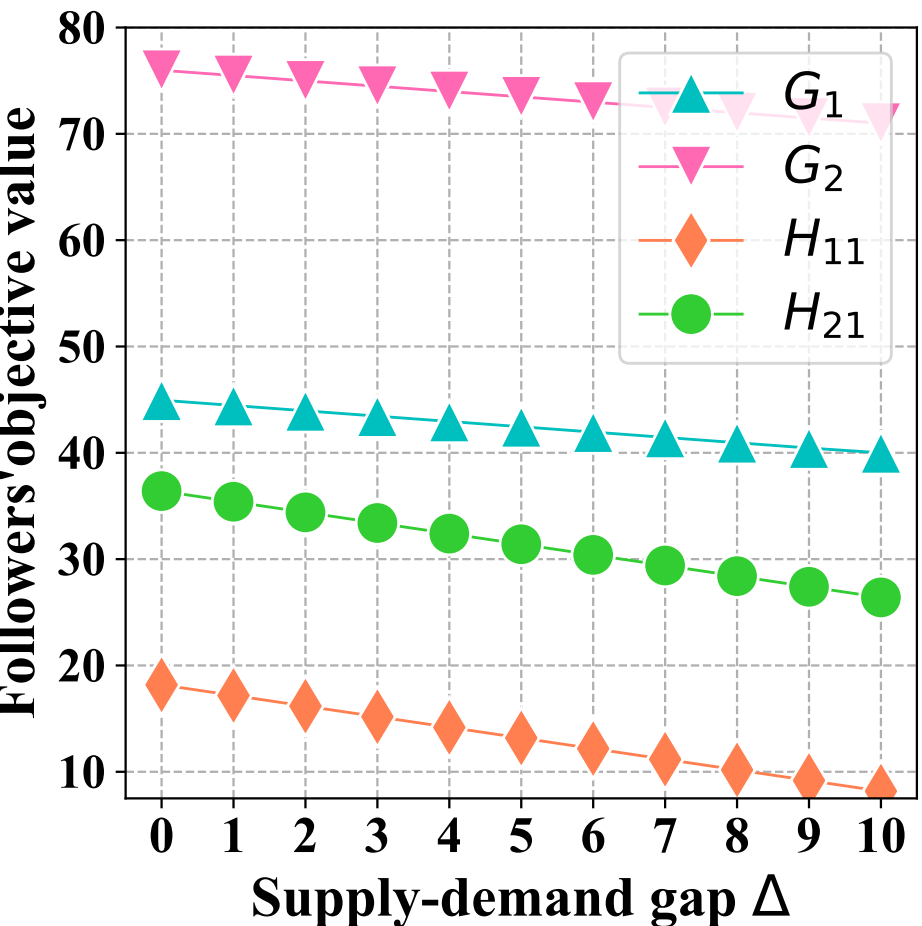}\label{N2c}}
\label{N2}
\end{figure}

\noindent \textbf{Example 1}. Consider a MaaS platform with two travelers $\mathcal{I}=\{1,2\}$, two travel modes $\mathcal{M}=\{1,2\}$ and one TSP per mode $\mathcal{N}_1 = \{1\}$ and $\mathcal{N}_2 = \{1\}$. For each traveler $i \in \mathcal{I}$, unit reserve prices are $r_{1}=\$2.5$ and $r_{2}=\$4.5$, unit waiting time costs are $\alpha_{1}=\$2$ and $\alpha_{2}=\$2.2$, the purchase-bids $\mathcal{B}_{i}=(D_{i},T_{i}, b_{i}, B_{i}, R_{i}, \Gamma_{i})$ are: $\mathcal{B}_{1}=$(40 km, 160 min, \$2,\$20, 50 min, \$200), $\mathcal{B}_{2}=$(60 km, 300 min,\$4, \$40, 30 min, \$300), then the weighted quantity of mobility resources are obtained: $Q_{1}=$ 10 km$^{2}$/min and $Q_{2}=$ 12km$^{2}$/min, Traveler $i$' perceived waiting time functions are $\Psi_{i}^{l}=-0.025(\Delta-Q_{i}x_{i})+10$ and $\Psi_{i}^{q}=0.025(\Delta-Q_{i}x_{i})^{2}-0.2(\Delta-Q_{i}x_{i})+10$. For each TSP $mn$, $m\in \mathcal{M}, n\in\mathcal{N}_m$, unit operating costs are $\gamma_{11}=\$0.1$ and $\gamma_{21}=\$0$, unit idle time costs are $\eta_{11}=\$0.15$ and $\eta_{21}=\$0.1$, unit reserve prices are $\rho_{11}=\$0.5$ and $\rho_{21}=\$1$, the sell-bids $\mathcal{B}_{mn}=(C_{mn},\beta_{mn},\bar{B}_{mn})$ are: $\mathcal{B}_{11}=$(25 km$^{2}$/min, \$1.5,\$3),  $\mathcal{B}_{21}=$(30 km$^{2}$/min, \$2, \$5.5), TSP $mn$'s perceived idle time functions are $\Phi_{mn}^{l}=0.2(\Delta+C_{mn}y_{mn})+1$ and $\Phi_{mn}^{q}=0.2(\Delta+C_{mn}y_{mn})^{2}+1$.
\begin{figure}[b!]
	\caption{The relationship between $\Delta$ and followers' participation level/objective value in \SLMFGL.}
		\vspace{-0.1in}
	\subfloat[$\Delta-\bm{x}$ ]{\includegraphics[width=0.33\textwidth]{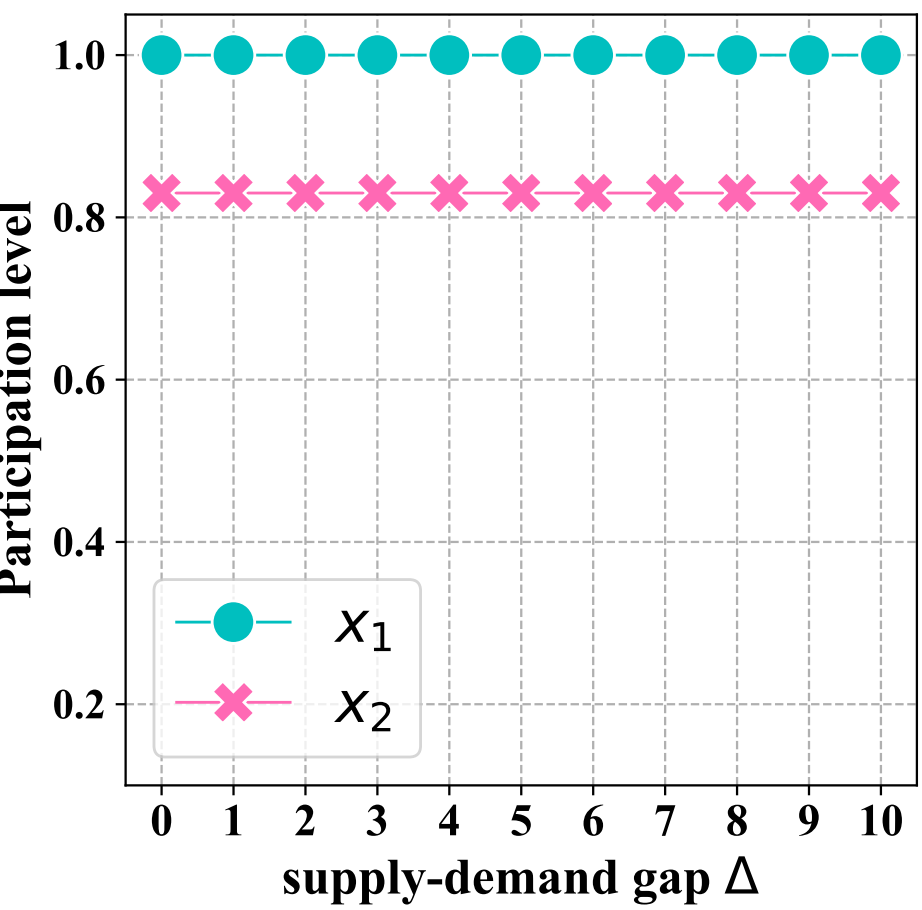}\label{N1a}}
    \subfloat[$\Delta-\bm{y}$]{\includegraphics[width=0.33\textwidth]{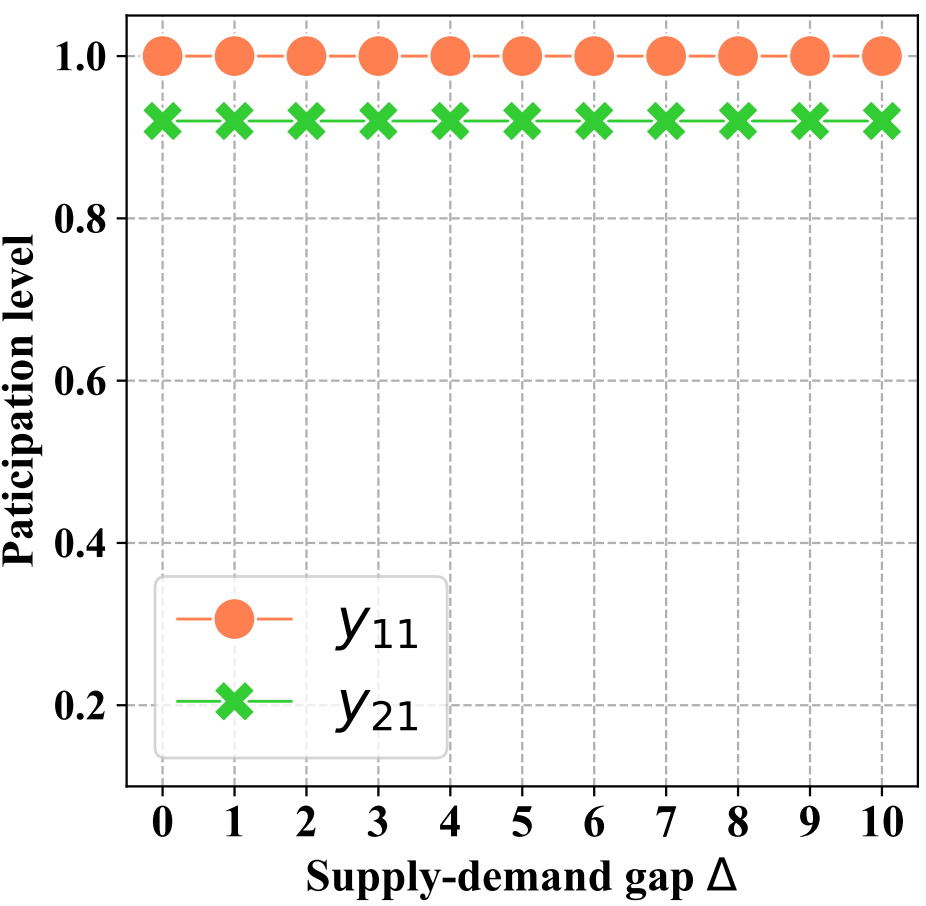}\label{N1b}}
    \subfloat[$\Delta-\bm{G}$ and $\Delta-\bm{H}$]{\includegraphics[width=0.33\textwidth]{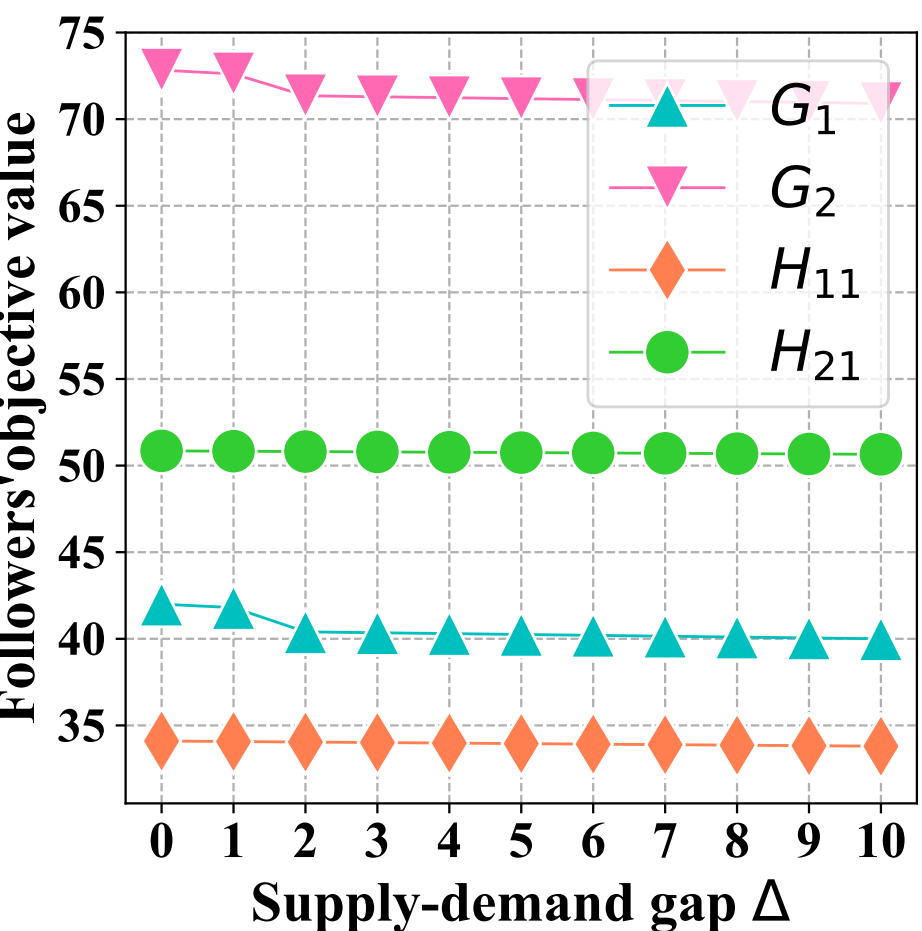}\label{N1c}}
\label{N1}
\end{figure}
\HX{For the convenience of analyzing the two-sided MaaS markets, we assume that followers at both sides of the MaaS market are accepted by the leader ($u_{1}=1$, $u_{2}=1$, $w_{11}=1$, $w_{21}=1$)}. Given the \HX{unit threshold price} $p=\$1$ and $q=\$2$ for travelers and TSPs, when the value of $\Delta$ is changed from 0 to 10, we solve the corresponding follower problem \eqref{BL1a}-\eqref{BL1d} and \eqref{BL2a}-\eqref{BL2d} to obtain the values of $x_{1}$, $x_{2}$, $y_{11}$ and $y_{21}$. We analyze the relationship between the supply-demand gap ($\Delta$) and followers' participation levels/objective values in \Fig\ref{N1}/\Fig\ref{N2}. As $\Delta$ increases in \SLMFGL, the participation level of each traveler ($\bm{x}$) and TSP ($\bm{y}$) will not change (\Fig\ref{N1a} and \ref{N1b}), and each traveler's objective value ($\bm{G}$) and TSP's objective value ($\bm{H}$) will slightly decrease (\Fig\ref{N1c}). As $\Delta$ increases in \SLMFGQ, each follower's participation level ($\bm{x}$ and $\bm{y}$) will increase (\Fig\ref{N2a} and \ref{N2b}), and each follower's objective value ($\bm{G}$ and $\bm{H}$) will considerably decrease (\Fig\ref{N2c}). 

\begin{figure}[ht!]
\centering
\caption{The relationship between the \HX{unit threshold price} $p,q$ and followers' participation level $\bm{x},\bm{y}$ in \SLMFGL.} 
\vspace{-0.1in}
\subfloat[$p-\bm{x}$ ]{\includegraphics[width=0.4\textwidth]{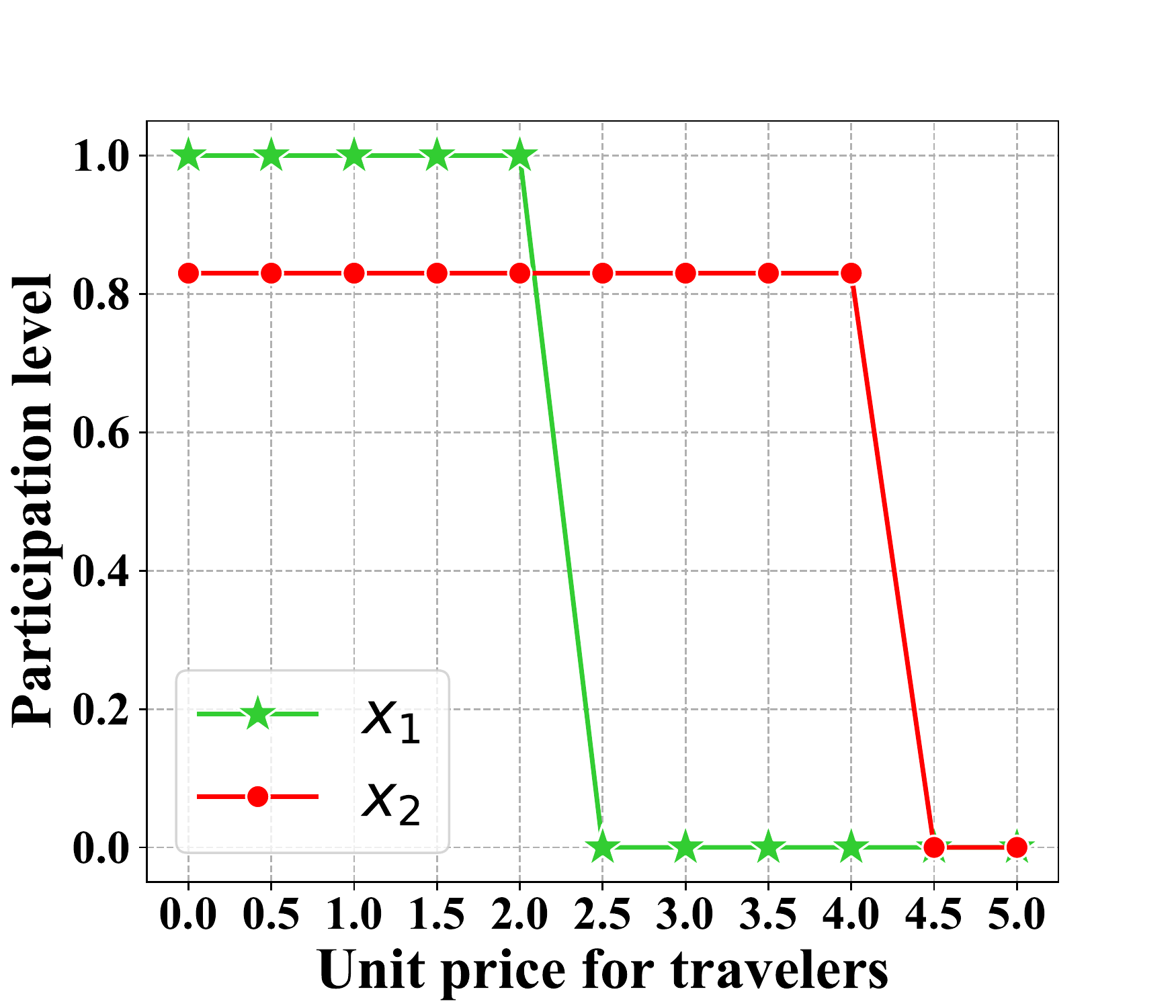}\label{N3a}}
\subfloat[$q-\bm{y}$ ]{\includegraphics[width=0.4\textwidth]{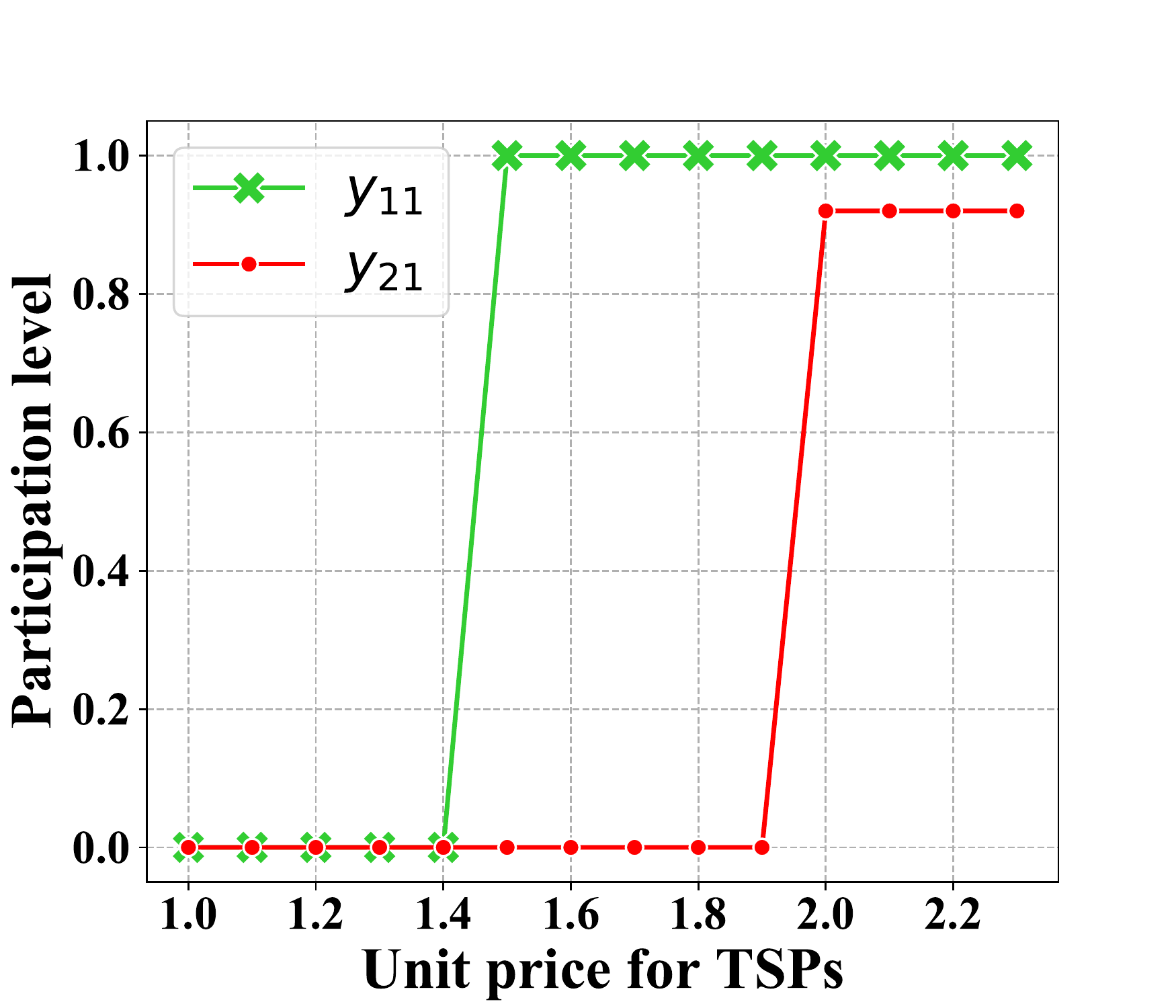}\label{N3b}}
\label{N3}
\end{figure}
Then we analyze the impacts of \HX{unit threshold price} for travelers and TSPs ($p$ and $q$) determined by the MaaS regulator on their participation levels ($\bm{x}$ and $\bm{y}$). As is shown in \Fig\ref{N3a}, if $p\leq b_{1}\leq b_{2}$, then $x_{1}=1$ and $x_{2}=0.83$; if $b_{1}\leq p\leq b_{2}$, then $x_{1}=0$, $x_{2}=0.83$; and if $b_{1}\leq b_{2}\leq p$, then $x_{1}=0$, $x_{2}=0$. Analogously, as the value of $q$ increases from 1 to 2.2, the value of $w_{11}$ and $u_{21}$ will change, leading to the change of TSPs' participation levels. In \Fig\ref{N3b}, when $ \beta_{11}\leq \beta_{21}\leq q$, $y_{11}=0$, $y_{21}=0$; when $\beta_{11}\leq q\leq \beta_{21}$, $y_{11}=1$, $y_{21}=0$; when $\beta_{11}\leq\beta_{21}\leq q$, $y_{11}=1$, $y_{21}=0.85$. 

The results of Example 1 are displayed in \Fig\ref{N1}-\ref{N3}. We summarize the obtained insights facilitating the MaaS regulator to manage the platform in two-sided markets. 

\begin{observation}
Based on the above analysis, the reciprocal interactions among the MaaS regulator, travelers and TSPs can be observed in the SLMFG with network effects (\SLMFGQ). As the supply-demand increases, travelers' participation level will increase and TSPs' participation level will decrease, and both travelers' costs and TSPs' profits will decrease. \HX{Moreover, the participation level of travelers (resp. TSPs) are endogenously dependent on the \HX{unit threshold price} for travelers (resp. TSPs) determined by the MaaS regulator. As the \HX{unit threshold price} for travelers increases, travelers' participation levels decrease. Conversely, as the \HX{unit threshold price} for TSPs increases, TSPs' participation levels increase.} 
\end{observation}

\section{Solution methods}
\label{Solution methods}
In Section \ref{Single}, we present two approaches to reformulate \SLMFGL\space and \SLMFGQ\space based on the optimality conditions of the follower problems: i) KKT conditions and ii) SD conditions. The classic B\&B algorithm proposed by \cite{bard1990branch}, \Bard, is based on the MPEC reformulation serving as the benchmark in our numerical experiments. In Section \ref{algorithm}, we propose a SD-based B\&B algorithm  to solve both MILBP and MIQBP problems. In Section \ref{llustrative example}, we present a numerical example to illustrate the proposed B\&B algorithm.

\subsection{Single-level reformulation}
\label{Single}
\DA{MPEC reformulations for \SLMFGL\ and \SLMFGQ\  are introduced in Section \ref{MPEC} while the SD-based reformulations for both problems are given in Section \ref{SD}.}   

\subsubsection{ MPEC reformulation.\label{MPEC}} \DA{\indent A classical numerical treatment of bilevel optimization problems is to replace the follower problems by the system of equations/inequations composed of the concatenation of their KKT conditions and then to solve the resulting MPEC. Our aim in this subsection is not only to describe the MPEC formulations associated with the MILBP and MIQBP problems but also to show that, using the specific structure of these problems, it is possible to prove that the bilevel formulations and the corresponding MPEC formulation are equivalent and generate the same solutions (decision variables). This quite rare fact is important to emphasize and will be respectively given in Theorem \ref{thm:equiv_MILBP} for MILBP problem and Theorem \ref{thm:equiv_MIQBP} for MIQBP problem. 
}

Let us first consider MILBP (\SLMFGL). According to Proposition \ref{prop:SLMFG_lin}, the associated KKT conditions of  Traveler $i$'s problem and TSP $mn$'s problem in  \SLMFGL\ are given in Online Appendix \ref{KKT-L}, i.e., Eqs. \eqref{29a}-\eqref{29d} and \eqref{31a}-\eqref{31d}, respectively. 

Then the MPEC formulation of \SLMFGL\ is denoted \MPECL\ and summarized in Model \III.\\
\textbf{Model 2.1 (MPEC-L)}
\begin{subequations}
\label{MPEC_lin}
\allowdisplaybreaks
\begin{align}
&\max_{p,q,\bm{l},\bm{u},\bm{w},\Delta,\bm{x},\bm{y},\bm{\lambda},\bm{\mu}} \quad\text{MaaS regulator's profits}\quad\eqref{B1a} \nonumber \\
&\text{subject to: }\nonumber\\
&\text{The MaaS regulator's constraints:}\quad\eqref{B1b}-\eqref{B1p},\nonumber\\
&\text{Traveler $i$'s KKT conditions:}\quad\eqref{29a}-\eqref{29d},\nonumber\\
&\text{TSP $mn$'s KKT conditions:}\quad\eqref{31a}-\eqref{31d},\nonumber
\end{align}
\end{subequations}
\noindent where $\bm{\lambda}$ and $\bm{\mu}$ stand respectively for $\bm{\lambda}=\left( \lambda^1_{i}, \lambda^2_{i}, \lambda^3_{i}\right)_{i\in \mathcal{I}}$ and $\bm{\mu}=\left( \mu^1_{mn}, \mu^2_{mn}, \mu^3_{mn}\right)_{\stackrel{m\in \mathcal{M}}{n\in \mathcal{N}}}$.

It is well known that usually SLMFG and their associated MPEC reformulations are not equivalent, roughly speaking, they don't share the same solutions, e.g. \cite{AusselSvenssonJNCA,Aussel2019pess}. As analyzed in \cite{AusselSvenssonChapter}, some convexity assumptions coupled with quite strong qualification conditions are needed to guarantee such an equivalence. Nevertheless, \SLMFGL\ is equivalent to the MPEC reformulation \MPECL under some mild hypothesis. We provide constraint qualifications for \SLMFGL\ in Theorem \ref{thm:equiv_MILBP}.
\begin{theorem} \label{thm:equiv_MILBP}
Assume that the perceived waiting/idle time functions of travelers/TSPs are linear.
\begin{itemize}
    \item [ a) ] If $(p,q,\bm{l},\bm{u},\bm{w},\Delta,\bm{x},\bm{y})$ is a global solution of \SLMFGL\space (or equivalently SLMFG without network effects given in Proposition \ref{prop:SLMFG_lin}),  then $\forall (\bm{\lambda},\bm{\mu})\in\Lambda(p,q,\bm{l},\bm{u},\bm{w},\Delta,\bm{x},\bm{y})$, $(p,q,\bm{l},\bm{u},\bm{w},\Delta,\bm{x},\bm{y},\bm{\lambda},\bm{\mu})$ is a global solution of the associated \MPECL.
    \item [ b) ] Assume additionally that 
    \begin{itemize}
    \item [ i) ] for all Traveler $i \in \mathcal{I}$, $b_i \geq p_{min}$;
    \item [ii) ] for all TSP $n \in \mathcal{N}_m, m \in \mathcal{M}$, $\beta_{mn} \leq q_{max}$.
    \end{itemize}
    If $(p,q,\bm{l},\bm{u},\bm{w},\Delta,\bm{x},\bm{y},\bm{\lambda},\bm{\mu})$ is a global solution of \MPECL, then  $(p,q,\bm{l},\bm{u},\bm{w},\Delta,\bm{x},\bm{y})$ is a global solution of \SLMFGL\space (or equivalently SLMFG without network effects).
\end{itemize}
\end{theorem}

In Theorem \ref{thm:equiv_MILBP}, the set $\Lambda(\cdot)$ denotes the set of Lagrange multipliers associated to the concatenated KKT conditions of the followers (travelers and TSPs). The proof of this theorem is fully described in Online Appendix \ref{ECAnnex_C}. Let us say thanks to the decoupled structure of the MILBP, to the linearity of the constraints and objective functions of the follower's problem and  finally to the binary character of the leader's variables $(u_i)_{i\in\mathcal{I}}$ and $(w_{mn})_{\stackrel{m\in\mathcal{M}}{n\in\mathcal{N}}}$, we will be in a position to prove the three main hypothesis of \cite{AusselSvenssonChapter}, i.e., the {\em joint convexity} of the constraint functions of the followers' problem, the {\em Joint Slater's qualification condition} and the {\em Guignard's qualification conditions for boundary opponent strategies} (see Proof in Online Appendix \ref{ECAnnex_C}). 

As for \SLMFGQ, The KKT conditions of Traveler $i$'s problem and TSP $mn$'s problem are given in Online Appendix \ref{KKT-Q}, i.e., Eqs. \eqref{48a}-\eqref{48d} and Eqs. \eqref{49a}-\eqref{49d}, respectively. The MPEC formulation associated to \SLMFGQ\ is denoted as \MPECQ\ and summarized in \M\IV :\\
\textbf{Model 2.2 (MPEC-Q)}
\begin{subequations}
\label{MPEC_quad}
\allowdisplaybreaks
\begin{align}
&\max_{p,q,\bm{l},\bm{u},\bm{w},\Delta,\bm{x},\bm{y},\bm{\lambda},\bm{\mu}} \quad\text{MaaS regulator's profits}\quad\eqref{B1a} \nonumber \\
&\text{subject to: }\nonumber\\
&\text{The MaaS regulator's constraints:}\quad\eqref{B1b}-\eqref{B1p},\nonumber\\
&\text{Traveler $i$'s KKT conditions:}\quad\eqref{48a}-\eqref{48d},\nonumber\\
&\text{TSP $mn$'s KKT conditions:}\quad\eqref{49a}-\eqref{49d},\nonumber
\end{align}
\end{subequations}

\DA{Mimicking the same arguments as in Theorem \ref{thm:equiv_MILBP}, one can easily establish that, in the case of the quadratic perceived waiting/idle time functions and under some mild hypothesis,  \SLMFGQ\ is equivalent to \MPECQ. We provide constraint qualifications for \SLMFGQ\ in Theorem \ref{thm:equiv_MIQBP}.
\begin{theorem} \label{thm:equiv_MIQBP}
Assume that the perceived waiting/idle time functions of followers are quadratic.
\begin{itemize}
    \item [ a) ] If $(p,q,\bm{l},\bm{u},\bm{w},\Delta,\bm{x},\bm{y})$ is a global solution of  \SLMFGQ\space (or equivalently MIQBP described in Proposition \ref{prop:SLMFG_quadra})  then, for any $(\bm{\lambda},\bm{\mu})\in\Lambda(p,q,\bm{l},\bm{u},\bm{w},\Delta,\bm{x},\bm{y})$, $(p,q,\bm{l},\bm{u},\bm{w},\Delta,\bm{x},\bm{y},\bm{\lambda},\bm{\mu})$ is a global solution of the associated \MPECQ.
    \item [ b) ] Assume additionally that 
    \begin{itemize}
    \item [ i) ] for all Traveler $i \in \mathcal{I}$, $b_i \geq p_{min}$;
    \item [ii) ] for all TSP $n \in \mathcal{N}_m, m \in \mathcal{M}$, $\beta_{mn} \leq q_{max}$.
    \end{itemize}
    If $(p,q,\bm{l},\bm{u},\bm{w},\Delta,\bm{x},\bm{y},\bm{\lambda},\bm{\mu})$ is a global solution of \MPECQ, then $(p,q,\bm{l},\bm{u},\bm{w},\Delta,\bm{x},\bm{y})$ is a global solution of \SLMFGQ\space (or equivalently SLMFG with network effects in Proposition \ref{prop:SLMFG_quadra}).
\end{itemize}
\end{theorem}
}

Each of \MPECL\ and \MPECQ\ contain eight nonlinear constraints (four for traveler's problems and four for TSP's problems) corresponding to  complementary slackness conditions.

\subsubsection{SD-based single-level formulation.\label{SD}} An alternative reformulation consists of replacing the follower problems by their strong-duality (SD) optimality conditions.

The SD conditions of \SLMFGL\space are given in \eqref{DU2aa}-\eqref{DU2f}:
\begin{subequations}
\allowdisplaybreaks
\begin{align}
&\text{Traveler $i$'s constraints}~\eqref{BL1b}-\eqref{BL1d}, &&  \forall i \in \mathcal{I},\tag{DL.1}\label{DU2aa}\\
&(r_{i}-b_{i})Q_{i}x_{i}\geq B_{i}\bar{\lambda}_{1}^{i}+u_{i}\bar{\lambda}_{2}^{i}, &&  \forall i \in \mathcal{I},\tag{DL.2}\label{DU2a} \\
&b_{i}Q_{i}\bar{\lambda}_{1}^{i}+\bar{\lambda}_{2}^{i}\geq (r_{i}-b_{i})Q_{i},&&  \forall i \in \mathcal{I},\tag{DL.3} \label{DU2b}\\
&\bar{\lambda}_{1}^{i},\bar{\lambda}_{2}^{i}\geq 0,  &&\forall i \in \mathcal{I},\tag{DL.4} \label{DU2c} \\
&\text{TSP $mn$'s constraints}~\eqref{BL2b}-\eqref{BL2d},&& \forall m \in \mathcal{M},n \in \mathcal{N}_{m},\tag{DL.5} \label{DU2dd}\\
&(\beta_{mn}-\rho_{mn})C_{mn}y_{mn}\geq \bar{B}_{mn}\bar{\mu}_{1}^{mn}+w_{mn}\bar{\mu}_{2}^{mn},&& \forall m \in \mathcal{M},n \in \mathcal{N}_{m},\tag{DL.6} \label{DU2d}\\
&C_{mn}\gamma_{mn}\bar{\mu}_{1}^{mn}+\bar{\mu}_{2}^{mn}\geq (\beta_{mn}-\rho_{mn})C_{mn},&& \forall m \in \mathcal{M},n \in \mathcal{N}_{m},\tag{DL.7}\label{DU2e}\\
&\bar{\mu}_{1}^{mn},\bar{\mu}_{2}^{mn}\geq 0 , && \forall m \in \mathcal{M},n \in \mathcal{N}_{m},\tag{DL.8} \label{DU2f}
\end{align}
\end{subequations}
where $\bar{\lambda}_{1}^{i},\bar{\lambda}_{2}^{i}$ denotes the dual variables corresponding to the primal constraint \eqref{BL1b} and \eqref{BL1c}, respectively. $\bar{\mu}_{1}^{mn},\bar{\mu}_{2}^{mn}$  denotes the dual variables corresponding to primal constraint \eqref{BL2b} and \eqref{BL2c}, respectively. \eqref{DU2a} compares the objective value of Traveler $i$'s primal problem on the left hand side with that of its dual problem on the right hand side, \eqref{DU2b} is the dual constraint of Traveler $i$; \eqref{DU2d} compares the objective value of the TSP $mn$'s primal problem on the left hand side with its dual problem on the right hand side and \eqref{DU2e} is the dual constraint of TSP $mn$. The resulting SD reformulation of \SLMFGL\space is denoted by \SDL\space and is summarized in Model \V.
\newpage
\textbf{Model 3.1 (SD-L)}
\begin{subequations}
\allowdisplaybreaks
\begin{align}
&F=\max_{p,q,\bm{l},\bm{u},\bm{w},\Delta,\bm{x},\bm{y},\bm{\bar{\lambda}},\bm{\bar{\mu}},}\quad\text{MaaS regulator's profits}\quad\eqref{B1a} \nonumber \\
&\text{subject to}  && \nonumber \\
&\text{MaaS regulator's constraints}\quad\eqref{B1b}-\eqref{B1p},\nonumber\\
&\text{Traveler $i$'s duality conditions}\quad\text{\eqref{DU2aa}-\eqref{DU2c}},\nonumber\\
&\text{TSP $mn$'s duality conditions}\quad\text{\eqref{DU2dd}-\eqref{DU2f}}.\nonumber
\end{align}
\end{subequations}

It is well-known that, for linear bilevel optimization problems, single-level reformulations based on SD conditions are equivalent to single-level reformulations based on KKT conditions \citep{fortuny1981representation,zare2019note,kleinert2020closing}. Hence, a direct extension of Theorem \ref{thm:equiv_MILBP} is that, assuming the same conditions as stated therein, the solutions of \SDL\ are also solutions of \SLMFGL\ and reciprocally. This result is summarized in Corollary \ref{thm:SD_MILBP} in Online Appendix \ref{corollary1}.

We next derive the SD reformulation of \SLMFGQ, which is more technical since follower problems therein are strictly convex quadratic (as opposed to linear in \SLMFGL). We give primal-dual formulations of \SLMFGQ\space in Lemma \ref{Lemma1} in Online Appendix \ref{LEMMA1}.

Based on Lemma \ref{Lemma1}, SD conditions of \SLMFGQ\ are summarized in Eqs.\eqref{DUa}-\eqref{DUf}. 
\begin{subequations}
\allowdisplaybreaks
\begin{align}
&\text{Traveler $i$'s constraints}~\eqref{BL1b}-\eqref{BL1d},\forall i \in \mathcal{I},\tag{DQ.1}\label{DUa}\\
&\frac{\zeta\alpha_{i}Q_{i}^{2}}{C^{2}}x_{i}^{2}-\left(\frac{2\zeta\alpha_{i}Q_{i}}{\bar{C}}-\frac{2\zeta\alpha_{i}\Delta}{C^{2}}+r_{i}-b_{i}\right)Q_{i}x_{i}\leq-\frac{(A_{i}\bar{\lambda}_{1}^{i})^{2}+(\bar{\lambda}_{2}^{i})^{2}}{2P_{i}}-\frac{U_{i}}{P_{i}}\left(b_{i}Q_{i}\bar{\lambda}_{1}^{i}+\bar{\lambda}_{2}^{i}\right)\nonumber\\
&-\frac{b_{i}Q_{i}\bar{\lambda}_{1}^{i}\bar{\lambda}_{2}^{i}}{P_{i}}-\frac{U_{i}^{2}}{2P_{i}}-B_{i}\bar{\lambda}_{1}^{i}-\bar{\lambda}_{2}^{i}u_{i},\forall i \in \mathcal{I},\tag{DQ.2}\label{DUb} \\
&\bar{\lambda}_{1}^{i},\bar{\lambda}_{2}^{i}\geq 0, \forall i \in \mathcal{I}, \tag{DQ.3} \label{DUc} \\
&\text{TSP $mn$'s constraints}~\eqref{BL2b}-\eqref{BL2d},\forall m \in \mathcal{M},n \in \mathcal{N}_{m},\tag{DQ.4}\label{DUd}\\
&\frac{\kappa\eta_{mn} C_{mn}^{2}}{C^{2}}y_{mn}^{2}-\left(\frac{2\kappa\eta_{mn}\Delta}{C^{2}}+\beta_{mn}-\rho_{mn}\right)C_{mn}y_{mn}\leq-\frac{(A_{mn}\bar{\mu}_{1}^{mn})^{2}+(\bar{\mu}_{2}^{mn})^{2}}{2P_{mn}}\nonumber\\
&-\frac{U_{mn}}{P_{mn}}\left(A_{mn}\bar{\mu}_{1}^{mn}+\bar{\mu}_{2}^{mn}\right)-\frac{A_{mn}\bar{\mu}_{1}^{mn}\bar{\mu}_{2}^{mn}}{P_{mn}}-\frac{U_{mn}^{2}}{2P_{mn}}-\bar{B}_{mn}\bar{\mu}_{1}^{mn}-w_{mn}\bar{\mu}_{2}^{mn},\forall m \in \mathcal{M},n \in \mathcal{N}_{m},\tag{DQ.5}\label{DUe} \\
&\bar{\mu}_{1}^{mn},\bar{\mu}_{2}^{mn}\geq 0, \forall m \in \mathcal{M},n \in \mathcal{N}_{m}, \tag{DQ.6}\label{DUf}
\end{align}
\end{subequations}
where \eqref{DUb} compares the objective value of Traveler $i$'s primal problem on the left hand side with that of Traveler $i$'s dual follower problem on the right hand side, and \eqref{DUe} compares the objective value of TSP $mn$'s primal problem on the left hand side with that of TSP $mn$'s dual follower problem on the right hand side.

Unlike the SD conditions of the follower problems of \SLMFGL, SD conditions of the follower problems of  \SLMFGQ\xspace involve bilinear terms which are product of continuous variables, namely, $\Delta \cdot x_{i}$ and $\Delta\cdot y_{mn}$ in \eqref{DUb} and \eqref{DUe}, respectively. We relax these bilinear terms by introducing auxiliary variables  $W_{i} \triangleq \Delta\cdot x_{i}$ and $W_{mn}\triangleq\Delta\cdot y_{mn}$  using their McCormick envelopes \citep{mccormick1976computability}. Accordingly, constraints \eqref{MKa}-\eqref{MKh} are added to the SD reformulation of \SLMFGQ.

The SD-based reformulation of \SLMFGQ\space is denoted \SDQ\space and summarized in Model \VI.\\
\textbf{Model 3.2 (SD-Q)}
\begin{subequations}
\allowdisplaybreaks
\begin{align}
&\max_{p,q,\bm{l},\bm{u},\bm{w},\Delta,\bm{x},\bm{y},\bm{\bar{\lambda}},\bm{\bar{\mu}}}\quad \text{MaaS regulator's profits}\quad \eqref{B1a},\nonumber\\
&\text{\emph{subject to}}  && \nonumber \nonumber\\
&\text{MaaS regulator's constraints}\quad\text{\eqref{B1b}-\eqref{B1p}},\nonumber\\
&\text{Traveler $i$'s duality conditions }\quad\eqref{DUa}-\eqref{DUc},\eqref{MKa}-\eqref{MKd},\nonumber\\
&\text{TSPs' duality conditions }\quad\eqref{DUd}-\eqref{DUf},\eqref{MKe}-\eqref{MKh}.\nonumber\nonumber
\end{align}
\end{subequations}

Since the duality gap for convex optimization problems is null,  SD conditions are equivalent to KKT conditions, Hence, analogously to Corollary \ref{thm:SD_MILBP}, a direct extension of Theorem \ref{thm:equiv_MIQBP} is given in Corollary \ref{thm:SD_MIQBP} in Online Appendix \ref{corollary2}.

Observe that \SDQ\xspace is a relaxation of \SLMFGQ\ whereas \SDL\ is an exact reformulation of \SLMFGL. We next present the customized B\&B algorithms based on formulations \SDL\xspace and \SDQ\xspace to solve \SLMFGL\xspace and \SLMFGQ, respectively.

\subsection{SD-based B\&B algorithm}
\label{algorithm}

In \SDL\xspace and \SDQ, the nonlinear terms $\bar{\lambda}_{2}^{i}u_{i}$ in \eqref{DU2a} and \eqref{DUb}, and $\bar{\mu}_{2}^{mn}w_{mn}$ in \eqref{DU2d} and \eqref{DUe} could be linearized using classic big-$M$ approaches, however, finding valid big-$M$ values is non-trivial \citep{kleinert2020there}. To circumvent this difficulty , we propose a customized B\&B algorithm, \SDBB, which branches on the binary variables $u_{i}$ and $w_{mn}$.

In \SDL, \SDBB\ algorithm starts by relaxing all nonlinear constraints \eqref{DU2a} and \eqref{DU2d}. This relaxed problem is thus a MILBP which ignores the optimality conditions of followers. This relaxation is also known as the High-Point-Relaxation (HPR). Upon solving the HPR, the optimality conditions of followers are examined and a sub-optimal follower is selected (otherwise the solution is bilevel-optimal). If the sub-optimal follower represents a traveler $i' \in \mathcal{I}$, two children nodes are generated with the constraints $u_{i'}=0$ and $u_{i'}=1$, respectively, and the corresponding constraint \eqref{DU2a} is added to both children nodes. Otherwise, the sub-optimal follower must represent a TSP $mn' \in \mathcal{MN}$ and two children nodes are generated with the constraints $w_{mn'}=0$ and $w_{mn'}=1$, respectively, and the corresponding constraint \eqref{DU2d} is added to both children nodes. Observe that the added constraints \eqref{DU2a} and \eqref{DU2d} in the generated children nodes are linear since the binary variables therein are fixed. Analogically, in \SDQ, \SDBB\ algorithm starts by solving the HPR problem, which is a MIQBP. If the sub-optimal follower represents a traveler $i' \in \mathcal{I}$, two children nodes are generated and the corresponding constraint \eqref{DUb} \HX{together with its McCormick envelopes \eqref{MKa}-\eqref{MKd}} are added to both children nodes. Otherwise, the sub-optimal follower must represent a TSP $mn' \in \mathcal{MN}$, two children nodes are generated and the corresponding constraint \eqref{DUe} \HX{together with its McCormick envelopes \eqref{MKe}-\eqref{MKh}} are added to both children nodes. Observe that the nonlinear terms $\bar{\lambda}_{2}^{i}u_{i}$ and $\bar{\mu}_{2}^{mn}w_{mn}$ of the constraints added in the generated children nodes are linearized since the binary variables are therein fixed. 

We next introduce some formal notations to present the proposed \SDBB\ algorithm.
\subsubsection{Follower optimality.}
Let $\overline{\mathcal{I}}_{0} =\left\{i \in \mathcal{I} : u_{i}=0\right\}$ and $\overline{\mathcal{I}}_{1}=\left\{i \in \mathcal{I} : u_{i}=1\right\}$ be the sets of accepted and rejected travelers, respectively. These sets are initialized to be empty sets at the root node of the B\&B tree. The set of undecided travelers is $\widehat{\mathcal{I}}=\mathcal{I} \setminus \{\mathcal{I}_{0} \cup \mathcal{I}_{1}\}$. Analogously, we define the sets of accepted/rejected TSPs as $\overline{\mathcal{MN}}_{0}=\left\{(m,n) \in \mathcal{MN} : w_{mn}=0\right\}$ and $\overline{\mathcal{MN}}_{1}=\left\{(m,n) \in \mathcal{MN} : w_{mn}=1\right\}$, respectively; and the set of undecided TSPs as $\widehat{\mathcal{MN}}=\mathcal{MN} \setminus \{\overline{\mathcal{MN}}_{0} \cup \overline{\mathcal{MN}}_{1}\}$.

To check the optimality of follower problems, we solve a series of single-variable linear programs (one per follower) parameterized by the leader solution $(\bm{u},\bm{w},\Delta)$. Specifically, at any node of the B\&B tree, for all travelers $i \in \widehat{\mathcal{I}}$, we solve \eqref{BL1a}-\eqref{BL1d} and compare the obtained objective value $G_{i}^{*}$ (optimal response) with the evaluation of \eqref{BL1a}, denoted $\hat{G}_{i}$, based on the relaxed solution obtained at this node. If $|G_{i}^{*}-\hat{G}_{i}|> \epsilon$, for a tolerance $\epsilon$, e.g., $\epsilon=1\times 10^{-4}$, then traveler $i$ is considered sub-optimal and added to the set $\mathcal{V}_{travel}$ of sub-optimal travelers. The same process is implemented for TSPs. The pseudocode of is summarized in \Alg\ref{alg1}.

\begin{algorithm}[ht!]
\footnotesize
\KwIn{$\Delta, \bm{u},\bm{w}$}
\KwOut{$ \mathcal{V}_{travel},\mathcal{V}_{TSP},O_{i},O_{mn}$}
	\For{$i \in \mathcal{I}$ }{		
		$G_{i}^{*}\leftarrow$  solve Traveler $i$'s follower problem 	\eqref{BL1a}-\eqref{BL1d}, given $u_{i},\Delta$\\
		$\hat{G}_{i}\leftarrow$ substitute $x_{i}$ and $\Delta$ into Traveler $i$'s objective function \eqref{BL1a}\\	
	   $O_{i}\leftarrow |G_{i}^{*}-\hat{G}_{i}|$\\
		\If{$|G_{i}^{*}-\hat{G}_{i}|> \epsilon$}{$\mathcal{V}_{travel}\leftarrow\mathcal{V}_{travel}\cup \left\{i \right\}$}
		}
     \For{(m,n) in $\mathcal{MN}$ }{		
		$H_{mn}^{*}\leftarrow$  solve TSP $mn$'s follower problem 	\eqref{BL2a}-\eqref{BL2d}, given $w_{mn}, \Delta$\\
		$\hat{H}_{mn}\leftarrow$ substitute $y_{mn}$ into TSP $mn$'s objective function \eqref{BL2a} \\
		 $O_{mn}\leftarrow |H_{mn}^{*}-\hat{H}_{mn}|$\\	
		\If{$|H_{mn}^{*}-\hat{H}_{mn}|>\epsilon$}{$\mathcal{V}_{TSP}\leftarrow\mathcal{V}_{TSP}\cup \left\{(m,n)\right\}$}	
		}
	\caption{Updating the non-optimal sets of followers}
	\label{alg1}
\end{algorithm}

\subsubsection{Branching rules.}

\label{branching}
We first branch on TSPs accept/reject variables $\bm{w}$ before branching on travelers accept/reject variables $\bm{u}$. This is motivated by empirical evidence and the fact that in realistic instances the number of TSPs is expected to be substantially smaller than the number of travelers. We customize three types of branching rules to select the branching variable among suboptimal followers:\\
\indent 1) \BP: If $\widehat{\mathcal{MN}}\neq \varnothing$, select a branching variable $w_{mn}$ corresponding to TSP $mn$ with the minimum sell-bidding price $\beta_{mn}$, $\forall (m,n)\in \mathcal{V}_{TSP}$. Otherwise, if $\widehat{\mathcal{I}}\neq \varnothing$, select a branching variable $u_{i}$ corresponding to Traveler $i$ with the maximum purchase-bidding price $b_{i}$, $\forall i\in \mathcal{V}_{travel}$. 

2) \Diff: If $\widehat{\mathcal{MN}}\neq \varnothing$, select a branching variable $w_{mn}$ corresponding to TSP $mn$ with the maximum optimality gap $|H_{mn}^{*}-\hat{H}_{mn}|$, $\forall (m,n)\in \mathcal{V}_{TSP}$. Otherwise, if $\widehat{\mathcal{I}}\neq \varnothing$, select a branching variable $u_{i}$ corresponding to Traveler $i$ with the maximum optimality gap  $|G_{i}^{*}-\hat{G}_{i}|, \forall i\in \mathcal{V}_{travel}$.

3) \Weight: combine the above two rules by introducing a weighting factor $\theta$. If $\widehat{\mathcal{MN}}\neq \varnothing$, select a branching variable corresponding to TSP $mn$ with the maximum weighted value $\theta|H_{mn}^{*}-\hat{H}_{mn}|+(1-\theta) \beta_{mn}$. Otherwise, if $\widehat{\mathcal{I}}\neq \varnothing$, select a branching variable corresponding to Traveler $i$ with the maximum weighted value $\theta|G_{i}^{*}-\hat{G}_{i}|+(1-\theta) b_{i}$.

The pseudocode of the branching rules is summarized in  \Alg\ref{alg3}.
\begin{algorithm}[ht!]
\footnotesize
\KwIn{$\widehat{\mathcal{MN}},\widehat{\mathcal{I}},\mathcal{V}_{travel},\mathcal{V}_{TSP},b_{i},\beta_{mn},O_{i},O_{mn}$}
\KwOut{$ \overline{\mathcal{I}}_{0}$, $ \overline{\mathcal{I}}_{1}, \overline{\mathcal{MN}}_{0}, \overline{\mathcal{MN}}_{1},\mathcal{L}$ }
		\If{$\widehat{\mathcal{MN}}\neq \varnothing$}{
			{\If{\emph{\text{branching rule} is \BP}}{
	$(m,n)\leftarrow argmin_{(m,n)\in\mathcal{V}_{TSP}} \left\{\beta_{mn}\right\}$ \COMMENT{\textcolor{OliveGreen}{select $(m,n)$ of $\mathcal{V}_{TSP}$ with minimum bidding price }} \\}
	\If{\emph{\text{branching rule} is \Diff}}{
	$(m,n)\leftarrow argmax_{(m,n)\in \mathcal{V}_{TSP}}\left\{O_{mn}\right\}$ \COMMENT{\textcolor{OliveGreen}{select $(m,n)$ of $\mathcal{V}_{TSP}$  with maximum optimality gap }}  \\}
	\If{\emph{\text{branching rule} is \Weight}}{
	$Z_{mn}\leftarrow \theta\cdot O_{mn}+(1-\theta)\cdot\beta_{mn}$\\
	$(m,n)\leftarrow argmax_{(m,n)\in\mathcal{V}_{TSP}} \left\{Z_{mn}\right\}$ \COMMENT{\textcolor{OliveGreen}{select $(m,n)$ of $\mathcal{V}_{TSP}$  with maximum weighted value}} \\}}
		$s_{2}\leftarrow(m,n) $\\
		 $\overline{\mathcal{MN}}_{0}\leftarrow\overline{\mathcal{MN}}_{0}\cup\{s_{2}\}$\COMMENT{\textcolor{OliveGreen}{add the new branching index to the fixed set $\overline{\mathcal{MN}}_{0}$}}\\
		 $\overline{\mathcal{MN}}_{1}\leftarrow\overline{\mathcal{MN}}_{1}\cup\{s_{2}\}$\COMMENT{\textcolor{OliveGreen}{add the new branching index to the fixed set $\overline{\mathcal{MN}}_{1}$}}\\
		 $n \leftarrow \text{len}(\mathcal{L})$\COMMENT{\textcolor{OliveGreen}{update the left child node}}\\
		  $\mathcal{L} \leftarrow \mathcal{L} \cup\{n\}$\COMMENT{\textcolor{OliveGreen}{update the set of active nodes}}\\
		 }
	\ElseIf{$\widehat{\mathcal{I}}\neq \varnothing$}
	{\If{\emph{\text{branching rule} is \BP}}{
	$i\leftarrow argmax_{i\in\mathcal{V}_{travel}} \left\{b_{i}\right\}$ \COMMENT{\textcolor{OliveGreen}{select $i$ of $\mathcal{V}_{travel}$ with the maximum  bidding price}}  \\}
	\If{\emph{\text{branching rule} is \Diff}}{
	$i\leftarrow argmax_{i\in \mathcal{V}_{travel}}\left\{O_{i}\right\}$ \COMMENT{\textcolor{OliveGreen}{select $i$ of $\mathcal{V}_{travel}$ with the maximum optimality gap}}  \\}
	\If{\emph{\text{branching rule} is \Weight}}{
	$Z_{i}\leftarrow \theta\cdot O_{i}+(1-\theta)\cdot b_{i}$\\
	$i\leftarrow argmax_{i\in \mathcal{V}_{travel}} \left\{Z_{i}\right\}$ \COMMENT{\textcolor{OliveGreen}{select $i$ of $\mathcal{V}_{travel}$ with the maximum weighted value}}  \\}
	$s_{1}\leftarrow i $\\
	$\overline{\mathcal{I}}_{0}\leftarrow \overline{\mathcal{I}}_{0}\cup\{s_{1}\}$\COMMENT{\textcolor{OliveGreen}{add new branching index to the fixed set $\overline{\mathcal{I}}_{0}$}}\\
	$\overline{\mathcal{I}}_{1}\leftarrow \overline{\mathcal{I}}_{1}\cup\{s_{1}\}$\COMMENT{\textcolor{OliveGreen}{add new branching index to the fixed set $ \overline{\mathcal{I}}_{1}$}}\\
			 $n\leftarrow n+1$\COMMENT{\textcolor{OliveGreen}{update the right child node}}\\
		  $\mathcal{L} \leftarrow \mathcal{L} \cup\{n\}$\COMMENT{\textcolor{OliveGreen}{update the set of active nodes}}\\
		} 
	\caption{Branching procedure}
	\label{alg3}
\end{algorithm}
\subsubsection{Algorithm summary.} The proposed \SDBB\space algorithm is summarized in \Alg\ref{alg4}. $\textit{SP}$ denotes a sub-problem at a node of the B\&B tree and index $k$ is a superscript to track the current B\&B node in the \textbf{while} loop, the formulations of the sub-problem $\textit{SP}^{k}$ of \SLMFGL\ and \SLMFGQ\ are provided in Online Appendix \ref{sub-problems formulation}.  The list of active sub-problem indices is denoted $\mathcal{L}$ and initialized with the root node sub-problem, i.e. the HPR of the \SDL\ or \SDQ. $\textit{UB}$ and $\textit{LB}$ are used to track upper and lower bounds throughout the process. 
\begin{algorithm}[ht!]
\KwIn{$\mathcal{B}_{i},\mathcal{B}_{mn}$}
\KwOut{$k,\overline{\textit{UB}},LB, Gap$}
\footnotesize
    $ Gap \leftarrow \infty$, $\text{LB}\leftarrow -\infty$, $n\leftarrow 0$ $\mathcal{L}\leftarrow \left\{n\right\}$ \COMMENT{\textcolor{OliveGreen}{initialization}}\\
     $ \overline{\mathcal{I}}_{0}^{k}\leftarrow \left\{\right\}$, $ \overline{\mathcal{I}}_{1}^{k}\leftarrow \left\{\right\}, \overline{\mathcal{MN}}_{0}^{k}\leftarrow \left\{\right\}, \overline{\mathcal{MN}}_{1}^{k}\leftarrow \left\{\right\}$ \COMMENT{\textcolor{OliveGreen}{initialize the fixed sets}}  \\
  $\textit{SP}^{0}\leftarrow \text{HPR}\left(\mathcal{B}_{i},\mathcal{B}_{mn}\right)$,\COMMENT{\textcolor{OliveGreen}{solve high-point-relaxation problem}}\\
    $\overline{\textit{UB}}\leftarrow F^{0}(p^{0},q^{0},\bm{l}^{0},\bm{u}^{0},\bm{w}^{0},\Delta^{0},\bm{x}^{0},\bm{y}^{0})$  \COMMENT{\textcolor{OliveGreen}{the objective value of HPR}} \\
    $k \leftarrow 1$\\
	\While{$\mathcal{L} \neq \varnothing$ }{
	\textit{parent node} $\leftarrow$ select a node from $\mathcal{L}$\COMMENT{\textcolor{OliveGreen}{best bound first search}}\\
		{$F^{k}$ $\leftarrow\textit{SP}^{k}$  $\left(\mathcal{B}_{i},\mathcal{B}_{mn},\textit{SP}^{k-1},\overline{\mathcal{I}}_{0}^{k}, \overline{\mathcal{I}}_{1}^{k}, \overline{\mathcal{MN}}_{0}^{k}, \overline{\mathcal{MN}}_{1}^{k}
\right)$\COMMENT{\textcolor{OliveGreen}{solve the $kth$ sub-problem}}\\	
	\If{\emph{$\textit{SP}^{k}$ is infeasible}}{ 
	   $\textit{UB}^{k}\leftarrow F^{k}(p^{k},q^{k},\bm{l}^{k},\bm{u}^{k},\bm{w}^{k},\Delta^{k},\bm{x}^{k},\bm{y}^{k})$\COMMENT{\textcolor{OliveGreen}{obtain \textit{UB}}}\\
	    $\mathcal{V}_{travel}^{k},\mathcal{V}_{TSP}^{k},O_{i}^{k},O_{mn}^{k}\leftarrow  \textbf{\Alg}\ref{alg1}$ $(\Delta^{k},\bm{u}^{k},\bm{w}^{k})$\\
	 { \If{all the followers are optimal ($\mathcal{V}_{travel}^{k}=\varnothing,\mathcal{V}_{TSP}^{k}=\varnothing$)}{ \If{$\textit{UB}^{k}>\textit{LB}$}
		{ $\textit{LB}\leftarrow \textit{UB}^{k}$\COMMENT{\textcolor{OliveGreen}{update LB }}\\	prune the nodes	of $\mathcal{L}$ with $\textit{UB}\leq\textit{LB}$ \COMMENT{\textcolor{OliveGreen}{fathoming }}
	} 
		}
		}      \Else{
           $\overline{\mathcal{I}}_{0}^{k}$, $ \overline{\mathcal{I}}_{1}^{k}, \overline{\mathcal{MN}}_{0}^{k}, \overline{\mathcal{MN}}_{1}^{k},\mathcal{L}\leftarrow$	\textbf{\Alg}\ref{alg3}  ($\widehat{\mathcal{MN}}^{k},\widehat{\mathcal{I}}^{k},\mathcal{V}_{TSP}^{k}$,$\mathcal{V}_{travel}^{k},b_{i},\beta_{mn},O_{i}^{k},O_{mn}^{k}$)\\
		}  }
		{{

              } }}
             $\mathcal{L}\leftarrow \mathcal{L}-\left\{\textit{parent node} \right\}$,\COMMENT{\textcolor{OliveGreen}{delete the selected node from $\mathcal{L}$}} \\
              \If{$\mathcal{L}=\varnothing$}{$\overline{\textit{UB}}\leftarrow LB$\\ $Gap\leftarrow 0$\COMMENT{\textcolor{OliveGreen}{converged}}}
              \Else{$\overline{\textit{UB}}\leftarrow$ select an active node of $\mathcal{L}$ with the maximum \textit{UB}\\ \If{($\overline{\textit{UB}}-\textit{LB})/\textit{LB}<Gap$}{$Gap=(\overline{\textit{UB}}-\textit{LB})/\textit{LB}$\COMMENT{\textcolor{OliveGreen}{update the gap}}}}
              \If{$Gap<1\times e^{-6}$}{break\COMMENT{\textcolor{OliveGreen}{exit by gap}}}
             $k\leftarrow k+1 $         
              }
	\caption{SD-based Branch and bound algorithm (\SDBB)}
	\label{alg4}
\end{algorithm}
\subsection{Illustration of Algorithm SD-B\&B}
\label{llustrative example}
\begin{figure}[ht!]
\caption{ Search tree for the example. The node expansion order is $\bm{0\rightarrow1\rightarrow 2 \rightarrow3\rightarrow 4\rightarrow 7\rightarrow 8\rightarrow 9\rightarrow 10}.$}
\includegraphics[width=0.95\textwidth]{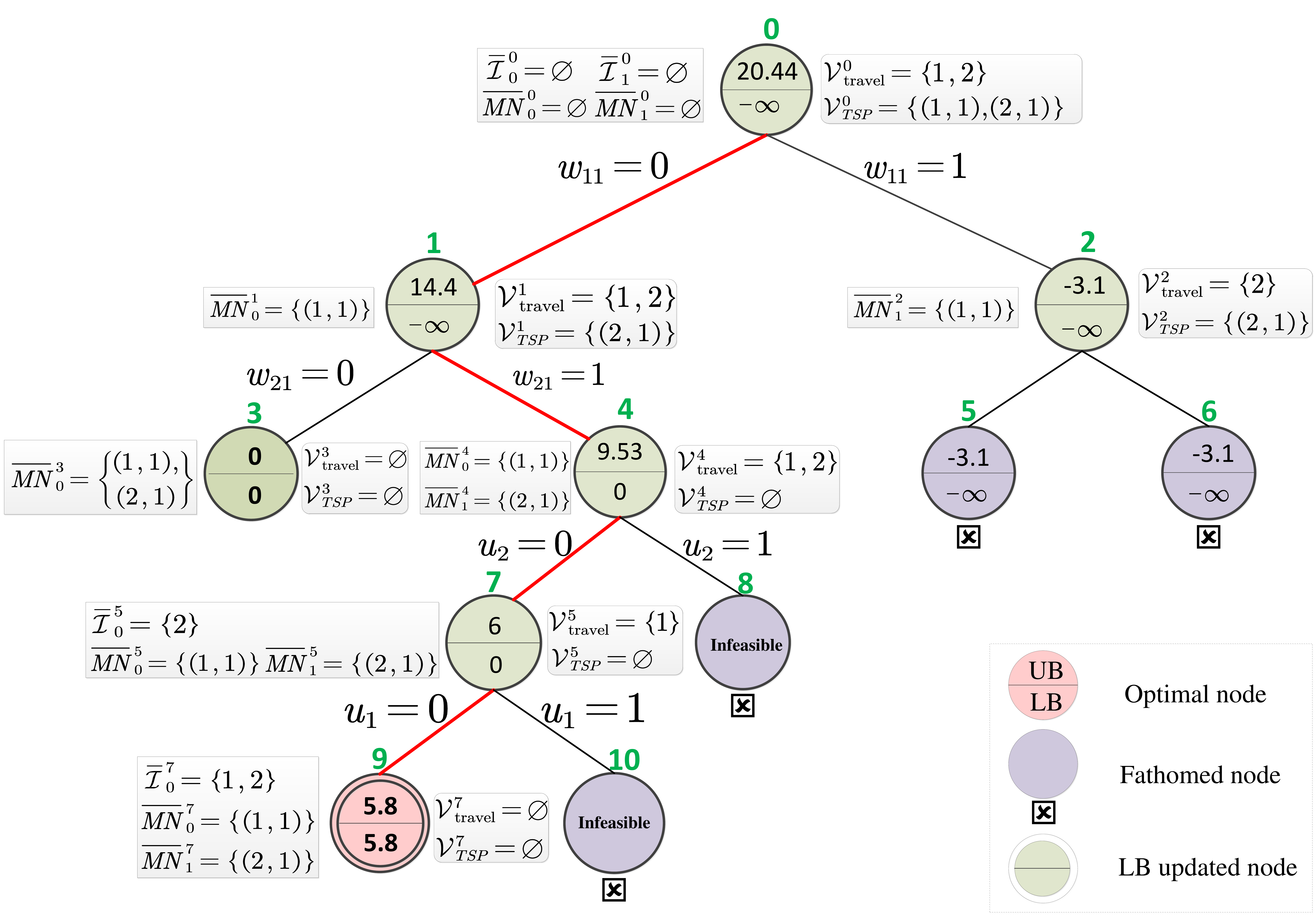}
\label{tree}
\end{figure}

We illustrate Algorithm \SDBB\ using the same data as in Example \textcolor{blue}{1}, we consider the linear case, i.e. \SLMFGL\space (MILBP) and solve it using \SDBB\space and branching rule \BP. The corresponding B\&B tree is shown in Figure \ref{tree}, and the formulation of the $kth$ sub-problem $\textit{SP}^{k}$ of \SLMFGL\ is provided in Online Appendix \ref{sub-problems formulation}.
At each iteration, the check is made to see if each user and each TSP is in an optimal condition, namely, $\mathcal{V}_{travel}=\varnothing$ and $\mathcal{V}_{TSP}=\varnothing$.\\

At the root node, we solve the HPR problem to obtain $\bm{u}^{0}=(1,1)$, $\bm{w}^{0}=(1,1)$, $\bm{x}^{0}=(1,0.68)$, $\bm{y}^{0}=(0.34,0.32)$ and an objective value  $F^{0}=20.44$. Since no follower is optimal, the accept/reject binary variable corresponding to the TSP with the minimum sell-bid is selected as the branching variable, i.e. $w_{11}$ and two children nodes are generated accordingly. Solving sub-problems $SP^1$ and $SP^2$ yields $F^1=14.4$ and $F^2=-3.1$, respectively. Follower-optimality is not satisfied, hence both sub-problems generate two children nodes by branching on the other TSP accept/reject variable, i.e. $w_{21}$. Solving $\textit{SP}^{3}$ yields $F^{3}=0$ and all followers are optimal, i.e. $\mathcal{V}_{travel}^{3}=\mathcal{V}_{TSP}^{3}=\varnothing$. Hence, \textit{LB} is updated to 0 and since $LB>F^{2}$, the children nodes of node 2 (nodes 5 and node 6) are pruned. Solving $\textit{SP}^{4}$ yields $F^{4}=9.53$, but travelers' problems are sub-optimal, hence the accept/reject variable of the traveler with the maximum purchase-bid, $u_2$, is selected for branching. Solving $\textit{SP}^{5}$ at node 7 yields $F^{5}=6$ and since the last undecided traveler is sub-optimal, two children nodes are generated by branching on $u_1$. The sub-problem at node 8 is infeasible, thus this node is fathomed. Solving $\textit{SP}^{7}$ at node 9 yields $F^{7}=5.8$ which improves on the current incumbent and thus $\textit{LB}$ is updated to 5.8. The sub-problem at node 10 is infeasible, thus this node is fathomed and the optimal solution is that corresponding to node 9.

\section{\HX{Computational study}}
\label{numerical experiments}
In this section, we conduct a series of numerical experiments to evaluate the performance of the proposed SLMFG models and \SDBB\space algorithm. The input data is introduced in Online Appendix \ref{input data}. We examine the impact of followers' bidding price and the capacity of mobility resources on the SLMFG models. We compare the computational performance of \SDBB\ algorithm using the three proposed branching rules against a benchmark. The benchmark is the B\&B algorithm proposed by \cite{bard1990branch} that branches on complementary slackness conditions of the MPEC reformulation, hereby referred to as \Bard.  All numerical experiments are conducted using Python 3.7.4 and CPLEX Python API on a Windows 10 machine with Intel(R) Core i7-8700 CPU $@$ 3.20 GHz, 3.19 GHz, 6 Core(s) and with 64 GB of RAM.
\subsection{Algorithms benchmarking}
\label{benchmarking}

\begin{figure}[!t]
 \centering
 \vspace{0.1in}
 \caption{Average Runtime under different branching rules in $\bm{\mathsf{SD}}$-$\bm{\mathsf{B\&B}}$ algorithm for \SLMFGL\space and \SLMFGQ.}
\includegraphics*[width=0.42\textwidth]{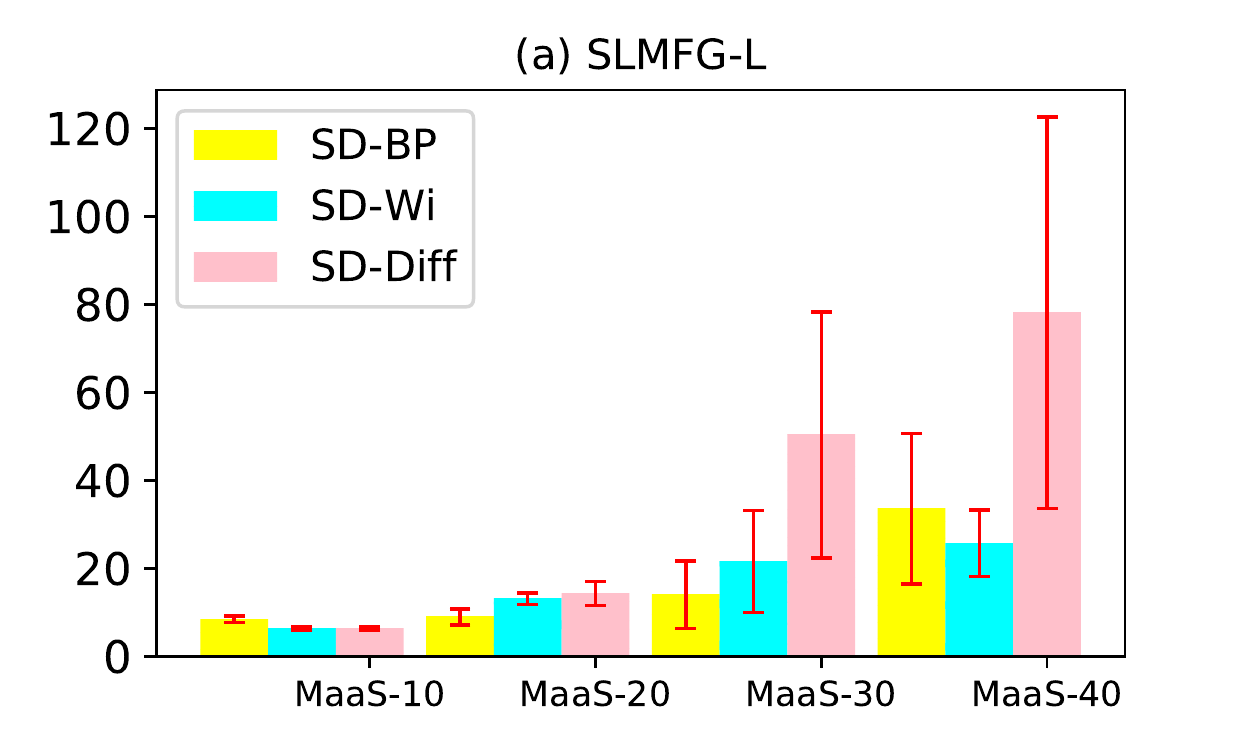}  
\includegraphics*[width=0.42\textwidth]{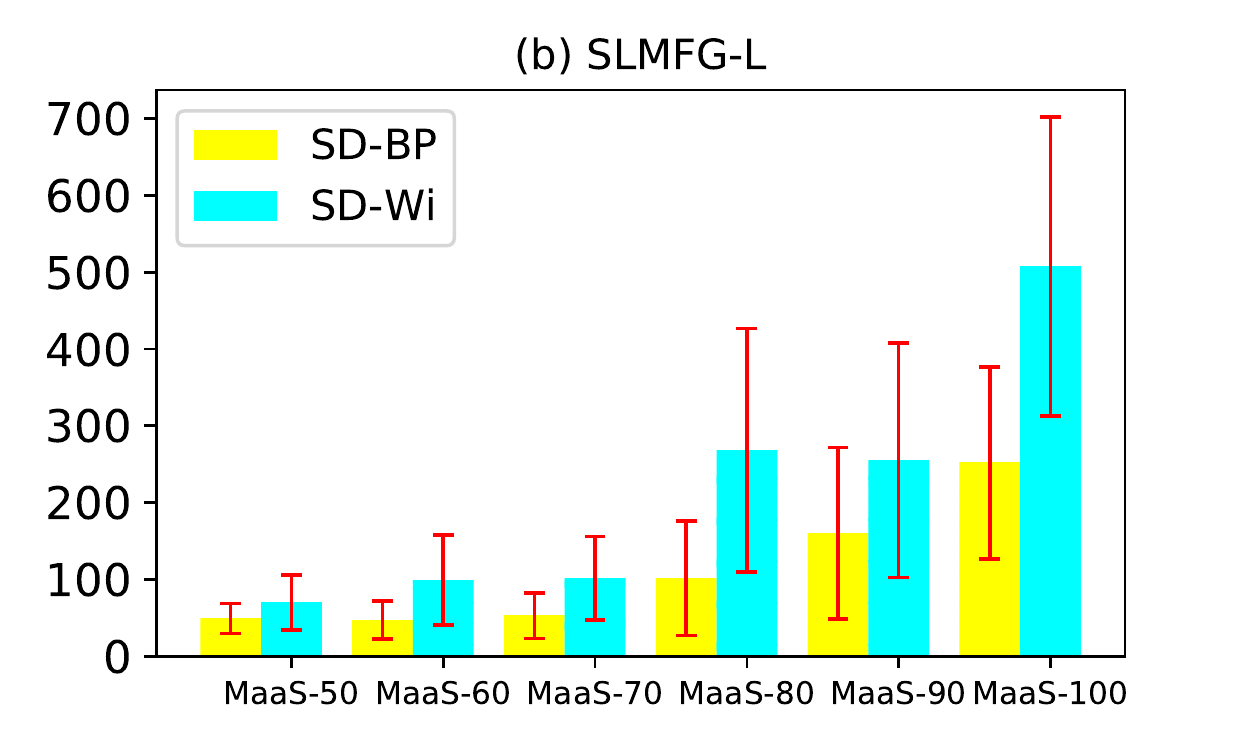} 
\includegraphics*[width=0.42\textwidth]{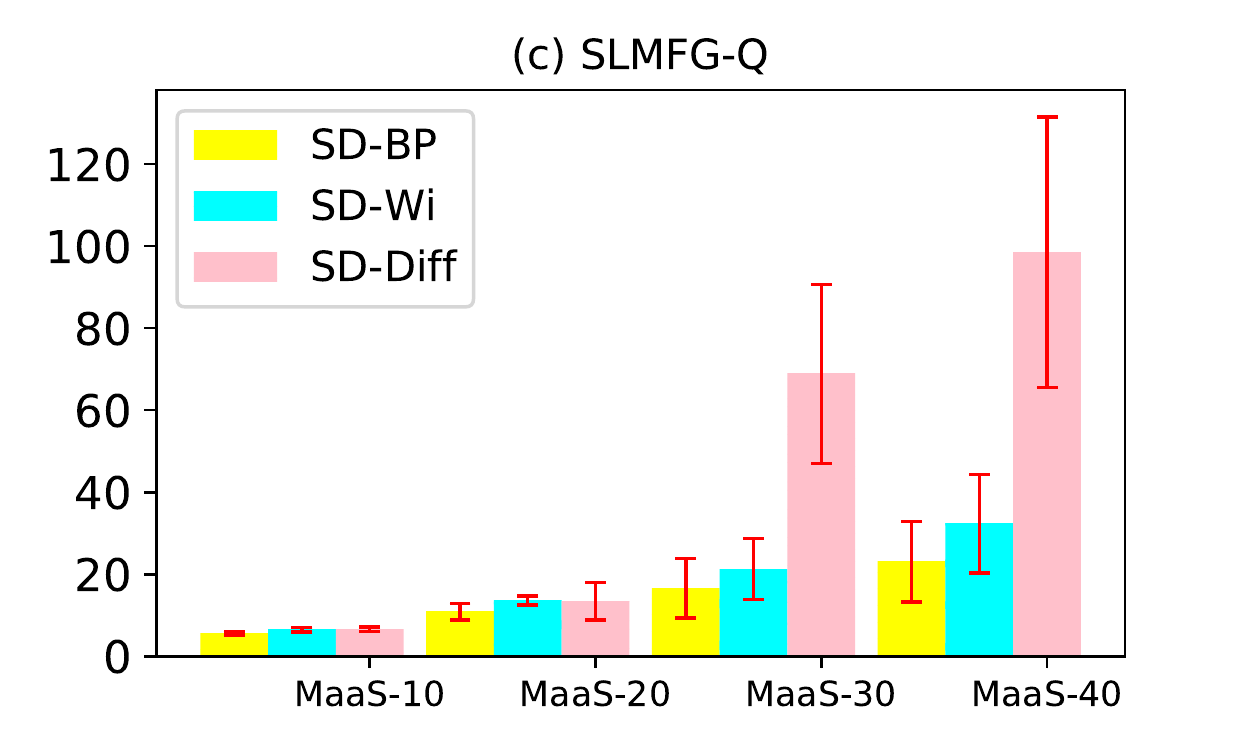}  
\includegraphics*[width=0.42\textwidth]{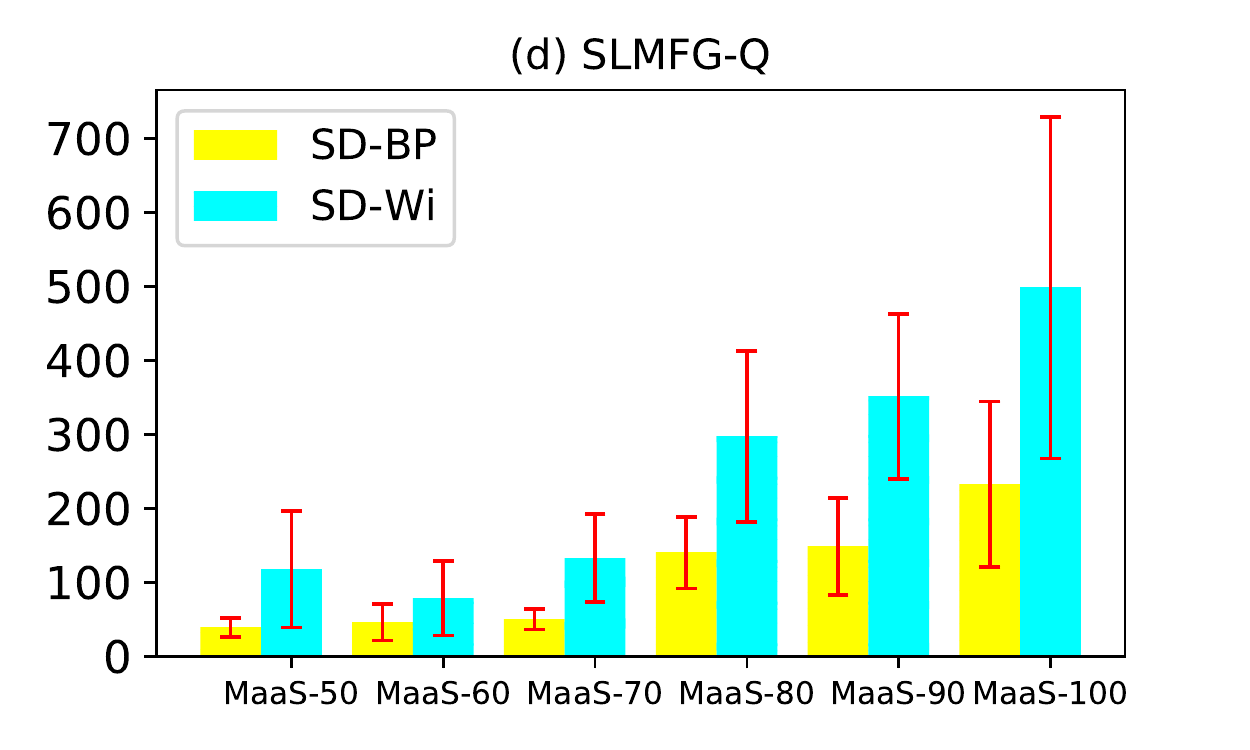} 
 \label{comparebranch}
\end{figure}

In this section, we present numerical results of the proposed \SDBB\space algorithm under three types of branching rules and compare its performance against the benchmark, i.e. Algorithm \Bard. We implement Algorithm \SDBB\xspace using the three proposed branching rule described in Section \ref{branching} and named accordingly: \BP\xspace, \Diff\xspace and \Weight\xspace. We consider 20 instances for each group, hence a total of 400 instances are respectively solved for the linear (\SLMFGL) and quadratic (\SLMFGQ) cases using all four methods (three variants of \SDBB\xspace and \Bard), here we only report detailed results of the first three instances of each group of instances for each method. The numerical results for \SLMFGL\ and \SLMFGQ\ are summarized in \T \ref{Ltable} and \T \ref{Qtable}, respectively, in which each row represents an instance. For each method, the results are organized in four columns: $k$ is the number of sub-problems solved, $\textit{LB}$ is the objective value of the best integer solution found, $\textit{Gap}$ is the optimality gap in percentage and $T$ is the total CPU runtime in seconds. We use a time limit of 10800 seconds (3 hours).
\setlength{\abovecaptionskip}{0pt}
\setlength{\belowcaptionskip}{0pt}
\begin{table}[ht!]
\label{Ltable}
\centering
\scriptsize
\caption{Performance of the proposed B\&B algorithms benchmarked with \Bard\ for \SLMFGL (MILBP)}
\begin{tabular}{@{\extracolsep{\fill}}cccccccccccccccccc@{}}
\toprule
Instances &
   &
  \multicolumn{4}{c}{\Bard} &
  \multicolumn{4}{c}{\BP} &
  \multicolumn{4}{c}{\Diff} &
  \multicolumn{4}{c}{\Weight} \\ \cmidrule(l){3-6} \cmidrule(l){7-10} \cmidrule(l){11-14} \cmidrule(l){15-18}
\multicolumn{1}{c}{\textit{}} &
  \multicolumn{1}{c}{\textit{}} &
  \multicolumn{1}{c}{\textit{k}} &
  \multicolumn{1}{c}{\textit{LB}} &
  \multicolumn{1}{c}{\textit{Gap(\%)}} &
  \multicolumn{1}{c}{\textit{$T (s)$ }} &
  \multicolumn{1}{c}{\textit{k}} &
  \multicolumn{1}{c}{\textit{LB}} &
  \multicolumn{1}{c}{\textit{Gap(\%)}} &
  \multicolumn{1}{c}{\textit{$T (s)$ }} &
  \multicolumn{1}{c}{\textit{k}} &
  \multicolumn{1}{c}{\textit{LB}} &
  \multicolumn{1}{c}{\textit{Gap(\%)}} &
  \multicolumn{1}{c}{\textit{$T (s)$ }} &
  \multicolumn{1}{c}{\textit{k}} &
  \multicolumn{1}{c}{\textit{LB}} &
  \multicolumn{1}{c}{\textit{Gap(\%)}} &
  \multicolumn{1}{c}{\textit{$T (s)$ }} \\ \cmidrule(l){3-6} \cmidrule(l){7-10} \cmidrule(l){11-14} \cmidrule(l){15-18}
\multirow{3}{*}{MaaS-10-5} &
  1 &
  52 &
  19 &
  0 &
  13 &
  32 &
  19 &
  0 &
  8 &
  26 &
  19 &
  0 &
  7 &
  26 &
  19 &
  0 &
  \textbf{6} \\
 &
  2 &
  63 &
  0 &
  0 &
  15 &
  34 &
  0 &
  0 &
  9 &
  26 &
  0 &
  0 &
  7 &
  26 &
  0 &
  0 &
  \textbf{6} \\
 &
  3 &
  47 &
  0 &
  0 &
  12 &
  40 &
  0 &
  0 &
  10 &
  28 &
  0 &
  0 &
  7 &
  28 &
  0 &
  0 &
  \textbf{7} \\ \hline
\multirow{3}{*}{MaaS-30-5} &
  1 &
  76 &
  310 &
  0 &
  24 &
  33 &
  310 &
  0 &
  14 &
  37 &
  310 &
  0 &
  22 &
  25 &
  310 &
  0 &
  \textbf{11} \\
 &
  2 &
  92 &
  315 &
  0 &
  28 &
  39 &
  315 &
  0 &
  17 &
  46 &
  315 &
  0 &
  23 &
  33 &
  315 &
  0 &
  \textbf{13} \\
 &
  3 &
  88 &
  314 &
  \textless{}$1e^{-6}$ &
  27 &
  41 &
  314 &
  0 &
  16 &
  43 &
  314 &
  0 &
  23 &
  37 &
  314 &
  0 &
  \textbf{14} \\ \hline
\multirow{3}{*}{MaaS-50-5} &
  1 &
  312 &
  127 &
  0 &
  128 &
  137 &
  127 &
  0 &
  \textbf{65} &
  155 &
  127 &
  \textless{}$1e^{-6}$ &
  86 &
  135 &
  127 &
  \textless{}$1e^{-6}$ &
  67 \\
 &
  2 &
  340 &
  243 &
  0 &
  136 &
  136 &
  243 &
  0 &
  \textbf{70} &
  216 &
  243 &
  \textless{}$1e^{-6}$ &
  115 &
  145 &
  243 &
  \textless{}$1e^{-6}$ &
  72 \\
 &
  3 &
  236 &
  578 &
  \textless{}$1e^{-6}$ &
  108 &
  58 &
  578 &
  0 &
  \textbf{30} &
  235 &
  578 &
  \textless{}$1e^{-6}$ &
  106 &
  88 &
  578 &
  \textless{}$1e^{-6}$ &
  40 \\\hline
\multirow{3}{*}{MaaS-70-6} &
  1 &
  11456 &
  499 &
  0 &
  6823 &
  262 &
  499 &
  0 &
  \textbf{202} &
  3890 &
  499 &
  0 &
  4569 &
  268 &
  499 &
  0 &
  221 \\
 &
  2 &
  11668 &
  743 &
  0 &
  7021 &
  78 &
  743 &
  0 &
  \textbf{59} &
  3011 &
  743 &
  0 &
  3895 &
  89 &
  743 &
  0 &
  68 \\
 &
  3 &
  11255 &
  603 &
  0 &
  6802 &
  133 &
  603 &
  \textless{}$1e^{-6}$ &
  \textbf{108} &
  2400 &
  603 &
  0 &
  2890 &
  160 &
  603 &
  0 &
  146 \\\hline
\multirow{3}{*}{MaaS-90-6} &
  1 &
  18128 &
  -$\infty$ &
  100\% &
  10800 &
  210 &
  832 &
  0 &
  \textbf{222} &
  595 &
  832 &
  0 &
  10800 &
  595 &
  832 &
  0 &
  611 \\
 &
  2 &
  16887 &
  -$\infty$ &
  100\% &
  10800 &
  217 &
  890 &
  0 &
  \textbf{242} &
  2281 &
  890 &
  0 &
  10800 &
  2281 &
  890 &
  0 &
  2172 \\
 &
  3 &
  17473 &
  -$\infty$ &
  100\% &
  10800 &
  294 &
  999 &
  \textless{}$1e^{-6}$ &
  \textbf{311} &
  636 &
  999 &
  0 &
  10800 &
  636 &
  999 &
  \textless{}$1e^{-6}$ &
  682 \\\hline
\multirow{3}{*}{MaaS-110-7} &
  1 &
  15436 &
  -$\infty$ &
  100\% &
  10800 &
  234 &
  716 &
  0 &
  \textbf{226} &
  3777 &
  -$\infty$ &
  100\% &
  10800 &
  1077 &
  716 &
  \textless{}$1e^{-6}$ &
  1314 \\
 &
  2 &
  15992 &
  -$\infty$ &
  100\% &
  10800 &
  329 &
  1212 &
  0 &
  \textbf{379} &
  3677 &
  -$\infty$ &
  100\% &
  10800 &
  967 &
  1212 &
  \textless{}$1e^{-6}$ &
  1228 \\
 &
  3 &
  15425 &
  -$\infty$ &
  100\% &
  10800 &
  276 &
  1102 &
  0 &
  \textbf{298} &
  3985 &
  -$\infty$ &
  100\% &
  10800 &
  1245 &
  1102 &
  \textless{}$1e^{-6}$ &
  1442 \\\hline
\multirow{3}{*}{MaaS-130-7} &
  1 &
  13582 &
  -$\infty$ &
  100\% &
  10800 &
  1030 &
  1034 &
  \textless{}$1e^{-6}$ &
  \textbf{1704} &
  3523 &
  -$\infty$ &
  100\% &
  10800 &
  1432 &
  1034 &
  \textless{}$1e^{-6}$ &
  2255 \\
 &
  2 &
  13559 &
  -$\infty$ &
  100\% &
  10800 &
  1342 &
  989 &
  0 &
  \textbf{2247} &
  3913 &
  -$\infty$ &
  100\% &
  10800 &
  1556 &
  989 &
  0 &
  2378 \\
 &
  3 &
  13902 &
  -$\infty$ &
  100\% &
  10800 &
  1021 &
  1120 &
  0 &
  \textbf{1689} &
  3803 &
  -$\infty$ &
  100\% &
  10800 &
  1641 &
  1120 &
  0 &
  2417 \\\hline
\multirow{3}{*}{MaaS-170-9} &
  1 &
  5776 &
  -$\infty$ &
  100\% &
  10800 &
  1516 &
  1671 &
  \textless{}$1e^{-6}$ &
  \textbf{3063} &
  3026 &
  -$\infty$ &
  100\% &
  10800 &
  3869 &
  1671 &
  \textless{}$1e^{-6}$ &
  7118 \\
 &
  2 &
  5694 &
  -$\infty$ &
  100\% &
  10800 &
  858 &
  1744 &
  0 &
  \textbf{1688} &
  2937 &
  -$\infty$ &
  100\% &
  10800 &
  2614 &
  1744 &
  0 &
  4174 \\
 &
  3 &
  5761 &
  -$\infty$ &
  100\% &
  10800 &
  602 &
  1648 &
  0 &
  \textbf{1224} &
  3141 &
  -$\infty$ &
  100\% &
  10800 &
  834 &
  1648 &
  0 &
  1516 \\\hline
\multirow{3}{*}{MaaS-190-10} &
  1 &
  4786 &
  -$\infty$ &
  100\% &
  10800 &
  438 &
  1784 &
  0 &
  \textbf{1028} &
  2240 &
  -$\infty$ &
  100\% &
  10800 &
  3917 &
  1784 &
  \textless{}$1e^{-6}$ &
  9800 \\
 &
  2 &
  4955 &
  -$\infty$ &
  100\% &
  10800 &
  519 &
  1702 &
  0 &
  \textbf{1141} &
  2421 &
  -$\infty$ &
  100\% &
  10800 &
  2233 &
  1702 &
  \textless{}$1e^{-6}$ &
  3882 \\
 &
  3 &
  4862 &
  -$\infty$ &
  100\% &
  10800 &
  560 &
  1787 &
  0 &
  \textbf{1325} &
  2322 &
  -$\infty$ &
  100\% &
  10800 &
  1641 &
  1787 &
  \textless{}$1e^{-6}$ &
  3186 \\\hline
\multirow{3}{*}{MaaS-200-10} &
  1 &
  4632 &
  -$\infty$ &
  100\% &
  10800 &
  1351 &
  1660 &
  0 &
  \textbf{3155} &
  2230 &
  -$\infty$ &
  100\% &
  10800 &
  2743 &
  1660 &
  \textless{}$1e^{-6}$ &
  9116 \\
 &
  2 &
  4616 &
  -$\infty$ &
  100\% &
  10800 &
  2894 &
  1643 &
  0 &
  \textbf{7572} &
  2688 &
  -$\infty$ &
  100\% &
  10800 &
 3656 &
  1643 &
  \textless{}$1e^{-6}$ &
  10262 \\
 &
  3 &
  4457 &
  -$\infty$ &
  100\% &
  10800 &
  877 &
  1847 &
  0 &
  \textbf{2212} &
  2867 &
  -$\infty$ &
  100\% &
  10800 &
  3923 &
  1847 &
  \textless{}$1e^{-6}$ &
  10703\\ \bottomrule
\end{tabular}
{Numbers in bold denotes the smallest CPU runtime for each instance}
\end{table}
\setlength{\abovecaptionskip}{0pt}
\setlength{\belowcaptionskip}{0pt}
\begin{table}[ht!]
\centering
\scriptsize
\caption{Performance of the proposed B\&B algorithms benchmarked with \Bard\ for \SLMFGQ\space(MIQBP)
\label{Qtable}}
{{\begin{tabular}{@{}cccccccccccccccccc@{}}
\toprule
Instances &
   &
  \multicolumn{4}{c}{Bard\&Moore} &
  \multicolumn{4}{c}{SD-BP} &
  \multicolumn{4}{c}{SD-Diffob} &
  \multicolumn{4}{c}{SD-Wi} \\ \cmidrule(l){3-6} \cmidrule(l){7-10} \cmidrule(l){11-14} \cmidrule(l){15-18}
 &
   &
  \textit{k} &
  \textit{LB} &
  \textit{Gap(\%) } &
  \textit{$T (s)$ } &
  \textit{k} &
  \textit{LB} &
  \textit{Gap(\%)} &
  \textit{$T (s)$ } &
  \textit{k} &
  \textit{LB} &
  \textit{Gap(\%)} &
  \textit{$T (s)$ } &
  \textit{k} &
  \textit{LB} &
  \textit{Gap(\%)} &
  \textit{$T (s)$ } \\\cmidrule(l){3-6} \cmidrule(l){7-10} \cmidrule(l){11-14} \cmidrule(l){15-18}
\multirow{3}{*}{MaaS-10-5} &
  1 &
  52 &
  19 &
  0 &
  9 &
  30 &
  19 &
  0 &
  5 &
  26 &
  19 &
  0 &
  7 &
  26 &
  19 &
  0 &
  \textbf{6} \\
 &
  2 &
  52 &
  -37 &
  0 &
  9 &
  29 &
  -37 &
  0 &
  5 &
  25 &
  -37 &
  0 &
  7 &
  25 &
  -37 &
  0 &
  \textbf{6} \\
 &
  3 &
  56 &
  -42 &
  0 &
  9 &
  33 &
  -42 &
  0 &
  6 &
  27 &
  -42 &
  0 &
  7 &
  27 &
  -42 &
  0 &
  \textbf{6} \\\hline
\multirow{3}{*}{MaaS-30-5} &
  1 &
  44 &
  310 &
  0 &
  10 &
  19 &
  310 &
  0 &
  \textbf{7} &
  17 &
  310 &
  0 &
  \textbf{7} &
  16 &
  310 &
  0 &
  9 \\
 &
  2 &
  392 &
  301 &
  0 &
  95 &
  38 &
  301 &
  0 &
  \textbf{11} &
  186 &
  301 &
  0 &
  67 &
  35 &
  301 &
  0 &
  15 \\
 &
  3 &
  216 &
  280 &
  \textless{}$1e^{-6}$ &
  51 &
  98 &
  280 &
  \textless{}$1e^{-6}$ &
  30 &
  96 &
  280 &
  \textless{}$1e^{-6}$ &
  39 &
  70 &
  280 &
  \textless{}$1e^{-6}$ &
  \textbf{28} \\\hline
\multirow{3}{*}{MaaS-50-5} &
  1 &
  300 &
  127 &
  0 &
  112 &
  72 &
  127 &
  0 &
  \textbf{37} &
  209 &
  127 &
  0 &
  94 &
  128 &
  127 &
  0 &
  63 \\
 &
  2 &
  120 &
  226 &
  0 &
  46 &
  67 &
  226 &
  0 &
  \textbf{27} &
  81 &
  226 &
  0 &
  40 &
  59 &
  226 &
  0 &
  30 \\
 &
  3 &
  76 &
  578 &
  \textless{}$1e^{-6}$ &
  29 &
  38 &
  578 &
  \textless{}$1e^{-6}$ &
  \textbf{15} &
  38 &
  578 &
  \textless{}$1e^{-6}$ &
  20 &
  38 &
  578 &
  \textless{}$1e^{-6}$ &
  20 \\\hline
\multirow{3}{*}{MaaS-70-6} &
  1 &
  21113 &
  -$\infty$ &
  100\% &
  10800 &
  122 &
  492 &
  0 &
  \textbf{96} &
  16804 &
  492 &
  0 &
  10714 &
  158 &
  492 &
  0 &
  138 \\
 &
  2 &
  22087 &
  -$\infty$ &
  100\% &
  10800 &
  178 &
  675 &
  \textless{}$1e^{-6}$ &
  \textbf{136} &
  10690 &
  675 &
  0 &
  6699 &
  305 &
  675 &
  \textless{}$1e^{-6}$ &
  296 \\
 &
  3 &
  20782 &
  -$\infty$ &
  100\% &
  10800 &
  243 &
  535 &
  \textless{}$1e^{-6}$ &
  \textbf{203} &
  12836 &
  535 &
  0 &
  10728 &
  216 &
  535 &
  \textless{}$1e^{-6}$ &
  208 \\\hline
\multirow{3}{*}{MaaS-90-6} &
  1 &
  16609 &
  -$\infty$ &
  100\% &
  10800 &
  104 &
  760 &
  0 &
  \textbf{104} &
  10031 &
  -$\infty$ &
  100\% &
  10800 &
  153 &
  760 &
  0 &
  144 \\
 &
  2 &
  18945 &
  -$\infty$ &
  100\% &
  10800 &
  126 &
  827 &
  0 &
  \textbf{166} &
  6390 &
  827 &
  \textless{}$1e^{-6}$ &
  6120 &
  234 &
  827 &
  0 &
  217 \\
 &
  3 &
  17824 &
  -$\infty$ &
  100\% &
  10800 &
  144 &
  917 &
  \textless{}$1e^{-6}$ &
  \textbf{160} &
  10756 &
  -$\infty$ &
  100\% &
  10800 &
  402 &
  917 &
  \textless{}$1e^{-6}$ &
  396 \\\hline
\multirow{3}{*}{MaaS-110-7} &
  1 &
  15194 &
  -$\infty$ &
  100\% &
  10800 &
  147 &
  668 &
  0 &
  \textbf{169} &
  8665 &
  -$\infty$ &
  100\% &
  10800 &
  715 &
  668 &
  0 &
  869 \\
 &
  2 &
  14945 &
  -$\infty$ &
  100\% &
  10800 &
  497 &
  1107 &
  0 &
  \textbf{694} &
  9237 &
  -$\infty$ &
  100\% &
  10800 &
  866 &
  1107 &
  0 &
  1189 \\
 &
  3 &
  14952 &
  -$\infty$ &
  100\% &
  10800 &
  128 &
  1055 &
  0 &
  \textbf{173} &
  9114 &
  -$\infty$ &
  100\% &
  10800 &
  928 &
  1055 &
  0 &
  1099 \\\hline
\multirow{3}{*}{MaaS-130-7} &
  1 &
  9219 &
  -$\infty$ &
  100\% &
  10800 &
  626 &
  945 &
  0 &
  \textbf{977} &
  5610 &
  -$\infty$ &
  100\% &
  10800 &
  1762 &
  945 &
  0 &
  2321 \\
 &
  2 &
  12685 &
  -$\infty$ &
  100\% &
  10800 &
  204 &
  949 &
  0 &
  \textbf{284} &
  8036 &
  -$\infty$ &
  100\% &
  10800 &
  1670 &
  949 &
  0 &
  2284 \\
 &
  3 &
  13396 &
  -$\infty$ &
  100\% &
  10800 &
  223 &
  1076 &
  0 &
  \textbf{336} &
  8154 &
  -$\infty$ &
  100\% &
  10800 &
  2045 &
  1076 &
  0 &
  2732 \\\hline
\multirow{3}{*}{MaaS-170-9} &
  1 &
  5343 &
  -$\infty$ &
  100\% &
  10800 &
  140 &
  1643 &
  0 &
  \textbf{267} &
  2712 &
  -$\infty$ &
  100\% &
  10800 &
  661 &
  1643 &
  0 &
  1340 \\
 &
  2 &
  5325 &
  -$\infty$ &
  100\% &
  10800 &
  190 &
  1694 &
  0 &
  \textbf{323} &
  2702 &
  -$\infty$ &
  100\% &
  10800 &
 888 &
  1694 &
  0 &
  1877 \\
 &
  3 &
  5564 &
  -$\infty$ &
  100\% &
  10800 &
  200 &
  1605 &
  0 &
  \textbf{411} &
  2709 &
  -$\infty$ &
  100\% &
  10800 &
  749 &
  1605 &
  0 &
  1481 \\\hline
\multirow{3}{*}{MaaS-190-10} &
  1 &
  4728 &
  -$\infty$ &
  100\% &
  10800 &
  364 &
  1724 &
  \textless{}$1e^{-6}$ &
  \textbf{855} &
  2093 &
  -$\infty$ &
  100\% &
  10800 &
  1036 &
  1724 &
  \textless{}$1e^{-6}$ &
  2363 \\
 &
  2 &
  5104 &
  -$\infty$ &
  100\% &
  10800 &
  231 &
  1673 &
  \textless{}$1e^{-6}$ &
  \textbf{491} &
  2143 &
  -$\infty$ &
  100\% &
  10800 &
  1374 &
  1673 &
  \textless{}$1e^{-6}$ &
  2965 \\
 &
  3 &
  4889 &
  -$\infty$ &
  100\% &
  10800 &
  270 &
  1730 &
  \textless{}$1e^{-6}$ &
  \textbf{626} &
  2246 &
  -$\infty$ &
  100\% &
  10800 &
  1297 &
  1730 &
  \textless{}$1e^{-6}$ &
  2736 \\\hline
\multirow{3}{*}{MaaS-200-10} &
  1 &
  4741 &
  -$\infty$ &
  100\% &
  10800 &
  565 &
  1621 &
  \textless{}$1e^{-6}$ &
  \textbf{1230} &
  2202 &
  -$\infty$ &
  100\% &
  10800 &
  2487 &
  1621 &
  \textless{}$1e^{-6}$ &
  5196\\
 &
  2 &
  4607 &
  -$\infty$ &
  100\% &
  10800 &
  256 &
  1584 &
  \textless{}$1e^{-6}$ &
  \textbf{561} &
  1986 &
  -$\infty$ &
  100\% &
  10800 &
  2474 &
  1584 &
  \textless{}$1e^{-6}$ &
 4783 \\
 &
  3 &
  4880 &
  -$\infty$ &
  100\% &
  10800 &
  509 &
  1814 &
  \textless{}$1e^{-6}$ &
  \textbf{1194} &
  1974 &
  -$\infty$ &
  100\% &
  10800 &
  2480 &
  1814 &
  \textless{}$1e^{-6}$ &
  4962 \\ \bottomrule
\end{tabular}}}	
{Numbers in bold denotes the smallest CPU runtime for each instance}
\end{table}
Globally, we find that Algorithm \BP\ is the best-performing method to solve both \SLMFGL\xspace and \SLMFGQ\xspace in terms of runtime. For \Bard\xspace and \Diff, as the number of followers increases, the average value of $k$ tends to first increase and then decrease in the following 13 groups where the time limit (10800 s) has been reached and the CPU runtime of solving each sub-problem increases; moreover, their runtime  exponentially grows with the size of the problem. As the size of the problem increases, the CPU runtime of \BP\xspace is significantly smaller than its counterparts, i.e., \BP\ can save at least 2 s $\sim$ 6 s on the instance `MaaS-10-5', 10633 s $\sim$ 10749 s on the instance `MaaS-80-6' and at least 3228 s $\sim$ 10239 s on the instance `MaaS-200-10', compared with the benchmark (\Bard). We further examine the behavior of the proposed SD-based algorithm by comparing the average CPU runtime with each of the three branching rules (\BP,\Weight,\Diff) over all 20 instances of each group of instances, for 10 groups. The results for both \SLMFGL\space and \SLMFGQ\ are depicted in Figure \ref{comparebranch}. 
The above numerical results along with those in Online Appendix \ref{Numerical results} demonstrate the following insights.

\begin{observation}
Compared with \Bard\xspace which exploits the explicit disjunctive structure of the complementarity conditions derived from a single violated complementarity condition, e.g., each of \MPECL\ or \MPECQ\ contains eight complementarity slackness conditions, the proposed SD-based B\&B algorithms only branch once per follower problem, leading to a huge computational improvements, i.e., the runtime of \BP\space is over 10 $\sim$ 100 times faster than that of \Bard.
\end{observation}
\begin{observation}
When the the scale of the instances is small, the CPU runtime of \Weight\xspace is slightly smaller than its counterparts. when the scale of the instances is large, \BP\space considerably outperforms the others for both SLMFGs.
\end{observation}
\begin{observation}
Since $\textit{LB}$ is only updated if all follower problems are optimal. Hence, observe that even though in the quadratic case where \SDQ\ is a relaxation of \SLMFGQ\ using McCormick envelop, the value of leader's objective value ($\textit{LB}$) is same for four solution methods in each instance, confirming the effectiveness of proposed algorithm.
\end{observation}
\begin{observation}
The profits of the MaaS regulator obtained in SLMFG with network effects (\SLMFGQ) is not greater than that in SLMFG without network effects (\SLMFGL). Moreover, these two formulations differ in the computational efforts required to solve the corresponding SLMFGs. Specifically, the average CPU runtime of \SLMFGQ\ is faster than that of \SLMFGL.
\end{observation}
\subsection{Sensitivity analysis}
\label{sa}
\begin{figure}[ht!]
\centering
\vspace{0.05in}
\caption{Sensitivity analysis on sell-bidding price range ratio $\bm{\beta_{max}/\beta_{min}}$}
\subfloat[\SLMFGL]{\includegraphics[width=0.4\textwidth]{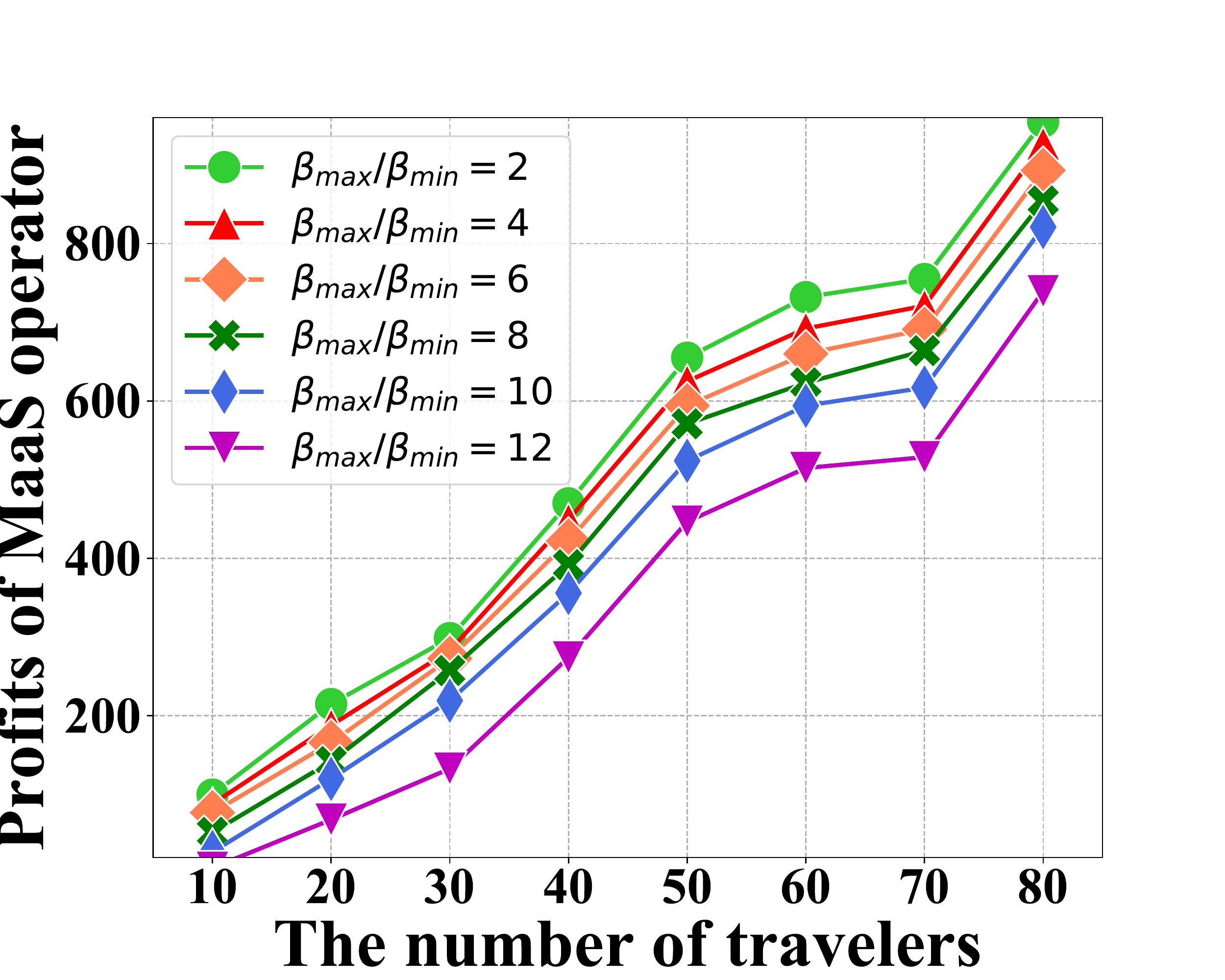}\label{Qb1}}
\subfloat[\SLMFGQ]{\includegraphics[width=0.4\textwidth]{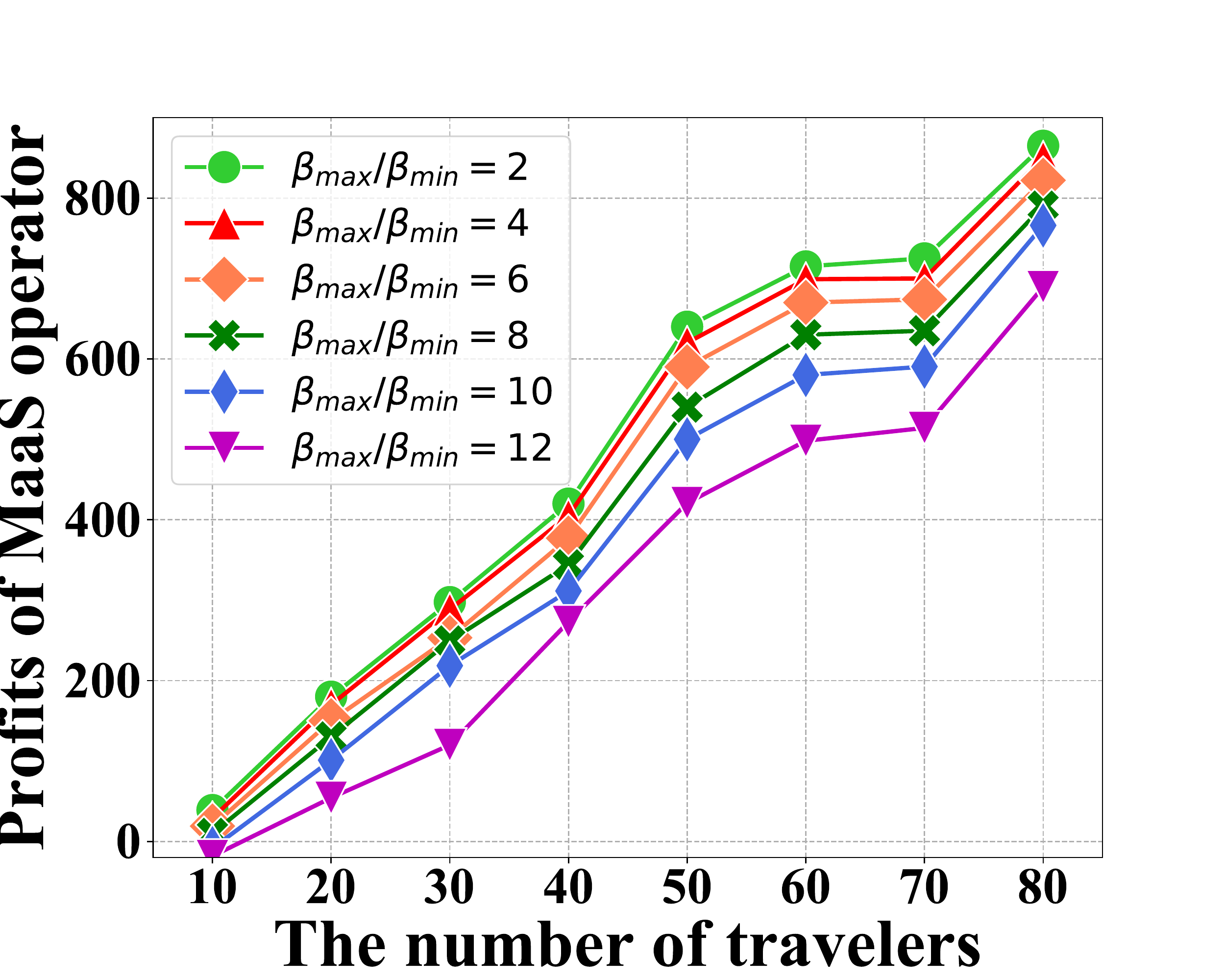}\label{Qb2}}
\label{Qb}
\end{figure}

\begin{figure}[t!]
\centering
\vspace{0.05in}
\caption{Sensitivity analysis on the lower and upper bound of $\Delta$ for \SLMFGQ.}
\subfloat[Lower-bound of $\Delta$ ($\underline{C}$) ]{\includegraphics[width=0.4\textwidth]{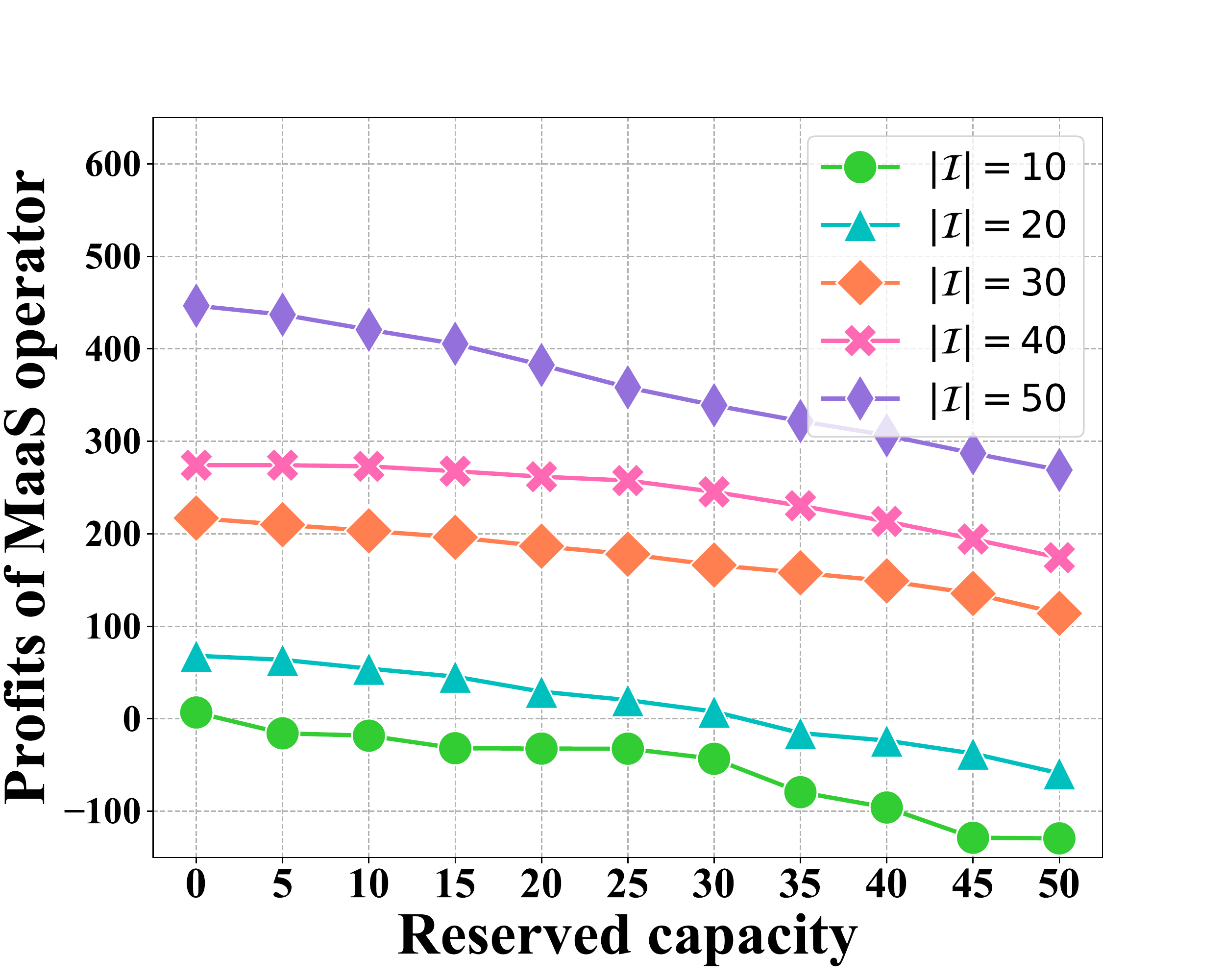}\label{LC1}}
\subfloat[Upper-bound of $\Delta$ ($\bar{C}$) ]{\includegraphics[width=0.4\textwidth]{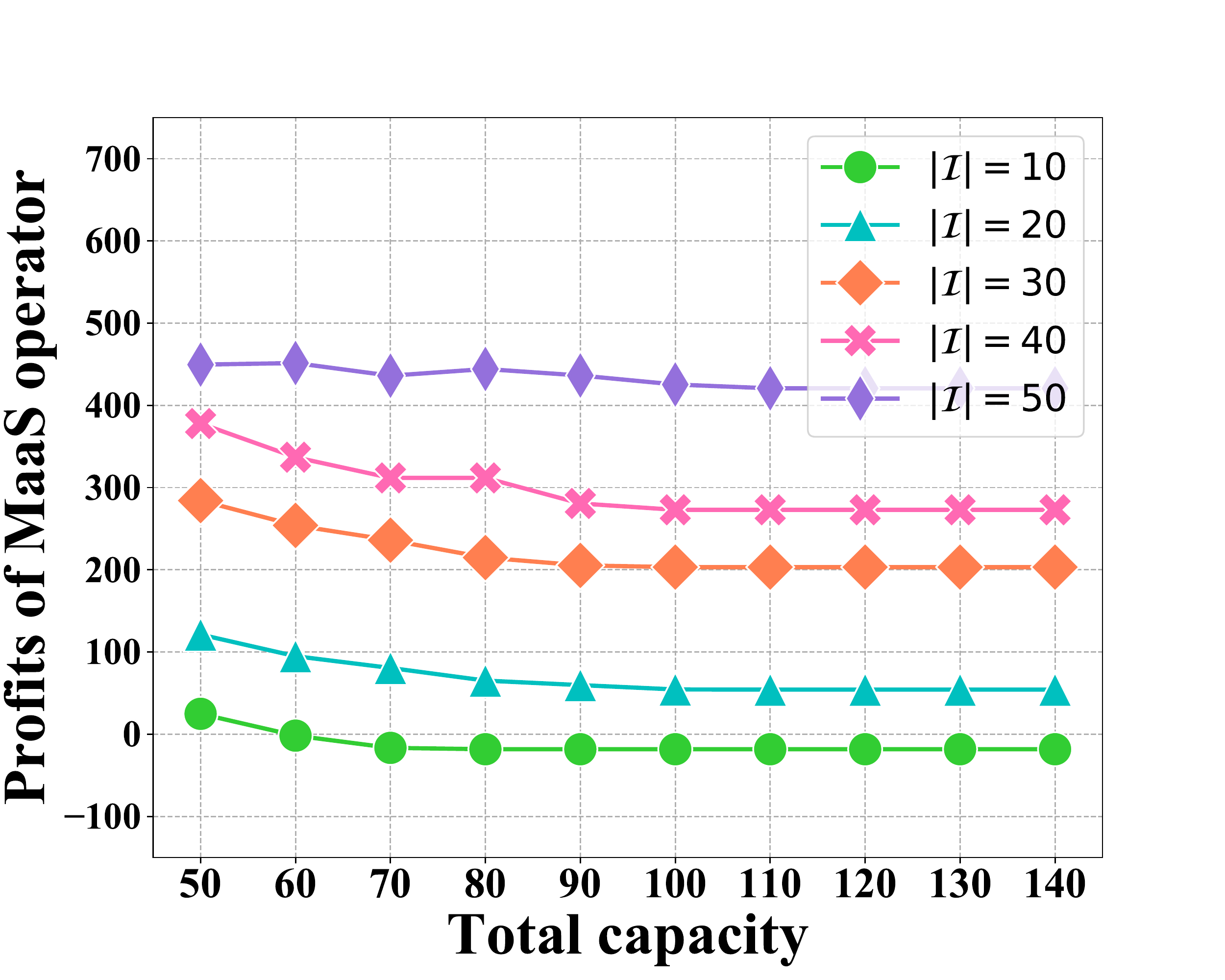}\label{UC1}}
\label{bound}
\label{Cap}
\end{figure}

In this subsection, we conduct a sensitivity analysis on the parameters of the proposed SLMFG formulations. Algorithm \SDBB\space with branching rule \BP\ is used to solve all instances. We define the purchase-bidding price range ratio as the ratio of travelers' maximum purchase-bidding price to the minimum purchase-bidding price $b_{max}/b_{min}$, and the sell-bidding price range ratio as the ratio of TSPs' maximum sell-bidding price to the minimum sell-bidding price $\beta_{max}/\beta_{min}$. We conduct a sensitivity analysis on the bounds of the supply-demand gap, $[\underline{C},\bar{C}]$, where the upper bound $\bar{C}$ is the total capacity of mobility resources supplied by all TSPs and the lower bound $\underline{C}$ is the reserved capacity of MaaS regulator. We examine the behavior of the linear and quadratic cases, \SLMFGL\ and \SLMFGQ, for varying  purchase-bid range ratio and sell-bid range ratio in \Fig\ref{Lb} and \Fig\ref{Qb}. We examine the behavior of the \SLMFGQ, for varying  reserved capacity and total capacity in \Fig\ref{Cap}. We summarize the following
insights which can
be used as practical managerial guidelines:
\begin{figure}[ht!]
\centering
\vspace{0.05in}
\caption{Sensitivity analysis on purchase-bidding price range ratio $\bm{b_{max}/b_{min}}$}
\subfloat[\SLMFGL ]{\includegraphics[width=0.4\textwidth]{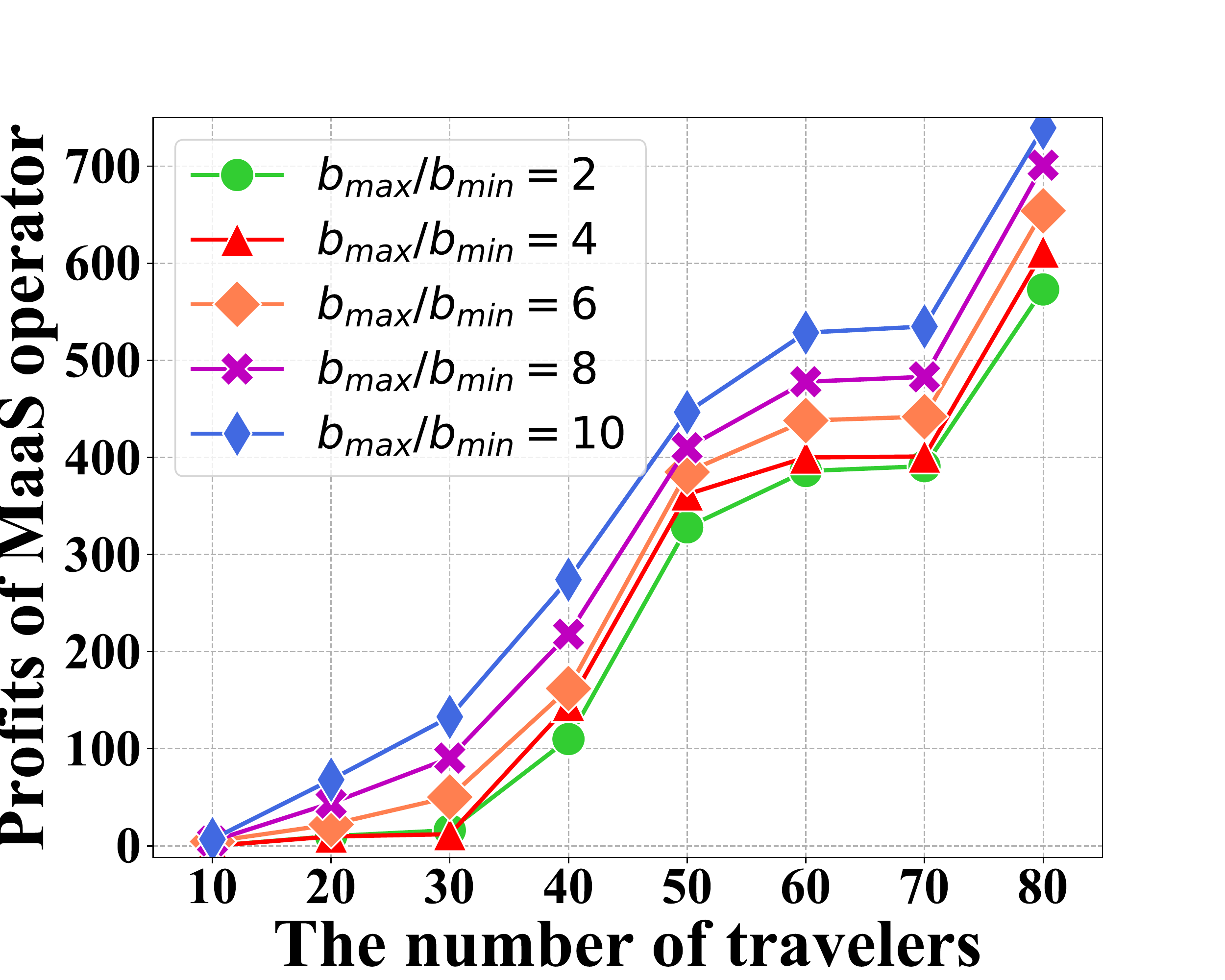}\label{Lb1}}
\subfloat[\SLMFGQ]{\includegraphics[width=0.4\textwidth]{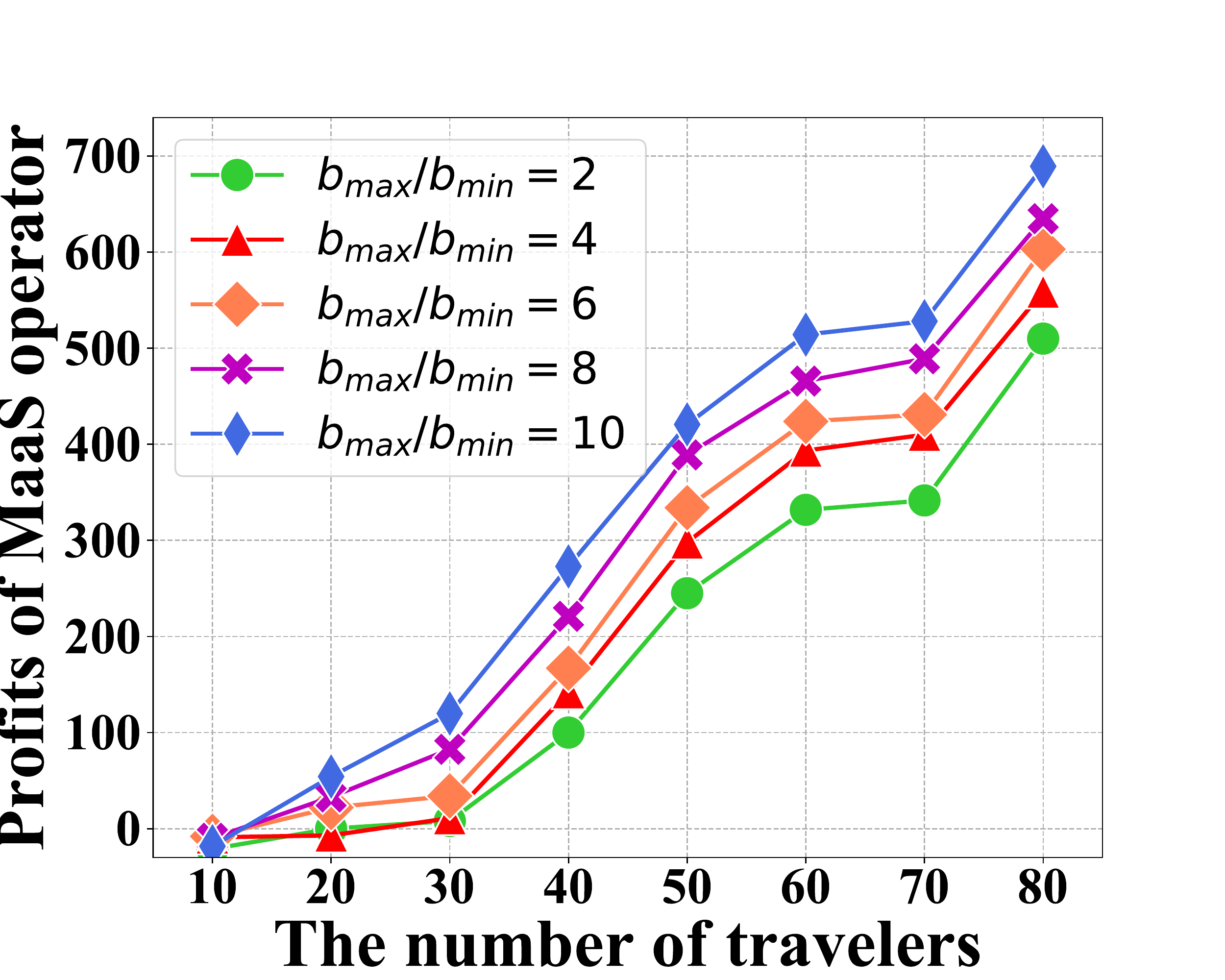}\label{Lb2}}
\label{Lb}
\end{figure}
\newpage
\begin{observation}
In SLMFG with/without network effects, an increase of the purchase-bidding price range ratio ($b_{\max}/b_{\min}$) leads to an increase in MaaS regulator's profits  for a fixed number of travelers, while an increase of the sell-bidding price range ratio ($\beta_{\max}/\beta_{\min}$) leads to a decrease in the MaaS regulator's profits  for a fixed number of TSPs. The profits of the MaaS regulator increase with the number of travelers (TSPs) under a fixed given purchase (sell)-bidding price range ratio.
\end{observation}

\begin{observation}
In SLMFG with network effects, MaaS regulator's profits will considerably decrease as $\underline{C}$ increases, and first slightly decrease and then keep unchanged as $\bar{C}$ increases. In SLMFG without network effects, the variation of $\underline{C}$ and $\bar{C}$ has no impacts on MaaS regulator's profits.
 \end{observation}
\section{Conclusion and remarks}
\label{Conclusion}
In this study, we proposed a novel modeling and optimization framework for the regulation of two-sided MaaS markets. We cast this problem as a SLMFG where the leader is the MaaS regulator and two groups of followers are travelers and TSPs. We propose a NYOP-auction mechanism which allows travelers submit purchase-bids to fulfill their travel demand via MaaS platform and TSPs submit sell-bids to supply mobility resources for the MaaS platform simultaneously. We capture the network effects through the supply-demand in the two-sided MaaS markets, and consider SLMFG with and without network effects corresponding to MILBP and MIQBP problem, respectively. We provide constraint qualifications for MPEC reformulations of the both SLMFGs and prove the equivalence between the MPEC reformulations and their original problems, which lays a basis to develop solution algorithms based on the single-level reformulation of the SLMFGs. Further, we propose an exact solution method, SD-based B\&B algorithm named \SDBB\xspace, which branches on accept/reject binary variables, and customize three rules to select branching variables. Our computational results show the performance of the proposed \SDBB\space algorithms based on SD reformulation relative to a benchmark B\&B algorithm based on MPEC reformulation and conclude that \BP\space considerably outperforms the others for both SLMFGs, i.e., the computational speed of \BP\space is over 10 $\sim$ 100 times faster than that of the benchmark when the scale of the instances is large. Our computational experiments uncover managerial insights into how the proposed SLMFGs behave in relation to the parameters, show that the MaaS regulator's profits of the SLMFG with network effects are not
greater than that of SLMFG without network effects, and indicate that both travelers' travel costs and TSPs' profits decrease with the increase of the supply-demand gap.

\HX{MaaS is a framework for delivering a portfolio of multi-modal mobility services that places user experience at the centre of the offer. Thus, instead of considering users' mode choice across a set of travel modes, the proposed MaaS system unifies the mobility services provided by different travel modes as \textit{mobility resources} and customize the optimal feasible \textit{MaaS bundles} according to their WTP and experience-relevant requirements. Therefore, the NYOP-auction is perfectly adapted to the proposed user-centric MaaS framework.} The results of this research provide managerial insights on how MaaS regulators make operational strategies (e.g., pricing and \textit{MaaS bundles}) to attract price-sensitive travelers and TSPs with heterogeneous preference, WTP or WTS. Our findings also highlight how different stakeholders (MaaS regulator, travelers and TSPs) interact with each other to maximize their own benefits, notably the impact of cross-network effects on the profits of the MaaS regulator. \HX{This study suggests that MaaS platforms should consider cross network effects between travelers and TSPs to design their operational strategies.}

This study can be extended in several research directions. The proposed SLMFG can be extended to multi-leader-follower problems to model the competition among different MaaS regulators. \HX{The proposed \SDBB\ algorithm is applicable for any MIBLP/MIBQP with linear/quadratic follower problems}. Future research is needed to explore stochastic formulations to capture the uncertainty in supply and demand, as well as the potential of time-varying pricing strategies.

\bibliographystyle{pomsref} 

 \let\oldbibliography\thebibliography
 \renewcommand{\thebibliography}[1]{%
 	\oldbibliography{#1}%
 	\baselineskip14pt 
 	\setlength{\itemsep}{10pt}
 }
\bibliography{ref1} 



\ECSwitch 

\ECHead{Online Appendix }
\vspace{12 pt}
\noindent \textbf{for ``Single-leader multi-follower games for the regulation
of two-sided Mobility-as-a-Service markets''}
\vspace{12 pt}
\section{Mathematical notations}
\label{EC1}
\setlength{\abovecaptionskip}{0pt}
\setlength{\belowcaptionskip}{0pt}
\begin{table}[ht!]
\label{A}
\centering
\footnotesize
\caption{Parameters of the mathematical models}
{{\begin{tabular}{ll}
\toprule
\multicolumn{2}{l}{\textbf{Parameters}}\\
 \hline
$b_{i}$ & Traveler $i$'s purchase-bidding price\\
$B_{i}$ & Traveler $i$'s travel expenditure budget\\
$\bar{B}_{mn}$ & TSP $mn$'s operating cost budget\\
$C_{mn}$& The capacity of the mobility resources provided by TSP $mn$ \\
$\bar{C}$& Capacity of the mobility resources in the MaaS system provided by all TSPs, $C=\sum_{m \in \mathcal{M}}\sum_{n \in \mathcal{N}_{m}}C_{mn}$\\
$D_{i}$ & Traveler $i$'s requested travel distance\\
$Q_{i}$ & Traveler $i$'s requested quantity of mobility resources\\
$r_{i}$ & Traveler $i$'s unit reserve price for mobility resources\\
${R}_{i}$ & Traveler $i$'s maximum acceptable delay time during a mobility service (travel delay budget) \\
$T_{i}$ &Traveler $i$'s requested total service time\\
$v_{m}$ & Average commercial speed of travel mode $m$ \\
$\alpha_{i}$ & Traveler $i$'s unit waiting time cost\\
$\beta_{mn}$ & TSP $mn$'s unit idle time cost\\
$\rho_{mn}$ &TSP $mn$'s unit reserve price\\
$\Gamma_{i}$ &Traveler $i$'s  maximum acceptable inconvenience cost during a mobility service (inconvenience tolerance)\\
$\delta_{m}$ & Inconvenience cost per unit of time for travel mode $m$\\
$\eta_{mn}$ & TSP $mn$'s unit idle time cost\\
$\gamma_{mn}$ & TSP $mn$'s unit operating cost\\
\hline
\multicolumn{2}{l}{\textbf{Decision variables}}\\
 \hline
$p$&Real varible denoting the \HX{unit threshold price} set for travelers using the MaaS platform\\
$q$&Real varible denoting the \HX{unit threshold price} set for TSPs participating in the MaaS platform\\
$u_{i}$& Binary variable indicating whether Traveler $i$ is accepted to register MaaS platform or not\\
$w_{mn}$& Binary variable indicating whether the MaaS regulator offers TSP $mn$ to join the MaaS platform\\
$l_{i}^{m}$& Real variable denoting the service time of travel mode $m$ in the MaaS bundle allocated to Traveler $i$\\
$\Delta$& Real variable denoting supply-demand gap\\
 $x_{i}$& Real varible representing the proportion mobility demand that $i$ chooses to use the MaaS platform\\
 $y_{mn}$& Real variable representing the proportion of mobility resources to supply the MaaS platform\\
\hline
\multicolumn{2}{l}{\textbf{Sets}}\\
 \hline
 $\mathcal{I}$& Set of travelers, $\mathcal{I}=\left\{1,2,\cdots,i,\cdots,|\mathcal{I}|\right\}$  \\
 $\mathcal{M}$& Set of travel modes, $\mathcal{M}=\left\{1,2,\cdots,m,\cdots,|\mathcal{M}|\right\}$ \\
  $\mathcal{N}_{m}$& Set of TSPs of travel mode $m$, $\mathcal{N}_{m}=\left\{1,2,\cdots,n,\cdots,|\mathcal{N}_{m}|\right\}, \forall m\in\mathcal{M}$\\
\bottomrule
\end{tabular}}}	
\end{table}

\section{Equivalent formulations}
\subsection{Proof of Proposition 2}
\label{ECProp2}
\proof{\textsc{Proof:}} 
\indent In \SLMFGL, since the objective function of travelers/TSPs has a separable structure and can be rewritten as:
\begin{equation}
G_{i}(x_{i},\Delta)=G_{i}^1(x_{i})+G_{i}^2(\Delta)
\quad\mbox{ and }\quad
H_{mn}(y_{mn},\Delta)=H_{mn}^1(y_{mn})+H_{mn}^2(\Delta),
\end{equation}
where 
\begin{subequations}
\begin{align}
&G_{i}^1(x_{i})=b_{i}Q_{i}x_{i}+r_{i}Q_{i}(1-x_{i}) -\frac{\alpha_{i}\zeta Q_{i}}{\bar{C}}x_{i}, \\
&G_{i}^2(\Delta)=-\alpha_{i}\left(\frac{\zeta}{\bar{C}}\Delta+\vartheta\right),\\
&H_{mn}^1(y_{mn})=\left (\beta_{mn}-\rho_{mn}\right)C_{mn}y_{mn}+(\rho_{mn}-\gamma_{mn})C_{mn}(1-y_{mn})+\frac{\eta_{mn}\kappa C_{mn}}{\bar{C}}y_{mn},\\ &H_{mn}^2(\Delta)=-\eta_{mn}\left(\frac{\kappa}{\bar{C}}\Delta+\xi\right).
\end{align}
\end{subequations}

\indent As an immediate consequence, the solution set of the follower's problem will not depend on the leader variable $\Delta$, we propose an equivalent formulation for \SLMFGL\ and state Proposition 1. \Halmos
\endproof
\subsection{Proof of Proposition 3}
\label{ECProp3}
\proof{\textsc{Proof:}} 
For any feasible value of the leader's variables $u_i$ and $w_{mn}$, the constraints set of the travelers and TSPs are nonempty, compact and their objective functions are continuous in their respective variables ($x_i$ and $y_{mn}$), therefore these follower problems admit at least a solution. Moreover, since the quadratic perceived waiting/idle functions ($\Psi_{i}^{q}(x_{i},\Delta)$ /($\Phi_{mn}^{q}(y_{mn},\Delta)$) are strictly convex, the objective functions of the travelers/TSPs are strictly convex/concave. Therefore, by classical uniqueness arguments, for any feasible value of the leader's variables $u_{i}$ and $w_{mn}$, \DR{these follower problems admit a unique solution which can be denoted by $x_i(u_{i},\Delta)$ for Traveler $i$'s problem and by $y_{mn}(w_{mn},\Delta)$ for the TSP $mn$'s problem}. Due to the uniqueness of the followers' solutions, the optimistic approach used in Model \ref{SLMFG} is no longer needed and \SLMFGQ\ is equivalent to SLMFG with network effects described in Proposition \ref{prop:SLMFG_quadra}.
\Halmos
\endproof
\section{MPEC formulations}
\subsection{MPEC-L}
\label{KKT-L}
\begin{subequations}
\begin{align}
&0\leq(B_{i}- b_{i}Q_{i}x_{i})\perp\lambda_{1}^{i}\geq 0 ,&&\forall i \in \mathcal{I},\tag{KL.1}\label{29a}\\
&0\leq\left(u_{i}-x_{i}\right)\perp\lambda_{2}^{i}\geq 0,&&\forall i \in \mathcal{I},\tag{KL.2}\label{29b} \\
&0\leq x_{i}\perp\lambda_{3}^{i}\geq 0,&&\forall i \in \mathcal{I},\tag{KL.3}\label{29c}\\
&0\leq \left[(b_{i}-r_{i})Q_{i}-\alpha_{i}\frac{\zeta}{\bar{C}}Q_{i}-\lambda_{1}^{i}b_{i}Q_{i}-\lambda_{2}^{i} +\lambda_{3}^{i}\right]\perp x_{i}\geq 0,&&\forall i \in \mathcal{I}, \tag{KL.4} \label{29d}  
\end{align}
\end{subequations}   

\noindent where $\lambda_{1}^{i}$, $\lambda_{2}^{i}$ and $\lambda_{3}^{i}$ denote the Lagrangian multipliers associated with the constraints \eqref{BL1b},\eqref{BL1c} and \eqref{BL1d}, respectively. The operator `$a\perp b$' denotes the complementarity
condition $a\cdot b = 0$. 

The KKT conditions of TSP $mn$'s problem in  \SLMFGQ\space are given in \eqref{31a}-\eqref{31d}.
\begin{subequations}
\begin{align}
&0\leq\left[\bar{B}_{mn}-C_{mn}\gamma_{mn}y_{mn}\right]\perp\mu_{1}^{mn}\geq 0,&&\forall m \in  \mathcal{M},n \in \mathcal{N}_{m},\tag{KL.5}\label{31a}\\
&0\leq (w_{mn}-y_{mn})\perp\mu_{2}^{mn}\geq 0, &&\forall m \in  \mathcal{M},n \in \mathcal{N}_{m},\tag{KL.6}\label{31b}\\
&0\leq y_{mn}\perp\mu_{3}^{mn}\geq 0,&&\forall m \in  \mathcal{M},n \in \mathcal{N}_{m},\tag{KL.7}\label{31c} \\
&0\leq\left[(\rho_{mn}-\beta_{mn})C_{mn}-\eta_{mn}\frac{\kappa}{\bar{C}}C_{mn}-\mu_{1}^{mn}C_{mn}\gamma_{mn}-\mu_{2}^{mn} +\mu_{3}^{mn}\right]\nonumber\\
&\perp y_{mn}\geq0,&&\forall m \in  \mathcal{M},n \in \mathcal{N}_{m},\tag{KL.8}\label{31d}
\end{align}
\end{subequations}  
where $\mu_{1}^{mn}$, $\mu_{2}^{mn}$ and $\mu_{3}^{mn}$ denote the Lagrangian multipliers associated with constraint \eqref{BL2b},\eqref{BL2c} and \eqref{BL2d}, respectively. 

\subsection{MPEC-Q}
\label{KKT-Q}
The KKT conditions associated to Traveler $i$'s problem are given in \eqref{48a}-\eqref{48d}:
\begin{subequations}
\begin{align}
&0\leq(B_{i}-b_{i}Q_{i}x_{i})\perp\lambda_{1}^{i}\geq 0,&&\forall i \in \mathcal{I},\tag{KQ.1}\label{48a}\\
&0\leq\left(u_{i}-x_{i}\right)\perp\lambda_{2}^{i}\geq 0,&&\forall i \in \mathcal{I},\tag{KQ.2}\label{48b} \\
&0\leq x_{i}\perp\lambda_{3}^{i}\geq 0,&&\forall i \in \mathcal{I},\tag{KQ.3}\label{48c}\\
&0\leq\left[\frac{2\zeta\alpha_{i}Q_{i}^{2}}{C^{2}}x_{i}-\left(\frac{2\zeta\alpha_{i}Q_{i}}{C}-\frac{2\zeta\alpha_{i}\Delta}{C^{2}}+r_{i}-b_{i}\right)Q_{i}-\lambda_{1}^{i}b_{i}Q_{i}-\lambda_{2}^{i}+\lambda_{3}^{i}\right]\perp x_{i}\geq 0,&&\forall i \in \mathcal{I},\tag{KQ.4} \label{48d} 
\end{align}
\end{subequations}    
\noindent where $\lambda_{1}^{i}$, $\lambda_{2}^{i}$ and $\lambda_{3}^{i}$ denotes the Lagrangian multipliers corresponding to the constraint \eqref{BL1b},\eqref{BL1c} and \eqref{BL1d}, respectively. 
\DA{On the other hand, $\forall m\in\mathcal{M},n\in\mathcal{N}_{m},$ the KKT conditions of the TSP $mn$'s follower- problem are given in
}\eqref{49a}-\eqref{49d}.
\begin{subequations}
\begin{align}
&0\leq\left[\bar{B}_{mn}-\beta_{mn}C_{mn}y_{mn}\right]\perp\mu_{1}^{mn}\geq 0,&&\forall m \in  \mathcal{M},n \in \mathcal{N}_{m},\tag{KQ.5}\label{49a}\\
&0\leq(w_{mn}-y_{mn})\perp\mu_{2}^{mn}\geq 0, &&\forall m \in  \mathcal{M},n \in \mathcal{N}_{m},\tag{KQ.6}\label{49b}\\
&0\leq y_{mn}\perp\mu_{3}^{mn}\geq 0,&&\forall m \in  \mathcal{M},n \in \mathcal{N}_{m},\tag{KQ.7}\label{49c} \\
&0\leq\Bigg[\frac{2\kappa\eta_{mn} C_{mn}^{2}}{C^{2}}y_{mn}-\left(\frac{2\kappa\eta_{mn}\Delta}{C^{2}}+\beta_{mn}-\rho_{mn}\right)C_{mn}\nonumber\\
&-\mu_{1}^{mn}\beta_{mn}C_{mn}-\mu_{2}^{mn}+\mu_{3}^{mn}\Bigg]\perp y_{mn}\geq 0,&&\forall m \in  \mathcal{M},n \in \mathcal{N}_{m},\tag{KQ.8}\label{49d}
\end{align}
\end{subequations}  
where $\mu_{1}^{mn}$, $\mu_{2}^{mn}$ and $\mu_{3}^{mn}$ denote the Lagrangian multipliers corresponding to the constraint \eqref{BL2b},\eqref{BL2c} and \eqref{BL2d}, respectively. 
\vspace{12 pt}
\section{Proof of Theorem 1 (constraint qualifications)}
\label{ECAnnex_C}
The proof of Theorem 1 is based on Theorems 3.3.8 and 3.3.11 of \cite{AusselSvenssonChapter}. Let us first observe that \SLMFGL\space as described in Proposition \ref{prop:SLMFG_lin} clearly satisfies the following basic properties:
\begin{itemize}
    \item [(H1)] (Follower's differentiability) $G_{i}(x_{i},\Delta)$ and $H_{mn}(y_{mn},\Delta)$ are differentiable;
    \item [(H2)] (Follower's convexity) $G_{i}(x_{i},\Delta)$ and $H_{mn}(y_{mn},\Delta)$ are convex with respect to variables $x_i$ and $y_{mn}$ respectively.
\end{itemize}

Thus item $a)$ of Theorem \ref{thm:equiv_MILBP} is a direct consequence of item i) of \cite{AusselSvenssonChapter}[Theorem 3.3.8].

\indent To prove item $b)$, let us now assume that for all traveler $i \in \mathcal{I}$, $b_i \geq p_{min}$ and for all TSP $n \in \mathcal{N}_m, m \in \mathcal{M}$, $\beta_{mn} \leq q_{max}$. First the MILBP clearly satisfies the so-called {Joint convexity assumption} defined in \cite{AusselSvenssonChapter} since all the constraint functions of the follower's problem are linear.

Second, and in order to prove the Joint Slater's qualification condition, let us first identify, in the forthcoming lemma, conditions under which a solution wherein all players of the market can participate. 

\DA{
\begin{lemma}
\label{lem:atleastone}
Let us assume that, for all Traveler $i \in \mathcal{I}$, $b_i \geq p_{min}$ and for all TSP $n \in \mathcal{N}_m, m \in \mathcal{M}$, $\beta_{mn} \leq q_{max}$. Then there exists at least one feasible point of  \SLMFGL\space such that $u_i = 1$, for all Traveler $i \in \mathcal{I}$, and $w_{mn} = 1$, for all TSP $n \in \mathcal{N}_m, m \in \mathcal{M}$.
\end{lemma}
\proof{\textsc{Proof:}} 
Observe that since $b_i \geq p_{min}$ for all Traveler $i \in \mathcal{I}$, if the leader chooses $p = p_{min}$, all travelers can participate, i.e. the couple $(p,u_i) = (p_{min},1)$ satisfies Eqs. \eqref{B1f} and \eqref{B1g}, for all $i \in \mathcal{I}$. Analogously, since $\beta_{mn} \leq q_{max}$ for all TSP $n \in \mathcal{N}_m, m \in \mathcal{M}$, then if the leader chooses $q = q_{max}$, all TSPs can participate, i.e. the couple $(q,w_{mn}) = (q_{max},1)$ satisfies Eqs. \eqref{B1h} and \eqref{B1i}, for all $n \in \mathcal{N}_m$, $m \in \mathcal{M}$. Since variables $p$ and $q$ only appear in constraints \eqref{B1f}-\eqref{B1i}, as well as in their bound constraints \eqref{B1m} and \eqref{B1n}, respectively, the choice $p = p_{min}$ and $q = q_{max}$ is always feasible. In addition, let us recall that for any feasible value of the leader's variables $u_i$, $w_{mn}$ and $\Delta$, the constraints set of the travelers and TSPs being nonempty compact and their objective function continuous in their respective variables ($x_i$ and $y_{mn}$), the followers' problem admits at least a solution, respectively $\bar{\bm{x}}$ for travelers and $\bar{\bm{y}}$ for TSP. As a consequence the vector $(p_{min}, q_{max}, l, \bm{1}, \bm{1}, \Delta, \bar{\bm{x}}, \bar{\bm{y}})$ is feasible for the \SLMFGL\space problem with $\bm{1}$ being vector $(1,\dots,1)$ and $(\bm{l}, \Delta)$ satisfying Eqs. \eqref{B1b}-\eqref{B1e}, \eqref{B1j}, \eqref{B1o} and \eqref{B1p}.
\endproof
Now taking into account Lemma \ref{lem:atleastone} and observing that, for any $i\in\mathcal{I}$, if $u_i=1$ then the constraint set of Traveler $i$'s problem is $[0,\min\{1,\frac{B_i}{b_iQ_i}\}]$ and  for any $(m,n)\in\mathcal{M}\times \mathcal{N}$, if $w_{mn}=1$ then the constraint set of the TSP$_{mn}$'s problem is $[0,\min\{1,\frac{\bar{B}_{mn}}{C_{mn}\omega_{mn}}\}]$ one immediately deduces that the {\em Joint Slater's qualification condition} needed in \cite{AusselSvenssonChapter}[Theorem 3.3.11] is satisfied.
}

\DA{It remains to show that {\em Guignard's qualification conditions for boundary opponent strategies} holds true. Let us precise that, according to \cite{AusselSvenssonChapter}[Definition 3.3.10] and taking into account that the constraint set of the Traveler $i$ (respectively the TSP $mn$) only depends on $u_i$ (respectively $w_{mn}$), a boundary opponent strategy for traveler $i$ (respectively for TSP $mn$) is a point $(p,q,\bm{l},\bm{u},\bm{w},\Delta,\bm{x}_{-i},\bm{y})$ (respectively $(p,q,\bm{l},\bm{u},\bm{w},\Delta,\bm{x},\bm{y}_{-{mn}})$) such that $u_i=0$ or $u_i=1$ (respectively $w_{mn}=0$ or $w_{mn}=1$). 
Let us recall also that Guignard's CQ is satisfied at $(p,q,\bm{l},\bm{u},\bm{w},\Delta,\bm{x},\bm{y})$ for follower $i$'s problem (respectively for follower $mn$'s problem) if 
\[T^o(x_i, (p,q,\bm{l},\bm{u},\bm{w},\Delta,\bm{x}_{-i},\bm{y}))=L^o(x_i, (p,q,\bm{l},\bm{u},\bm{w},\Delta,\bm{x}_{-i},\bm{y}))\]
where $T^o(x_i, (p,q,\bm{l},\bm{u},\bm{w},\Delta,\bm{x}_{-i},\bm{y}))$ stands for the Bouligand tangent cone at $x_i$ to the feasible set of traveler $i$'s problem parameterized by $(p,q,\bm{l},\bm{u},\bm{w},\Delta,\bm{x}_{-i},\bm{y})$ (actually parameterized by $u_i$) and $L^o(x_i, (p,q,\bm{l},\bm{u},\bm{w},\Delta,\bm{x}_{-i},\bm{y}))$ denotes the polar to the linearized cone at $x_i$ of the feasible set of traveler $i$'s problem, which  is defined by
\[\begin{array}{rcll}
L(x_i, (p,q,\bm{l},\bm{u},\bm{w},\Delta,\bm{x}_{-i},\bm{y})) & = & \left\{d~:~ \right.& \langle\nabla g_i^k(x_i, (p,q,\bm{l},\bm{u},\bm{w},\Delta,\bm{x}_{-i},\bm{y})), d\rangle\le 0,\\[3pt]
&& &\left.\forall\,k\mbox{ such that } g_i^k(x_i, (p,q,\bm{l},\bm{u},\bm{w},\Delta,\bm{x}_{-i},\bm{y}))=0\right\}
\end{array}
\]
and the feasible set of Traveler i's problem is defined as $\{x_i~:~g_i^k(x_i, (p,q,\bm{l},\bm{u},\bm{w},\Delta,\bm{x}_{-i},\bm{y}))\le 0,~k=1,\dots\}$. In our specific case the feasible set is defined by the following three linear functions $g_i^1(x_i,u_i)=b_iQ_ix_i-B_i$, $g_i^2(x_i,u_i)=-x_i$ and $g_i^3(x_i,u_i)=x_i-u_i$. Now if $u_i=0$ then the only feasible point of traveler $i$'s problem is $x_i=0$. In this case $\nabla g_i^1(0,0)=1$, $\nabla g_i^2(0,0)=-1$ and $\nabla g_i^3(0,0)=1$ therefore $L^o(x_i, (p,q,\bm{l},(0,\bm{u}_{-i}),\bm{w},\Delta,\bm{x}_{-i},\bm{y}))=\mathbb{R}=T^o((x_i, (p,q,\bm{l},(0,\bm{u}_{-i}),\bm{w},\Delta,\bm{x}_{-i},\bm{y}))$. Let us now consider the case $u_i=1$. In this case the feasible set of Traveler $i$'s problem is the interval $I=[0,\min\left\{1,\frac{B_i}{b_iQ_i}\right\}]$. For any $x\in ]0,\min\left\{1,\frac{B_i}{b_iQ_i}\right\}[$, one clearly have that $g_i^k(x,(p,q,\bm{l},(0,\bm{u}_{-i}),\bm{w},\Delta,\bm{x}_{-i},\bm{y}))\neq 0$, for $k=1,2,3$ and thus $L((x_i, (p,q,\bm{l},\bm{u},\bm{w},\Delta,\bm{x}_{-i},\bm{y}))=\mathbb{R}$. On the other hand $T^o(x_i, (p,q,\bm{l},(0,\bm{u}_{-i}),\bm{w},\Delta,\bm{x}_{-i},\bm{y}))=\mathbb{R}$ and thus the desired equality holds true. Now if $x_i=0$ then only $g_i^2$ is null at $(0,1)$ and $\nabla g_i^2(0,1)=-1$ which implies that $L^o(0, (p,q,\bm{l},(1,\bm{u}_{-i}),\bm{w},\Delta,\bm{x}_{-i},\bm{y}))=\mathbb{R^-}=T^o(0, (p,q,\bm{l},(1,\bm{u}_{-i}),\bm{w},\Delta,\bm{x}_{-i},\bm{y}))$. Finally if $x_i=\min\left\{1,\frac{B_i}{b_iQ_i}\right\}$ then one (and only one) of the functions $g_i^1$ or $g_i^3$ is null at $(x_i,1)$. Since $\nabla g_i^1(x_i,1)=1=\nabla g_i^3(x_i,1)=1$, then one has again $L^o(x_i, (p,q,\bm{l},(1,\bm{u}_{-i}),\bm{w},\Delta,\bm{x}_{-i},\bm{y}))=\mathbb{R^-}=T^o(x_i, (p,q,\bm{l},(1,\bm{u}_{-i}),\bm{w},\Delta,\bm{x}_{-i},\bm{y}))$. The same reasoning will work for the TSP $mn$'s problem. We can thus conclude that Guignard's CQ is satisfied at any boundary opponent strategy. All assumptions of \cite{AusselSvenssonChapter}[Theorem 3.3.11] are fulfilled and thus item $b)$ of Theorem \ref{thm:equiv_MILBP} is proved.
}\Halmos
\endproof
\section{SD reformulation}
\subsection{Proof of Corollary 1}
\label{corollary1}
\begin{corollary}
\label{thm:SD_MILBP}
Assume that the perceived waiting/idle time functions of followers are linear.
\begin{itemize}
    \item [ a) ] If $(p,q,\bm{l},\bm{u},\bm{w},\Delta,\bm{x},\bm{y})$ is a global solution of \SLMFGL, then $\forall (\bm{\bar{\lambda}},\bm{\bar{\mu}})\in\bar{\Lambda}(p,q,\bm{l},\bm{u},\bm{w},\Delta,\bm{x},\bm{y})$, $(p,q,\bm{l},\bm{u},\bm{w},\Delta,\bm{x},\bm{y},\bm{\bar{\lambda}},\bm{\bar{\mu}})$ is a global solution of the associated SD-based reformulation \SDL.
    \item [ b) ] Assume additionally that 
    \begin{itemize}
    \item [ i) ] for all Traveler $i \in \mathcal{I}$, $b_i \geq p_{min}$;
    \item [ii) ] for all TSP $n \in \mathcal{N}_m, m \in \mathcal{M}$, $\beta_{mn} \leq q_{max}$.
    \end{itemize}
    If $(p,q,\bm{l},\bm{u},\bm{w},\Delta,\bm{x},\bm{y},\bm{\bar{\lambda}},\bm{\bar{\mu}})$ is a global solution of \SDL, then  $(p,q,\bm{l},\bm{u},\bm{w},\Delta,\bm{x},\bm{y})$ is a global solution of \SLMFGL.
\end{itemize}
\end{corollary}
\subsection{primal-dual formulations of \SLMFGQ}
\label{LEMMA1}
To  simplify  the derivation, we convert each TSP's maximization problem into a minimization problem and give primal-dual formulations of SLMFG-Q as follows.
\label{AppenC}
\begin{lemma}
\label{Lemma1}
The dual problem of $G_{i}$ can be written as:
\begin{subequations}
\allowdisplaybreaks
\begin{align}
&\quad\max_{\boldsymbol{\lambda}}\quad\
-\frac{1}{2P_{i}}\left[(A_{i}\bar{\lambda}_{1}^{i})^{2}+(\bar{\lambda}_{2}^{i})^{2}\right]-\frac{U_{i}}{P_{i}}\left(A_{i}\bar{\lambda}_{1}^{i}+\bar{\lambda}_{2}^{i}\right)-\frac{1}{P_{i}}\left(A_{i}\bar{\lambda}_{1}^{i}\bar{\lambda}_{2}^{i}\right)-\frac{U_{i}^{2}}{2P_{i}}-B_{i}\bar{\lambda}_{1}^{i}-\bar{\lambda}_{2}^{i}u_{i}
\label{Dual1}\\
&\text{\emph{subject to}} \nonumber \\
&\bar{\lambda}_{1}^{i},\bar{\lambda}_{2}^{i} \geq 0,\label{DL1c}
\end{align}
\end{subequations}

where $A_{i}= b_{i}Q_{i}$, $U_{i}=-\left(\frac{2\zeta\alpha_{i}Q_{i}}{\bar{C}}-\frac{2\zeta\alpha_{i}\Delta}{C^{2}}+r_{i}-b_{i}\right)Q_{i}$, $1/P_{i}=\frac{C^{2}}{2\zeta\alpha_{i}Q_{i}^{2}}$, and $\bar{\lambda}_{1}^{i}$ and $\bar{\lambda}_{2}^{i}$ are the dual variables corresponding to constraints \eqref{BL1b} and \eqref{BL1c}, respectively.

The dual problem of $H_{mn}$ can be written as:
\begin{subequations}
\begin{align}
&\max_{\bm{\mu}}\quad -\frac{1}{2P_{mn}}\left[(A_{mn}\bar{\mu}_{1}^{mn})^{2}+(\bar{\mu}_{2}^{mn})^{2}\right]-\frac{U_{mn}}{P_{mn}}\left(A_{mn}\bar{\mu}_{1}^{mn}
+\bar{\mu}_{2}^{mn}\right)-\frac{1}{P_{mn}}(A_{mn}\bar{\mu}_{1}^{mn}\bar{\mu}_{2}^{mn})\nonumber\\
&-\frac{U_{mn}^{2}}{2P_{mn}}-\bar{B}_{mn}\bar{\mu}_{1}^{mn}-w_{mn}\bar{\mu}_{2}^{mn}\label{Dual2}\\
&\text{\emph{subject to}} \nonumber\\
&\bar{\mu}_{1}^{mn},\bar{\mu}_{2}^{mn}\geq 0,
\end{align}
\end{subequations}
\noindent where  $A_{mn}=C_{mn}\gamma_{mn}$, $U_{mn}=-\left(\frac{2\kappa\eta_{mn}\Delta}{C^{2}}+\beta_{mn}-\rho_{mn}\right)C_{mn}$, and $1/P_{mn}=\frac{C^{2}}{2\kappa\eta_{mn}C_{mn}^{2}}$, and $\bar{\mu}_{1}^{mn}$ and $\bar{\mu}_{2}^{mn}$ are the dual variables corresponding to constraints \eqref{BL2b} and \eqref{BL2c}, respectively.
\end{lemma}
\proof{\textsc{Proof:}} Followers' primal problem are summarized in Eqs. \eqref{QD11}-\eqref{QD3}:
\begin{subequations}
\label{MPrimal}
\begin{align}
&(\text{Primal problem}): \quad \min_{\bm{x}}\quad \frac{1}{2} \bm{x}^{T} \mathbf{P} \bm{x}+\mathbf{U}^{T} \bm{x},\label{QD11}\\
&\quad\quad \hspace*{2.7cm}\text{subject to:}  && \nonumber \\
&\quad\quad \hspace*{2.7cm}\mathbf{A}\bm{x}\leq \mathbf{b},\label{QD1}\\
&\quad\quad \hspace*{2.7cm}\bm{x}\leq\boldsymbol{u},\label{QD2}\\
&\quad\quad \hspace*{2.7cm}\bm{x}\geq \mathbf{0}.\label{QD3}
\end{align}
\end{subequations}
\indent The Lagrangian function is given in \eq\eqref{Lar}:
\begin{equation}
\mathcal{L}_{1}(\bm{x},\boldsymbol{\lambda_{1}},\boldsymbol{\lambda_{2}},\boldsymbol{\lambda_{3}})=\inf _{\bm{x}}\left\{\frac{1}{2} \bm{x}^{T}\boldsymbol{ P} \bm{x}+\boldsymbol{U}^{T} \bm{x}+\boldsymbol{\lambda_{1}}^{T}(\boldsymbol{A}\bm{x}-\boldsymbol{b})+\boldsymbol{\lambda_{2}}^{T}(\bm{x}-\boldsymbol{\mu})-\boldsymbol{\lambda_{3}}^{T}\bm{x}\right\},
\label{Lar}
\end{equation}
where $\bm{\lambda_{1}}$, $\bm{\lambda_{2}}$ and $\bm{\lambda_{3}}$ denotes the Lagrangian multipliers corresponding to constraint \eqref{QD1}-\eqref{QD3}, respectively.

The dual problem is then given in Eqs. \eqref{46a}-\eqref{46c}:
\begin{subequations}
\begin{align}
&\max _{\bm{\lambda_{1}},\bm{\lambda_{2}},\bm{\lambda_{3}}}\quad\left\{\min _{\bm{x}}\mathcal{L}_{1}(\boldsymbol{x},\boldsymbol{\lambda_{1}},\boldsymbol{\lambda_{2}},\boldsymbol{\lambda_{3}})\right\},\label{46a}\\
&\text{subject to:}\nonumber\\
&\frac{\partial \mathcal{L}_{1}(\boldsymbol{x},\boldsymbol{\lambda_{1}},\boldsymbol{\lambda_{2}},\boldsymbol{\lambda_{3}})}{\partial \bm{x}}=\mathbf{P} \bm{x}+\mathbf{U}+\mathbf{A}^{T} \boldsymbol{\lambda_{1}}+\boldsymbol{\lambda_{2}}-\boldsymbol{\lambda_{3}}=\mathbf{0},\label{46b}\\
&\bm{\lambda_{1}},\bm{\lambda_{2}},\bm{\lambda_{3}}\geq 0.\label{46c}
\end{align}
\end{subequations}
\indent The dual constraint \eqref{46b} implies
that $\bm{x}^{T}[\mathbf{P} \bm{x}+\mathbf{U}+\mathbf{A}^{T} \boldsymbol{\lambda_{1}}+\boldsymbol{\lambda_{2}}-\boldsymbol{\lambda_{3}}]=\bm{x}^{T}\mathbf{0}$, namely,
\begin{equation}
  \bm{x}^{T}\mathbf{P} \bm{x}+\bm{x}^{T}\mathbf{U}+\boldsymbol{\lambda_{1}}^{T}\mathbf{A}\bm{x} +\boldsymbol{\lambda_{2}}^{T}\bm{x}-\boldsymbol{\lambda_{3}}^{T}\bm{x}= \mathbf{0}. 
\end{equation}

Thus the dual objective \eqref{46a} can be written as:
\begin{equation}
-\frac{1}{2}\bm{x}^{T}\mathbf{P}\bm{x}-\boldsymbol{\lambda_{1}}^{T}\mathbf{b}-\boldsymbol{\lambda_{2}}^{T}\boldsymbol{u}+(\bm{x}^{T}\mathbf{P} \bm{x}+\bm{x}^{T}\mathbf{U}+\boldsymbol{\lambda_{1}}^{T}\mathbf{A}\bm{x} +\boldsymbol{\lambda_{2}}^{T}\bm{x}-\boldsymbol{\lambda_{3}}^{T}\bm{x}).    
\end{equation}

Then the dual problem of Eqs.\eqref{QD11}-\eqref{QD3} is given in Eqs. \eqref{QDD1}-\eqref{QDD3}:
\begin{subequations}
\label{Mdual}
\begin{align}
&(\text{Dual problem}): \quad \max_{\bm{\lambda_{1},\bm{\lambda_{2}}}}\quad -\frac{1}{2}\bm{x}^{T}\mathbf{P}\bm{x}-\boldsymbol{\lambda_{1}}^{T}\mathbf{b}-\boldsymbol{\lambda_{2}}^{T}\boldsymbol{u},\label{QDD1}\\
&\quad\quad \hspace*{2.4cm}\text{subject to:}  && \nonumber \\
&\quad\quad \hspace*{2.4cm}\mathbf{P} \bm{x}+\mathbf{U}+\mathbf{A}^{T} \boldsymbol{\lambda_{1}}+\boldsymbol{\lambda_{2}}=\mathbf{0},\\
&\quad\quad \hspace*{2.4cm}\boldsymbol{\lambda_{1}},\boldsymbol{\lambda_{2}}\geq \mathbf{0}.\label{QDD3}
\end{align}
\end{subequations}

Since $\mathbf{P}^{-1}$ exists, from the dual constraints we have:
\begin{equation}
\mathbf{P} \bm{x}+\mathbf{U}+\mathbf{A}^{T} \boldsymbol{\lambda_{1}}+\boldsymbol{\lambda_{2}}=\mathbf{0}\Longrightarrow \bm{x}+\mathbf{P}^{-1} \mathbf{U}+\mathbf{P}^{-1} \mathbf{A}^{T} \boldsymbol{\lambda_{1}}+\mathbf{P}^{-1} \boldsymbol{\lambda_{2}}=\mathbf{0},
\end{equation}
\indent i.e.,
\begin{equation}
\bm{x}=-\mathbf{P}^{-1}\left( \boldsymbol{\mathbf{U}}+\mathbf{A}^{T}\boldsymbol{\lambda_{1}}+\boldsymbol{\lambda_{2}}\right). \label{51}
\end{equation}
\indent Then, substituting \eq\eqref{51}  into the dual objective
function \eqref{QDD1} and rearranging, we obtain:
\begin{subequations}
\begin{align}
&-\frac{1}{2}\bm{x}^{T}\mathbf{P}\bm{x}-\boldsymbol{\lambda_{1}}^{T}\mathbf{b}-\boldsymbol{\lambda_{2}}^{T}\boldsymbol{u}\\
&=-\frac{1}{2}[\mathbf{P}^{-1}\left( \boldsymbol{\mathbf{U}}+\mathbf{A}^{T}\boldsymbol{\lambda_{1}}+\boldsymbol{\lambda_{2}}\right) ^{T}] \cdot  \mathbf{P} \cdot [\mathbf{P}^{-1}\left( \boldsymbol{\mathbf{U}}+\mathbf{A}^{T}\boldsymbol{\lambda_{1}}+\boldsymbol{\lambda_{2}}\right)] -\mathbf{b}^{T} \boldsymbol{\lambda_{1}}-\mathbf{u}^{T} \boldsymbol{\lambda_{2}} \\
&=-\frac{1}{2}\left( \boldsymbol{\mathbf{U}}^{T}+\boldsymbol{\lambda_{1}}^{T}\mathbf{A}+\boldsymbol{\lambda_{2}}^{T}\right) \cdot  \mathbf{P}^{-1} \cdot \left( \boldsymbol{\mathbf{U}}+\mathbf{A}^{T}\boldsymbol{\lambda_{1}}+\boldsymbol{\lambda_{2}}\right)-\mathbf{b}^{T} \boldsymbol{\lambda_{1}}-\mathbf{u}^{T} \boldsymbol{\lambda_{2}} \\
&=-\frac{1}{2}\left( \boldsymbol{\mathbf{U}}^{T}\mathbf{P}^{-1}+\boldsymbol{\lambda_{1}}^{T}\mathbf{A}\mathbf{P}^{-1}+\boldsymbol{\lambda_{2}}^{T}\mathbf{P}^{-1}\right) \cdot \left( \boldsymbol{\mathbf{U}}
+\mathbf{A}^{T}\boldsymbol{\lambda_{1}}
+\boldsymbol{\lambda_{2}}\right) -\mathbf{b}^{T} \boldsymbol{\lambda_{1}}-\mathbf{u}^{T} \boldsymbol{\lambda_{2}} \\
&=-\frac{1}{2}(\boldsymbol{\lambda_{1}}^{T}\mathbf{A}\mathbf{P}^{-1}\mathbf{A}^{T}\boldsymbol{\lambda_{1}}+\boldsymbol{\lambda_{2}}^{T}\mathbf{P}^{-1}\boldsymbol{\lambda_{2}})-\frac{1}{2}(2\boldsymbol{\mathbf{U}}^{T}\mathbf{P}^{-1} \mathbf{A}^{T}\boldsymbol{\lambda_{1}}
+2\boldsymbol{\mathbf{U}}^{T}\mathbf{P}^{-1}\boldsymbol{\lambda_{2}}+2\boldsymbol{\lambda_{1}}^{T}\mathbf{A}\mathbf{P}^{-1}{\boldsymbol{\lambda_{2}}}) \nonumber
\\
&-\frac{1}{2}\boldsymbol{\mathbf{U}}^{T}\mathbf{P}^{-1}\boldsymbol{\mathbf{U}}-\mathbf{b}^{T} \boldsymbol{\lambda_{1}}-\mathbf{u}^{T}\boldsymbol{\lambda_{2}}  \\
&=-\frac{1}{2}(\boldsymbol{\lambda_{1}}^{T}\mathbf{A}\mathbf{P}^{-1}\mathbf{A}^{T}\boldsymbol{\lambda_{1}}+\boldsymbol{\lambda_{2}}^{T}\mathbf{P}^{-1}\boldsymbol{\lambda_{2}})-(\boldsymbol{\mathbf{U}}^{T}\mathbf{P}^{-1} \mathbf{A}^{T}\boldsymbol{\lambda_{1}}
+\boldsymbol{\mathbf{U}}^{T}\mathbf{P}^{-1}\boldsymbol{\lambda_{2}}+\boldsymbol{\lambda_{1}}^{T}\mathbf{A}\mathbf{P}^{-1}{\boldsymbol{\lambda_{2}}}) \nonumber
\\
&-\frac{1}{2}\boldsymbol{\mathbf{U}}^{T}\mathbf{P}^{-1}\boldsymbol{\mathbf{U}}-\mathbf{b}^{T} \boldsymbol{\lambda_{1}}-\mathbf{u}^{T}\boldsymbol{\lambda_{2}}.  
\end{align}
\end{subequations}

The dual problem can be written as:
\begin{subequations}
\begin{align}
& \min_{\bm{\lambda_{1}},\bm{\lambda_{2}}}\quad -\frac{1}{2}\mathbf{P}^{-1}(\boldsymbol{\lambda_{1}}^{T}\mathbf{A}\mathbf{A}^{T}\boldsymbol{\lambda_{1}}+\boldsymbol{\lambda_{2}}^{T}\boldsymbol{\lambda_{2}})-\mathbf{U}^{T}\mathbf{P}^{-1}( \mathbf{A}^{T}\boldsymbol{\lambda_{1}}
+\boldsymbol{\lambda_{2}})-\boldsymbol{\lambda_{1}}^{T}\mathbf{A}\mathbf{P}^{-1}{\boldsymbol{\lambda_{2}}}-\frac{1}{2}\boldsymbol{\mathbf{U}}^{T}\mathbf{P}^{-1}\boldsymbol{\mathbf{U}}-\mathbf{b}^{T} \boldsymbol{\lambda_{1}}-\mathbf{u}^{T}\boldsymbol{\lambda_{2}},\\
&\text{subject to:}  && \nonumber \\
&\boldsymbol{\lambda_{1}},\boldsymbol{\lambda_{2}}\geq \mathbf{0}.
\end{align}
\end{subequations} 
\indent Thus concluding the proof. \Halmos
\endproof
\subsection{Proof of Corollary 2}
\label{corollary2}
In Corollary \ref{thm:SD_MILBP}, the set $\bar{\Lambda}(\cdot)$ stands for the set of dual variables associated to the SD conditions of the followers (travelers and TSPs).
\begin{corollary}
\label{thm:SD_MIQBP}
Assume that the perceived waiting/idle time functions of followers are quadratic.
\begin{itemize}
    \item [ a) ] If $(p,q,\bm{l},\bm{u},\bm{w},\Delta,\bm{x},\bm{y})$ is a global solution of \SLMFGQ, then $\forall (\bm{\bar{\lambda}},\bm{\bar{\mu}})\in\bar{\Lambda}(p,q,\bm{l},\bm{u},\bm{w},\Delta,\bm{x},\bm{y})$, $(p,q,\bm{l},\bm{u},\bm{w},\Delta,\bm{x},\bm{y},\bm{\bar{\lambda}},\bm{\bar{\mu}})$ is a global solution of the associated SD-based reformulation \SDQ.
    \item [ b) ] Assume additionally that 
    \begin{itemize}
    \item [ i) ] for all Traveler $i \in \mathcal{I}$, $b_i \geq p_{min}$;
    \item [ii) ] for all TSP $n \in \mathcal{N}_m, m \in \mathcal{M}$, $\beta_{mn} \leq q_{max}$.
    \end{itemize}
    If $(p,q,\bm{l},\bm{u},\bm{w},\Delta,\bm{x},\bm{y},\bm{\bar{\lambda}},\bm{\bar{\mu}})$ is a global solution of \SDQ, then  $(p,q,\bm{l},\bm{u},\bm{w},\Delta,\bm{x},\bm{y})$ is a global solution of \SLMFGQ.
\end{itemize}
\end{corollary}
\subsection{McCormick envelope relaxation}
We relax these bilinear terms by introducing auxiliary variables  $W_{i} \triangleq \Delta \cdot x_{i}$ and $W_{mn}\triangleq\Delta \cdot y_{mn}$  using their McCormick envelopes \citep{mccormick1976computability}. Accordingly, constraints \eqref{MKa}-\eqref{MKh} are added to the SD reformulation of \SLMFGQ.
\begin{subequations}
\allowdisplaybreaks
\begin{align}
&0\leq W_{i}\leq C, &&\forall i \in \mathcal{I},\tag{MK.1}\label{MKa}\\
&W_{i}\geq \Delta+Cx_{i}-C, &&\forall i \in \mathcal{I},\tag{MK.2}\label{MKb}\\
&W_{i}\leq \Delta+\underline{C}x_{i}-\underline{C}, &&\forall i \in \mathcal{I},\tag{MK.3}\label{MKc}\\
&W_{i}\leq Cx_{i}, &&\forall i \in \mathcal{I},\tag{MK.4}\label{MKd}\\
&0\leq W_{mn}\leq C,  && m \in \mathcal{M},n \in \mathcal{N}_{m},\tag{MK.5}\label{MKe}\\
&W_{mn}\geq \Delta+Cy_{mn}-C,&& m \in \mathcal{M},n \in \mathcal{N}_{m},\tag{MK.6}\label{MKf}\\
&W_{mn}\leq \Delta+\underline{C}y_{mn}-\underline{C},&& m \in \mathcal{M},n \in \mathcal{N}_{m},\tag{MK.7}\label{MKg}\\
&W_{mn}\leq C y_{mn},&& m \in \mathcal{M},n \in \mathcal{N}_{m},\tag{MK.8}\label{MKh}
\end{align}
\end{subequations}
\vspace{12pt}
\section{Formulations of the sub-problems}
\label{sub-problems formulation}
\noindent\textbf{\noindent (The $\bm{kth}$ sub-problem $\textit{SP}^{k}$ of \SLMFGL)}
\begin{subequations}
\label{SP1}
\allowdisplaybreaks
\begin{align}
&\max_{p,q,\bm{l},\bm{u},\bm{\tau},\Delta,\bm{x},\bm{y}}F= \sum_{i \in \mathcal{I}} b_{i}Q_{i}x_{i}-\sum_{m \in \mathcal{M}}\sum_{n \in \mathcal{N}_{m}} \beta_{mn}c_{mn} y_{mn},\tag{M.1}\label{E1a}\\
&\text{\emph{subject to:}}  && \nonumber \\
&\text{MaaS regulator's constraints:}&& \nonumber\\
&\sum_{m \in \mathcal{M}} v_{m} l_{i}^{m} = D_{i} x_{i},&& \forall i \in \mathcal{I},\tag{M.2}\\
& 0\leq\sum_{m \in \mathcal{M}}l_{i}^{m}-T_{i}x_{i}\leq {R}_{i},&& \forall i \in \mathcal{I},\tag{M.3}\\
&\sum_{m \in  \mathcal{M}}\sigma_{m} l_{i}^{m} \leq  \Gamma_{i},&& \forall i \in \mathcal{I}, \tag{M.4}\\
&\sum_{i \in \mathcal{I}} v_{m}^{2}l_{i}^{m}\leq \sum_{n \in \mathcal{N}_{m}} C_{mn}y_{mn},&& \forall m \in  \mathcal{M},\tag{M.5}\\
&b_{i}-p\geq(1-u_{i})(b_{i}-p_{max}),&& \forall i \in \mathcal{I},\tag{M.6}\\
&b_{i}-p\leq u_{i}(b_{i}-p_{min}),&& \forall i \in \mathcal{I},\tag{M.7}\\
&q-\beta_{mn}\geq(1-w_{mn})(q_{min}-\beta_{mn}),&& \forall m \in \mathcal{M},n \in \mathcal{N}_{m},\tag{M.8}\\
&q-\beta_{mn}\leq w_{mn}(q_{max}-\beta_{mn}),&& \forall m \in \mathcal{M},n \in \mathcal{N}_{m},\tag{M.9}\\
&\Delta= \sum_{m \in \mathcal{M}}\sum_{n \in \mathcal{N}_{m}}C_{mn} y_{mn} - \sum_{i \in \mathcal{I}} Q_{i} x_{i},\tag{M.10}\\
&u_{i}\in \left\{0,1\right\},&& \forall i \in \mathcal{I},\tag{M.11} \\
&w_{mn}\in \left\{0,1\right\},&& \forall m \in \mathcal{M},n \in \mathcal{N}_{m},\tag{M.12}\\
&l_{i}^{m}\geq 0, && \forall i \in \mathcal{I},m \in \mathcal{M},\tag{M.13}\\
&p_{min}\leq p\leq p_{max},\tag{M.14}\\
&q_{min}\leq q\leq q_{max},\tag{M.15}\\
&\underline{C}\leq\Delta\leq \overline{C},\tag{M.16}\\
& \text{SD conditions of travelers:}\nonumber\\
&b_{i}Q_{i}x_{i}\leq B_{i},&&  \forall i \in\mathcal{I},\tag{T.2}\\
&x_{i} \leq u_{i},&&\forall i \in\mathcal{I},\tag{T.3}\\
&x_{i}\geq 0,&&\forall i \in\mathcal{I},\tag{T.4}\\
&u_{i}=0, &&  \forall i \in \overline{\mathcal{I}}_{0}^{k},\label{E1c} \\
& u_{i}=1, &&  \forall i \in \overline{\mathcal{I}}_{1}^{k},\label{E1d}  \\
&(r_{i}-b_{i})Q_{i}x_{i}\geq B_{i}\bar{\lambda}_{1}^{i}+u_{i}\bar{\lambda}_{2}^{i}, &&  \forall i \in \overline{\mathcal{I}}_{1}^{k},\label{E1e} \\
&b_{i}Q_{i}\bar{\lambda}_{1}^{i}+\bar{\lambda}_{2}^{i}\geq (r_{i}-b_{i})Q_{i},&&  \forall i \in \overline{\mathcal{I}}_{1}^{k},\label{E1f}\\
&\bar{\lambda}_{1}^{i},\bar{\lambda}_{2}^{i}\geq 0,  &&\forall i \in \overline{\mathcal{I}}_{1}^{k},\label{E1g}\\
&\text{SD conditions of TSPs:}\nonumber\\
&C_{mn}\gamma_{mn}y_{mn}\leq \bar{B}_{mn},&&\forall (m,n) \in\mathcal{MN},\tag{P.2}\\
&y_{mn}\leq w_{mn},&&\forall (m,n) \in\mathcal{MN},\tag{P.3}\\
&y_{mn}\geq 0,&&\forall (m,n) \in\mathcal{MN},\tag{P.4}\\
&\tau_{mn}=0,&&\forall (m,n) \in \overline{\mathcal{MN}}_{0}^{k},\label{E1i}\\
&\tau_{mn}=1,&&\forall (m,n) \in \overline{\mathcal{MN}}_{1}^{k},\label{E1j}\\
&(\beta_{mn}-\rho_{mn})c_{mn}y_{mn}\geq \bar{B}_{mn}\bar{\mu}_{1}^{mn}+\tau_{mn}\bar{\mu}_{2}^{mn},&&\forall (m,n) \in \overline{\mathcal{MN}}_{1}^{k},\label{E1k}\\
& c_{mn}\gamma_{mn}\bar{\mu}_{1}^{mn}+\bar{\mu}_{2}^{mn}\geq (\beta_{mn}-\rho_{mn})c_{mn},&&\forall (m,n) \in \overline{\mathcal{MN}}_{1}^{k},\label{E1l}\\
&\bar{\mu}_{1}^{mn},\bar{\mu}_{2}^{mn}\geq 0 , &&\forall (m,n) \in \overline{\mathcal{MN}}_{1}^{k},\label{E1m}
\end{align}
\end{subequations}\\

\textbf{\noindent (The $\bm{kth}$ sub-problem $\textit{SP}^{k}$ of \SLMFGQ)}
\begin{subequations}
\label{SP2}
\allowdisplaybreaks
\begin{align}
&\max_{p,q,\bm{l},\bm{u},\bm{\tau},\Delta,\bm{x},\bm{y}}F= \sum_{i \in \mathcal{I}} b_{i}Q_{i}x_{i}-\sum_{m \in \mathcal{M}}\sum_{n \in \mathcal{N}_{m}} \beta_{mn}c_{mn} y_{mn},&&\label{F1a} \\
&\text{\emph{subject to:}}  && \nonumber \\
&\text{MaaS regulator's constraints:}&& \nonumber\\
&\sum_{m \in \mathcal{M}} v_{m} l_{i}^{m} = D_{i} x_{i},&& \forall i \in \mathcal{I},\tag{M.2}\\
& 0\leq\sum_{m \in \mathcal{M}}l_{i}^{m}-T_{i}x_{i}\leq {R}_{i},&& \forall i \in \mathcal{I},\tag{M.3}\\
&\sum_{m \in  \mathcal{M}}\sigma_{m} l_{i}^{m} \leq  \Gamma_{i},&& \forall i \in \mathcal{I}, \tag{M.4}\\
&\sum_{i \in \mathcal{I}} v_{m}^{2}l_{i}^{m}\leq \sum_{n \in \mathcal{N}_{m}} C_{mn}y_{mn},&& \forall m \in  \mathcal{M},\tag{M.5}\\
&b_{i}-p\geq(1-u_{i})(b_{i}-p_{max}),&& \forall i \in \mathcal{I},\tag{M.6}\\
&b_{i}-p\leq u_{i}(b_{i}-p_{min}),&& \forall i \in \mathcal{I},\tag{M.7}\\
&q-\beta_{mn}\geq(1-w_{mn})(q_{min}-\beta_{mn}),&& \forall m \in \mathcal{M},n \in \mathcal{N}_{m},\tag{M.8}\\
&q-\beta_{mn}\leq w_{mn}(q_{max}-\beta_{mn}),&& \forall m \in \mathcal{M},n \in \mathcal{N}_{m},\tag{M.9}\\
&\Delta= \sum_{m \in \mathcal{M}}\sum_{n \in \mathcal{N}_{m}}C_{mn} y_{mn} - \sum_{i \in \mathcal{I}} Q_{i} x_{i},\tag{M.10}\\
&u_{i}\in \left\{0,1\right\},&& \forall i \in \mathcal{I},\tag{M.11} \\
&w_{mn}\in \left\{0,1\right\},&& \forall m \in \mathcal{M},n \in \mathcal{N}_{m},\tag{M.12}\\
&l_{i}^{m}\geq 0, && \forall i \in \mathcal{I},m \in \mathcal{M},\tag{M.13}\\
&p_{min}\leq p\leq p_{max},\tag{M.14}\\
&q_{min}\leq q\leq q_{max},\tag{M.15}\\
&\underline{C}\leq\Delta\leq \overline{C},\tag{M.16}\\
& \text{SD conditions of travelers:}\nonumber\\
&b_{i}Q_{i}x_{i}\leq B_{i},&&  \forall i \in\mathcal{I},\tag{T.2}\\
&x_{i} \leq u_{i},&&\forall i \in\mathcal{I},\tag{T.3}\\
&x_{i}\geq 0,&&\forall i \in\mathcal{I},\tag{T.4}\\
&u_{i}=0, &&  \forall i \in \overline{\mathcal{I}}_{0}^{k},\label{F1c} \\
&u_{i}=1, &&  \forall i \in \overline{\mathcal{I}}_{1}^{k},\label{F1d}  \\
&\frac{\zeta\alpha_{i}Q_{i}^{2}}{C^{2}}x_{i}^{2}-\left(\frac{2\zeta\alpha_{i}Q_{i}}{\bar{C}}-\frac{2\zeta\alpha_{i}\Delta}{C^{2}}+r_{i}-b_{i}\right)Q_{i}x_{i}\leq-\frac{(A_{i}\bar{\lambda}_{1}^{i})^{2}+(\bar{\lambda}_{2}^{i})^{2}}{2P_{i}}\nonumber\\
&-\frac{U_{i}}{P_{i}}\left(b_{i}Q_{i}\bar{\lambda}_{1}^{i}+\bar{\lambda}_{2}^{i}\right)-\frac{b_{i}Q_{i}\bar{\lambda}_{1}^{i}\bar{\lambda}_{2}^{i}}{P_{i}}-\frac{U_{i}^{2}}{2P_{i}}-B_{i}\bar{\lambda}_{1}^{i}-\bar{\lambda}_{2}^{i}u_{i}, &&  \forall i \in \overline{\mathcal{I}}_{1}^{k},\label{F1e} \\
&\bar{\lambda}_{1}^{i},\bar{\lambda}_{2}^{i}\geq 0,  &&\forall i \in \overline{\mathcal{I}}_{1}^{k},\label{F1f}\\
&\text{McCormick envelopes Relaxation:} \nonumber\\
&0\leq W_{i}\leq C, &&\forall i \in \mathcal{I},\tag{MK.1}\\
&W_{i}\geq \Delta+Cx_{i}-C, &&\forall i \in \overline{\mathcal{I}}_{1}^{k},\tag{MK.2}\\
&W_{i}\leq \Delta+\underline{C}x_{i}-\underline{C}, &&\forall i \in \overline{\mathcal{I}}_{1}^{k},\tag{MK.3}\\
&W_{i}\leq Cx_{i}, &&\forall i \in \overline{\mathcal{I}}_{1}^{k},\tag{MK.4}\\
&\text{SD conditions of TSPs:}\nonumber\\
&C_{mn}\gamma_{mn}y_{mn}\leq \bar{B}_{mn},&&\forall (m,n) \in\mathcal{MN},\tag{P.2}\\
&y_{mn}\leq w_{mn},&&\forall (m,n) \in\mathcal{MN},\tag{P.3}\\
&y_{mn}\geq 0,&&\forall (m,n) \in\mathcal{MN},\tag{P.4}\\
&\tau_{mn}=0,&&\forall (m,n) \in \overline{\mathcal{MN}}_{0}^{k},\label{F1i}\\
&\tau_{mn}=1,&&\forall (m,n) \in \overline{\mathcal{MN}}_{1}^{k},\label{F1j}\\
&\frac{\kappa\eta_{mn} C_{mn}^{2}}{C^{2}}y_{mn}^{2}-\left(\frac{2\kappa\eta_{mn}\Delta}{C^{2}}+\beta_{mn}-\rho_{mn}\right)C_{mn}y_{mn}\leq\nonumber\\
&-\frac{(A_{mn}\bar{\mu}_{1}^{mn})^{2}+(\bar{\mu}_{2}^{mn})^{2}}{2P_{mn}}-\frac{U_{mn}}{P_{mn}}\left(A_{mn}\bar{\mu}_{1}^{mn}+\bar{\mu}_{2}^{mn}\right)-\frac{A_{mn}\bar{\mu}_{1}^{mn}\bar{\mu}_{2}^{mn}}{P_{mn}}\nonumber \\
&-\frac{U_{mn}^{2}}{2P_{mn}}-\bar{B}_{mn}\bar{\mu}_{1}^{mn}-w_{mn}\bar{\mu}_{2}^{mn}, &&\forall (m,n) \in \overline{\mathcal{MN}}_{1}^{k},\label{F1k}\\
&\bar{\mu}_{1}^{mn},\bar{\mu}_{2}^{mn}\geq 0, &&\forall (m,n) \in \overline{\mathcal{MN}}_{1}^{k},\label{F1l}\\
&\text{McCormick envelopes Relaxation:}\nonumber\\
&0\leq W_{mn}\leq C,  &&\forall (m,n) \in \overline{\mathcal{MN}}_{1}^{k},\tag{MK.5}\\
&W_{mn}\geq \Delta+Cy_{mn}-C,&&\forall (m,n) \in \overline{\mathcal{MN}}_{1}^{k},\tag{MK.6}\\
&W_{mn}\leq \Delta+\underline{C}y_{mn}-\underline{C},&&\forall (m,n) \in \overline{\mathcal{MN}}_{1}^{k},\tag{MK.7}\\
&W_{mn}\leq C y_{mn},&&\forall (m,n) \in \overline{\mathcal{MN}}_{1}^{k},\tag{MK.8}
\end{align}
\end{subequations}
\vspace{12pt}
\section{Computation study}
\label{Numerical results}
\subsection{Input data}
\label{input data}
MaaS requests of both travelers and TSPs are input information to the proposed two-sided MaaS market. To obtain travelers' and TSPs' request data, we conduct stochastic simulations on travelers and TSPs in the two-sided MaaS market. We consider five types of travel modes with different commercial speed (km/min) and inconvenience cost per unit of time (\$/min) given in \T\ref{T1}. 
\begin{table}[ht]
\setlength{\abovecaptionskip}{0pt}
\setlength{\belowcaptionskip}{0pt}
\centering
\footnotesize
\caption{The commercial speed and inconvenience cost per unit time  of different travel modes}
\begin{tabular}{cccccc}
\toprule
\multirow{2}{*}{Modes} & $\mathrm{m}=1$ & $\mathrm{m}=2$ & $\mathrm{m}=3$ & $\mathrm{m}=4$ & $\mathrm{m}=5$ \\
&Taxi&Ride sharing with 2 riders&Ride sharing with 3 riders&Public transit&Bicycle-sharing\\ 
\midrule
Average speed & $v_{1}$ & $v_{2}$ & $v_{3}$ & $v_{4}$ & $v_{5}$ \\
$(\mathrm{km} / \mathrm{min})$ & 0.5 & 0.27 & 0.24 & 0.21&0.18 \\
\hline Inconvenience & $\delta_{1}$ & $\delta_{2}$ & $\delta_{3}$ & $\delta_{4}$ & $\delta_{5}$ \\
cost (\$)/min& 0 & 0.1 & 0.2&0.5&1.5 \\
\bottomrule
\end{tabular}
\label{T1}
\end{table}

Recall that Traveler $i$'s MaaS request is the tuple $\mathcal{B}_{i}=(D_{i},T_{i},b_{i}, \Phi_{i}, \Gamma_{i})$, $\forall i\in\mathcal{I}$. Traveler $i$'s requested distance $D_{i}$ (km) is randomly generated within [1, 18]. Traveler $i$'s requested total service time $T_{i}$ (min) is randomly generated within $\left[\frac{D_{i}}{v_{5}},\frac{D_{i}}{v_{1}}\right]$. Travelers' purchase bidding price ($b_{i}$) is set based on the tariff of the transportation system in Sydney, Australia including: Uber, metro, bus, Tram and Lime. $b_{min}$ is set based on the price of public transit in Sydney\footnote{Opal Trip Planner, \url{https://transportnsw.info/trip##/trip}, is used to estimate the fare of different public travel modes in NSW, Australia. } and $b_{max}$ is set based on the price of UberX in Sydney, which varies over the time during one day\footnote{Real-time Uber Estimator provides real-time fare estimates on each trip.\url{https://uberestimator.com} }; thus $b_{i}$ (\$) is randomly generated within $[b_{\min},b_{\max}]$. Traveler $i$'s unit reserve cost $r_{i}$ (\$) is randomly generated within $[b_{i}+1.5,b_{i}+4]$, her travel cost budget $B_{i}$ (\$) is randomly generated within $[b_{i}Q_{i}-0.5,b_{i}Q_{i}+1.5]$, her unit waiting time cost $\alpha_{i}$ (\$) is randomly generated within $[0.02,0.05]$. Considering that a traveler's inconvenience tolerance and travel delay budget will decrease as her bidding price increases,  $\Phi_{i}$ (\$) is randomly generated within $\left[0,\frac{100}{b_{i}}\right]$ and $\Gamma_{i}$ (\$) is randomly generated within $\left[0,\frac{100 D_{i}}{b_{i}}\right]$. 

For TSPs' data, recall that TSP $mn$'s MaaS request is the tuple $\mathcal{B}_{mn}=(\beta_{mn},C_{mn},\bar{B}_{mn} )$, $\forall m\in \mathcal{M}, n\in \mathcal{N}_{m}$. Let $\mathcal{M}=\left\{1,2,3,4,5\right\}$ denote the set of travel modes shown in \T\ref{T1}, TSP $mn$'s sell-bidding price $\beta_{mn}$ (\$) is generated within $[\underline{\beta_{m}},\overline{\beta_{m}}]$, TSP $mn$'s capacity $C_{mn} (\text{km}^{2}/min)$ is randomly generated within $[\underline{C_{m}},\overline{C_{m}}]$, where the values of $\underline{C_{m}},\overline{C_{m}}$, $\underline{\beta_{m}}$ and $\overline{\beta_{m}}$ are given in \T\ref{T2}. TSP $mn$'s operating cost budget $\bar{B}_{mn}$ (\$) is randomly generated within $[\beta_{mn}C_{mn}-0.75, \beta_{mn}C_{mn}+2.75]$. TSP $mn$'s unit idle time cost $\eta_{mn}$ (\$) is randomly generated within $[1,2]$, and TSP $mn$'s reserve price $\rho_{mn}$ (\$) is randomly generated within $[\beta_{mn}-0.5,\beta_{mn}-1]$. 

\begin{table}[ht]
\setlength{\abovecaptionskip}{0pt}
\setlength{\belowcaptionskip}{0pt}
\centering
\footnotesize
\caption{The values of $\bm{\underline{\beta_{m}}}(\$)$, $\bm{\overline{\beta_{m}}} (\$)$, $\bm{\underline{C_{m}}} (km^{2}/min)$ and $\bm{\overline{C_{m}}} (km^{2}/min)$}
\begin{tabular}{@{}cccccccccccccccccccc@{}}
\toprule
\multicolumn{4}{c}{m=1}  & \multicolumn{4}{c}{m=2}  & \multicolumn{4}{c}{m=3}  & \multicolumn{4}{c}{m=4}  & \multicolumn{4}{c}{m=5} \\
\cmidrule(l){1-4} \cmidrule(l){5-8} \cmidrule(l){9-12} \cmidrule(l){13-16}\cmidrule(l){17-20}  
$\underline{\beta_{1}}$& $\overline{\beta_{1}}$ & $\underline{C_{1}}$ & $\overline{C_{1}}$  & $\underline{\beta_{2}}$& $\overline{\beta_{2}}$ & $\underline{C_{2}}$ & $\overline{C_{2}}$  & $\underline{\beta_{3}}$& $\overline{\beta_{3}}$ & $\underline{C_{3}}$ & $\overline{C_{3}}$  & $\underline{\beta_{4}}$& $\overline{\beta_{4}}$ & $\underline{C_{4}}$ & $\overline{C_{4}}$  & $\underline{\beta_{5}}$& $\overline{\beta_{5m}}$ & $\underline{C_{5}}$ & $\overline{C_{5}}$ \\
8     & 12    & 50 & 100 & 7     & 8     & 50 & 100 & 5     & 7     & 30 & 100 & 3     & 5     & 30 & 100 & 1     & 2     & 10 & 50\\
\bottomrule
\end{tabular}
\label{T2}
\end{table}

Based on the above parameters setting, we randomly generate 20 group of instances. Each group of instances corresponds to a fixed number of travelers $|\mathcal{I}|$ and TSPs $|\mathcal{MN}|$ and is named `MaaS-$|\mathcal{I}|$-$|\mathcal{MN}|'$, e.g. `MaaS-10-5' corresponds to the group of instances with 10 travelers and 5 TSPs, $|\mathcal{N}_{1}|$$\sim$$|\mathcal{N}_{5}|=1$; `MaaS-200-10' corresponds to the group of instances with 200 travelers and 10 TSPs, $|\mathcal{N}_{1}|$$\sim$$|\mathcal{N}_{5}|=2$.  We generate 20 random instances for each group of instances, hence a total of 400 instances are generated. 
\newpage
\subsection{Computational Performance of B\&B algorithms}
\setlength{\abovecaptionskip}{0pt}
\setlength{\belowcaptionskip}{0pt}
\begin{table}[ht!]
\label{ALtable}
\centering
\scriptsize
\caption{Performance of the Branch and Bound algorithms for \SLMFGL (MILBP)}
\begin{tabular}{@{\extracolsep{\fill}}cccccccccccccccccc@{}}
\toprule
Instances &
   &
  \multicolumn{4}{c}{\Bard} &
  \multicolumn{4}{c}{\BP} &
  \multicolumn{4}{c}{\Diff} &
  \multicolumn{4}{c}{\Weight} \\ \cmidrule(l){3-6} \cmidrule(l){7-10} \cmidrule(l){11-14} \cmidrule(l){15-18}
\multicolumn{1}{c}{\textit{}} &
  \multicolumn{1}{c}{\textit{}} &
  \multicolumn{1}{c}{\textit{k}} &
  \multicolumn{1}{c}{\textit{LB}} &
  \multicolumn{1}{c}{\textit{Gap(\%)}} &
  \multicolumn{1}{c}{\textit{$T (s)$ }} &
  \multicolumn{1}{c}{\textit{k}} &
  \multicolumn{1}{c}{\textit{LB}} &
  \multicolumn{1}{c}{\textit{Gap(\%)}} &
  \multicolumn{1}{c}{\textit{$T (s)$ }} &
  \multicolumn{1}{c}{\textit{k}} &
  \multicolumn{1}{c}{\textit{LB}} &
  \multicolumn{1}{c}{\textit{Gap(\%)}} &
  \multicolumn{1}{c}{\textit{$T (s)$ }} &
  \multicolumn{1}{c}{\textit{k}} &
  \multicolumn{1}{c}{\textit{LB}} &
  \multicolumn{1}{c}{\textit{Gap(\%)}} &
  \multicolumn{1}{c}{\textit{$T (s)$ }} \\ \cmidrule(l){3-6} \cmidrule(l){7-10} \cmidrule(l){11-14} \cmidrule(l){15-18}
\multirow{3}{*}{MaaS-10-5} &
  1 &
  52 &
  19 &
  0 &
  13 &
  32 &
  19 &
  0 &
  8 &
  26 &
  19 &
  0 &
  7 &
  26 &
  19 &
  0 &
  \textbf{6} \\
 &
  2 &
  63 &
  0 &
  0 &
  15 &
  34 &
  0 &
  0 &
  9 &
  26 &
  0 &
  0 &
  7 &
  26 &
  0 &
  0 &
  \textbf{6} \\
 &
  3 &
  47 &
  0 &
  0 &
  12 &
  40 &
  0 &
  0 &
  10 &
  28 &
  0 &
  0 &
  7 &
  28 &
  0 &
  0 &
  \textbf{7} \\ \hline
\multirow{3}{*}{MaaS-20-5} &
  1 &
  40 &
  226 &
  0 &
  11 &
  32 &
  226 &
  0 &
  8 &
  20 &
  226 &
  0 &
  8 &
  20 &
  226 &
  0 &
  \textbf{7} \\
 &
  2 &
  88 &
  40 &
$<1e^{-6}$ &
  23 &
  34 &
  40 &
$<1e^{-6}$ &
  \textbf{9} &
  44 &
  40 &
$<1e^{-6}$ &
  14 &
  44 &
  40 &
$< 1e^{-6}$ &
  12 \\
 &
  3 &
  104 &
  39 &
  0 &
  27 &
  40 &
  39 &
  0 &
  \textbf{10} &
  52 &
  39 &
  0 &
  16 &
  46 &
  39 &
  0 &
  13 \\\hline
\multirow{3}{*}{MaaS-30-5} &
  1 &
  76 &
  310 &
  0 &
  24 &
  33 &
  310 &
  0 &
  14 &
  37 &
  310 &
  0 &
  22 &
  25 &
  310 &
  0 &
  \textbf{11} \\
 &
  2 &
  92 &
  315 &
  0 &
  28 &
  39 &
  315 &
  0 &
  17 &
  46 &
  315 &
  0 &
  23 &
  33 &
  315 &
  0 &
  \textbf{13} \\
 &
  3 &
  88 &
  314 &
  \textless{}$1e^{-6}$ &
  27 &
  41 &
  314 &
  0 &
  16 &
  43 &
  314 &
  0 &
  23 &
  37 &
  314 &
  0 &
  \textbf{14} \\ \hline
\multirow{3}{*}{MaaS-40-5} &
  1 &
  144 &
  343 &
  0 &
  50 &
  71 &
  343 &
  \textless{}$1e^{-6}$ &
  32 &
  72 &
  343 &
  \textless{}$1e^{-6}$ &
  42 &
  63 &
  343 &
  \textless{}$1e^{-6}$ &
  \textbf{28} \\
 &
  2 &
  192 &
  216 &
  0 &
  65 &
  73 &
  216 &
  0 &
  34 &
  96 &
  216 &
  0 &
  54 &
  63 &
  216 &
  0 &
  \textbf{27} \\
 &
  3 &
  196 &
  272 &
  \textless{}$1e^{-6}$ &
  66 &
  87 &
  272 &
  0 &
  39 &
  98 &
  272 &
  0 &
  55 &
  74 &
  272 &
  0 &
  \textbf{32} \\\hline
\multirow{3}{*}{MaaS-50-5} &
  1 &
  312 &
  127 &
  0 &
  128 &
  137 &
  127 &
  0 &
  \textbf{65} &
  155 &
  127 &
  \textless{}$1e^{-6}$ &
  86 &
  135 &
  127 &
  \textless{}$1e^{-6}$ &
  67 \\
 &
  2 &
  340 &
  243 &
  0 &
  136 &
  136 &
  243 &
  0 &
  \textbf{70} &
  216 &
  243 &
  \textless{}$1e^{-6}$ &
  115 &
  145 &
  243 &
  \textless{}$1e^{-6}$ &
  72 \\
 &
  3 &
  236 &
  578 &
  \textless{}$1e^{-6}$ &
  108 &
  58 &
  578 &
  0 &
  \textbf{30} &
  235 &
  578 &
  \textless{}$1e^{-6}$ &
  106 &
  88 &
  578 &
  \textless{}$1e^{-6}$ &
  40 \\\hline
\multirow{3}{*}{MaaS-60-5} &
  1 &
  588 &
  523 &
  \textless{}$1e^{-6}$ &
  280 &
  58 &
  523 &
  0 &
  \textbf{39} &
  294 &
  523 &
  0 &
  243 &
  62 &
  523 &
  0 &
  46 \\
 &
  2 &
  644 &
  482 &
  0 &
  299 &
  89 &
  482 &
  0 &
  53 &
  322 &
  482 &
  0 &
  167 &
  83 &
  482 &
  0 &
  \textbf{51} \\
 &
  3 &
  440 &
  382 &
  0 &
  211 &
  95 &
  382 &
  0 &
  59 &
  219 &
  382 &
  \textless{}$1e^{-6}$ &
  127 &
  91 &
  382 &
  0 &
  \textbf{57} \\  \hline
\multirow{3}{*}{MaaS-70-6} &
  1 &
  11456 &
  499 &
  0 &
  6823 &
  262 &
  499 &
  0 &
  \textbf{202} &
  3890 &
  499 &
  0 &
  4569 &
  268 &
  499 &
  0 &
  221 \\
 &
  2 &
  11668 &
  743 &
  0 &
  7021 &
  78 &
  743 &
  0 &
  \textbf{59} &
  3011 &
  743 &
  0 &
  3895 &
  89 &
  743 &
  0 &
  68 \\
 &
  3 &
  11255 &
  603 &
  0 &
  6802 &
  133 &
  603 &
  \textless{}$1e^{-6}$ &
  \textbf{108} &
  2400 &
  603 &
  0 &
  2890 &
  160 &
  603 &
  0 &
  146 \\\hline
\multirow{3}{*}{MaaS-80-6} &
  1 &
  17456 &
  -$\infty$ &
  100\% &
  10800 &
  58 &
  728 &
  \textless{}$1e^{-6}$ &
  51 &
  5409 &
  728 &
  0 &
  5076 &
  51 &
  728 &
  \textless{}$1e^{-6}$ &
  \textbf{46} \\
 &
  2 &
  17668 &
  -$\infty$ &
  100\% &
  10800 &
  136 &
  621 &
  0 &
  \textbf{138} &
  595 &
  621 &
  0 &
  611 &
  196 &
  621 &
  \textless{}$1e^{-6}$ &
  190 \\
 &
  3 &
  18455 &
  -$\infty$ &
  100\% &
  10800 &
  250 &
  782 &
  0 &
  \textbf{261} &
  2281 &
  782 &
  0 &
  2172 &
  482 &
  782 &
  0 &
  484 \\\hline
\multirow{3}{*}{MaaS-90-6} &
  1 &
  18128 &
  -$\infty$ &
  100\% &
  10800 &
  210 &
  832 &
  0 &
  \textbf{222} &
  595 &
  832 &
  0 &
  10800 &
  595 &
  832 &
  0 &
  611 \\
 &
  2 &
  16887 &
  -$\infty$ &
  100\% &
  10800 &
  217 &
  890 &
  0 &
  \textbf{242} &
  2281 &
  890 &
  0 &
  10800 &
  2281 &
  890 &
  0 &
  2172 \\
 &
  3 &
  17473 &
  -$\infty$ &
  100\% &
  10800 &
  294 &
  999 &
  \textless{}$1e^{-6}$ &
  \textbf{311} &
  636 &
  999 &
  0 &
  10800 &
  636 &
  999 &
  \textless{}$1e^{-6}$ &
  682 \\\hline
\multirow{3}{*}{MaaS-100-6} &
  1 &
  17225 &
  -$\infty$ &
  100\% &
  10800 &
  236 &
  812 &
  0 &
  \textbf{246} &
  3332 &
  -$\infty$ &
  100\% &
  10800 &
  4286 &
  812 &
  0 &
  1431 \\
 &
  2 &
  17645 &
  -$\infty$ &
  100\% &
  9800 &
  261 &
  857 &
  0 &
  \textbf{255} &
  4286 &
  -$\infty$ &
  100\% &
  10800 &
  417 &
  857 &
  \textless{}$1e^{-6}$ &
  482 \\
 &
  3 &
  17331 &
  -$\infty$ &
  100\% &
  10800 &
  296 &
  922 &
  \textless{}$1e^{-6}$ &
  \textbf{306} &
  4017 &
  -$\infty$ &
  100\% &
  10800 &
  332 &
  922 &
  \textless{}$1e^{-6}$ &
  391 \\\hline
\multirow{3}{*}{MaaS-110-7} &
  1 &
  15436 &
  -$\infty$ &
  100\% &
  10800 &
  234 &
  716 &
  0 &
  \textbf{226} &
  3777 &
  -$\infty$ &
  100\% &
  10800 &
  1077 &
  716 &
  \textless{}$1e^{-6}$ &
  1314 \\
 &
  2 &
  15992 &
  -$\infty$ &
  100\% &
  10800 &
  329 &
  1212 &
  0 &
  \textbf{379} &
  3677 &
  -$\infty$ &
  100\% &
  10800 &
  967 &
  1212 &
  \textless{}$1e^{-6}$ &
  1228 \\
 &
  3 &
  15425 &
  -$\infty$ &
  100\% &
  10800 &
  276 &
  1102 &
  0 &
  \textbf{298} &
  3985 &
  -$\infty$ &
  100\% &
  10800 &
  1245 &
  1102 &
  \textless{}$1e^{-6}$ &
  1442 \\\hline
\multirow{3}{*}{MaaS-120-7} &
  1 &
  11475 &
  -$\infty$ &
  100\% &
  10800 &
  346 &
  1073 &
  0 &
  \textbf{456} &
  3716 &
  -$\infty$ &
  100\% &
  10800 &
  716 &
  1073 &
  \textless{}$1e^{-6}$ &
  1011 \\
 &
  2 &
  13576 &
  -$\infty$ &
  100\% &
  10800 &
  362 &
  1057 &
  \textless{}$1e^{-6}$ &
  \textbf{419} &
  3995 &
  -$\infty$ &
  100\% &
  10800 &
  3995 &
  1057 &
  \textless{}$1e^{-6}$ &
  5232 \\
 &
  3 &
  14247 &
  -$\infty$ &
  100\% &
  10800 &
  435 &
  944 &
  \textless{}$1e^{-6}$ &
  \textbf{540} &
  4079 &
  -$\infty$ &
  100\% &
  10800 &
  4079 &
  944 &
  \textless{}$1e^{-6}$ &
  5097 \\\hline
\multirow{3}{*}{MaaS-130-7} &
  1 &
  13582 &
  -$\infty$ &
  100\% &
  10800 &
  1030 &
  1034 &
  \textless{}$1e^{-6}$ &
  \textbf{1704} &
  3523 &
  -$\infty$ &
  100\% &
  10800 &
  1432 &
  1034 &
  \textless{}$1e^{-6}$ &
  2255 \\
 &
  2 &
  13559 &
  -$\infty$ &
  100\% &
  10800 &
  1342 &
  989 &
  0 &
  \textbf{2247} &
  3913 &
  -$\infty$ &
  100\% &
  10800 &
  1556 &
  989 &
  0 &
  2378 \\
 &
  3 &
  13902 &
  -$\infty$ &
  100\% &
  10800 &
  1021 &
  1120 &
  0 &
  \textbf{1689} &
  3803 &
  -$\infty$ &
  100\% &
  10800 &
  1641 &
  1120 &
  0 &
  2417 \\\hline
\multirow{3}{*}{MaaS-140-8} &
  1 &
  6441 &
  -$\infty$ &
  100\% &
  10800 &
  1114 &
  1492 &
  0 &
  \textbf{2273} &
  3272 &
  -$\infty$ &
  100\% &
  10800 &
  3408 &
  1492 &
  0 &
  5421 \\
 &
  2 &
  6712 &
  -$\infty$ &
  100\% &
  10800 &
  847 &
  1187 &
  \textless{}$1e^{-6}$ &
  \textbf{1607} &
  3843 &
  -$\infty$ &
  100\% &
  10800 &
  2900 &
  1187 &
  \textless{}$1e^{-6}$ &
  4300 \\
 &
  3 &
  6381 &
  -$\infty$ &
  100\% &
  10800 &
  586 &
  1216 &
  \textless{}$1e^{-6}$ &
  \textbf{1092} &
  3616 &
  -$\infty$ &
  100\% &
  10800 &
  2862 &
  1216 &
  \textless{}$1e^{-6}$ &
  4233 \\\hline
\multirow{3}{*}{MaaS-150-8} &
  1 &
  6300 &
  -$\infty$ &
  100\% &
  10800 &
  438 &
  1396 &
  0 &
  \textbf{882} &
  3425 &
  -$\infty$ &
  100\% &
  10800 &
  823 &
  1396 &
  0 &
  1381 \\
 &
  2 &
  6442 &
  -$\infty$ &
  100\% &
  10800 &
  738 &
  1159 &
  0 &
  \textbf{1531} &
  3433 &
  -$\infty$ &
  100\% &
  10800 &
  790 &
  1159 &
  0 &
  1401 \\
 &
  3 &
  6351 &
  -$\infty$ &
  100\% &
  10800 &
  1080 &
  1116 &
  \textless{}$1e^{-6}$ &
  \textbf{2362} &
  3512 &
  -$\infty$ &
  100\% &
  10800 &
  2692 &
  1116 &
  \textless{}$1e^{-6}$ &
  4622 \\\hline
\multirow{3}{*}{MaaS-160-8} &
  1 &
  5625 &
  -$\infty$ &
  100\% &
  10800 &
  625 &
  1341 &
  0 &
  \textbf{1233} &
  3246 &
  -$\infty$ &
  100\% &
  10800 &
  3022 &
  1341 &
  0 &
  5571 \\
 &
  2 &
  5352 &
  -$\infty$ &
  100\% &
  10800 &
  558 &
  1390 &
  0 &
  \textbf{1086} &
  3214 &
  -$\infty$ &
  100\% &
  10800 &
  5962 &
  1390 &
  0 &
  10249 \\
 &
  3 &
  5314 &
  -$\infty$ &
  100\% &
  10800 &
  532 &
  1281 &
  0 &
  \textbf{956} &
  3278 &
  -$\infty$ &
  100\% &
  10800 &
2945&
  1281 &
  0 &
  4493 \\\hline
\multirow{3}{*}{MaaS-170-9} &
  1 &
  5776 &
  -$\infty$ &
  100\% &
  10800 &
  1516 &
  1671 &
  \textless{}$1e^{-6}$ &
  \textbf{3063} &
  3026 &
  -$\infty$ &
  100\% &
  10800 &
  3869 &
  1671 &
  \textless{}$1e^{-6}$ &
  7118 \\
 &
  2 &
  5694 &
  -$\infty$ &
  100\% &
  10800 &
  858 &
  1744 &
  0 &
  \textbf{1688} &
  2937 &
  -$\infty$ &
  100\% &
  10800 &
  2614 &
  1744 &
  0 &
  4174 \\
 &
  3 &
  5761 &
  -$\infty$ &
  100\% &
  10800 &
  602 &
  1648 &
  0 &
  \textbf{1224} &
  3141 &
  -$\infty$ &
  100\% &
  10800 &
  834 &
  1648 &
  0 &
  1516 \\\hline
\multirow{3}{*}{MaaS-180-9} &
  1 &
  4953 &
  -$\infty$ &
  100\% &
  10800 &
  892 &
  1953 &
  0 &
  \textbf{2088} &
  2451 &
  -$\infty$ &
  100\% &
  10800 &
  1019 &
  1953 &
  0 &
  2016 \\
 &
  2 &
  5380 &
  -$\infty$ &
  100\% &
  10800 &
  648 &
  1709 &
  0 &
  \textbf{1390} &
  2795 &
  -$\infty$ &
  100\% &
  10800 &
  4041 &
  1709 &
  \textless{}$1e^{-6}$ &
  7337 \\
 &
  3 &
  5338 &
  -$\infty$ &
  100\% &
  10800 &
  895 &
  1607 &
  0 &
  \textbf{1916} &
  3311 &
  -$\infty$ &
  100\% &
  10800 &
  2140 &
  1607 &
  \textless{}$1e^{-6}$ &
  3522 \\\hline
\multirow{3}{*}{MaaS-190-10} &
  1 &
  4786 &
  -$\infty$ &
  100\% &
  10800 &
  438 &
  1784 &
  0 &
  \textbf{1028} &
  2240 &
  -$\infty$ &
  100\% &
  10800 &
  3917 &
  1784 &
  \textless{}$1e^{-6}$ &
  9800 \\
 &
  2 &
  4955 &
  -$\infty$ &
  100\% &
  10800 &
  519 &
  1702 &
  0 &
  \textbf{1141} &
  2421 &
  -$\infty$ &
  100\% &
  10800 &
  2233 &
  1702 &
  \textless{}$1e^{-6}$ &
  3882 \\
 &
  3 &
  4862 &
  -$\infty$ &
  100\% &
  10800 &
  560 &
  1787 &
  0 &
  \textbf{1325} &
  2322 &
  -$\infty$ &
  100\% &
  10800 &
  1641 &
  1787 &
  \textless{}$1e^{-6}$ &
  3186 \\\hline
\multirow{3}{*}{MaaS-200-10} &
  1 &
  4632 &
  -$\infty$ &
  100\% &
  10800 &
  1351 &
  1660 &
  0 &
  \textbf{3155} &
  2230 &
  -$\infty$ &
  100\% &
  10800 &
  2743 &
  1660 &
  \textless{}$1e^{-6}$ &
  9116 \\
 &
  2 &
  4616 &
  -$\infty$ &
  100\% &
  10800 &
  2894 &
  1643 &
  0 &
  \textbf{7572} &
  2688 &
  -$\infty$ &
  100\% &
  10800 &
 3656 &
  1643 &
  \textless{}$1e^{-6}$ &
  10262 \\
 &
  3 &
  4457 &
  -$\infty$ &
  100\% &
  10800 &
  877 &
  1847 &
  0 &
  \textbf{2212} &
  2867 &
  -$\infty$ &
  100\% &
  10800 &
  3923 &
  1847 &
  \textless{}$1e^{-6}$ &
  10703\\ \bottomrule
\end{tabular}
{Numbers in bold denotes the smallest CPU runtime for each instance }
\end{table}
\setlength{\abovecaptionskip}{0pt}
\setlength{\belowcaptionskip}{0pt}
\begin{table}
\centering
\scriptsize
{Performance of the Branch and Bound algorithms for \SLMFGQ\space(MIQBP)\label{AQtable}}
{{\begin{tabular}{@{}cccccccccccccccccc@{}}
\toprule
Instances &
   &
  \multicolumn{4}{c}{Bard\&Moore} &
  \multicolumn{4}{c}{SD-BP} &
  \multicolumn{4}{c}{SD-Diffob} &
  \multicolumn{4}{c}{SD-Wi} \\ \cmidrule(l){3-6} \cmidrule(l){7-10} \cmidrule(l){11-14} \cmidrule(l){15-18}
 &
   &
  \textit{k} &
  \textit{LB} &
  \textit{Gap(\%) } &
  \textit{$T (s)$ } &
  \textit{k} &
  \textit{LB} &
  \textit{Gap(\%)} &
  \textit{$T (s)$ } &
  \textit{k} &
  \textit{LB} &
  \textit{Gap(\%)} &
  \textit{$T (s)$ } &
  \textit{k} &
  \textit{LB} &
  \textit{Gap(\%)} &
  \textit{$T (s)$ } \\\cmidrule(l){3-6} \cmidrule(l){7-10} \cmidrule(l){11-14} \cmidrule(l){15-18}
\multirow{3}{*}{MaaS-10-5} &
  1 &
  52 &
  19 &
  0 &
  9 &
  30 &
  19 &
  0 &
  5 &
  26 &
  19 &
  0 &
  7 &
  26 &
  19 &
  0 &
  \textbf{6} \\
 &
  2 &
  52 &
  -37 &
  0 &
  9 &
  29 &
  -37 &
  0 &
  5 &
  25 &
  -37 &
  0 &
  7 &
  25 &
  -37 &
  0 &
  \textbf{6} \\
 &
  3 &
  56 &
  -42 &
  0 &
  9 &
  33 &
  -42 &
  0 &
  6 &
  27 &
  -42 &
  0 &
  7 &
  27 &
  -42 &
  0 &
  \textbf{6} \\\hline
\multirow{3}{*}{MaaS-20-5} &
  1 &
  100 &
  226 &
  0 &
  19 &
  22 &
  226 &
  0 &
  6 &
  20 &
  226 &
  0 &
  \textbf{6} &
  20 &
  226 &
  0 &
  7 \\
 &
  2 &
  136 &
  40 &
  0 &
  25 &
  44 &
  40 &
  0 &
  \textbf{10} &
  86 &
  40 &
  0 &
  25 &
  44 &
  40 &
  0 &
  14 \\
 &
  3 &
  112 &
  39 &
  0 &
  21 &
  53 &
  39 &
  0 &
  \textbf{13} &
  52 &
  39 &
  0 &
  14 &
  46 &
  39 &
  0 &
  14 \\\hline
\multirow{3}{*}{MaaS-30-5} &
  1 &
  44 &
  310 &
  0 &
  10 &
  19 &
  310 &
  0 &
  \textbf{7} &
  17 &
  310 &
  0 &
  \textbf{7} &
  16 &
  310 &
  0 &
  9 \\
 &
  2 &
  392 &
  301 &
  0 &
  95 &
  38 &
  301 &
  0 &
  \textbf{11} &
  186 &
  301 &
  0 &
  67 &
  35 &
  301 &
  0 &
  15 \\
 &
  3 &
  216 &
  280 &
  \textless{}$1e^{-6}$ &
  51 &
  98 &
  280 &
  \textless{}$1e^{-6}$ &
  30 &
  96 &
  280 &
  \textless{}$1e^{-6}$ &
  39 &
  70 &
  280 &
  \textless{}$1e^{-6}$ &
  \textbf{28} \\\hline
\multirow{3}{*}{MaaS-40-5} &
  1 &
  268 &
  343 &
  \textless{}$1e^{-6}$ &
  92 &
  65 &
  343 &
  \textless{}$1e^{-6}$ &
  \textbf{22} &
  114 &
  343 &
  \textless{}$1e^{-6}$ &
  48 &
  59 &
  343 &
  \textless{}$1e^{-6}$ &
  27 \\
 &
  2 &
  192 &
  216 &
  0 &
  66 &
  65 &
  216 &
  0 &
  \textbf{24} &
  96 &
  216 &
  0 &
  39 &
  59 &
  216 &
  0 &
  27 \\
 &
  3 &
  216 &
  272 &
  0 &
  75 &
  81 &
  272 &
  0 &
  \textbf{28} &
  114 &
  272 &
  0 &
  48 &
  68 &
  272 &
  0 &
  30 \\\hline
\multirow{3}{*}{MaaS-50-5} &
  1 &
  300 &
  127 &
  0 &
  112 &
  72 &
  127 &
  0 &
  \textbf{37} &
  209 &
  127 &
  0 &
  94 &
  128 &
  127 &
  0 &
  63 \\
 &
  2 &
  120 &
  226 &
  0 &
  46 &
  67 &
  226 &
  0 &
  \textbf{27} &
  81 &
  226 &
  0 &
  40 &
  59 &
  226 &
  0 &
  30 \\
 &
  3 &
  76 &
  578 &
  \textless{}$1e^{-6}$ &
  29 &
  38 &
  578 &
  \textless{}$1e^{-6}$ &
  \textbf{15} &
  38 &
  578 &
  \textless{}$1e^{-6}$ &
  20 &
  38 &
  578 &
  \textless{}$1e^{-6}$ &
  20 \\\hline
\multirow{3}{*}{MaaS-60-5} &
  1 &
  2900 &
  482 &
  0 &
  1508 &
  186 &
  482 &
  0 &
  \textbf{110} &
  618 &
  482 &
  0 &
  421 &
  220 &
  482 &
  0 &
  172 \\
 &
  2 &
  2700 &
  468 &
  0 &
  1080 &
  126 &
  468 &
  0 &
  \textbf{64} &
  1152 &
  468 &
  0 &
  721 &
  200 &
  468 &
  0 &
  139 \\
 &
  3 &
  628 &
  382 &
  0 &
  279 &
  136 &
  382 &
  0 &
  \textbf{87} &
  215 &
  382 &
  0 &
  127 &
  186 &
  382 &
  0 &
  150 \\\hline
\multirow{3}{*}{MaaS-70-6} &
  1 &
  21113 &
  -$\infty$ &
  100\% &
  10800 &
  122 &
  492 &
  0 &
  \textbf{96} &
  16804 &
  492 &
  0 &
  10714 &
  158 &
  492 &
  0 &
  138 \\
 &
  2 &
  22087 &
  -$\infty$ &
  100\% &
  10800 &
  178 &
  675 &
  \textless{}$1e^{-6}$ &
  \textbf{136} &
  10690 &
  675 &
  0 &
  6699 &
  305 &
  675 &
  \textless{}$1e^{-6}$ &
  296 \\
 &
  3 &
  20782 &
  -$\infty$ &
  100\% &
  10800 &
  243 &
  535 &
  \textless{}$1e^{-6}$ &
  \textbf{203} &
  12836 &
  535 &
  0 &
  10728 &
  216 &
  535 &
  \textless{}$1e^{-6}$ &
  208 \\\hline
\multirow{3}{*}{MaaS-80-6} &
  1 &
  14358 &
  693 &
  \textless{}$1e^{-6}$ &
  9800 &
  114 &
  693 &
  0 &
  \textbf{119} &
  8594 &
  693 &
  0 &
  9801 &
  146 &
  693 &
  0 &
  168 \\
 &
  2 &
  15874 &
  -$\infty$ &
  100\% &
  10800 &
  152 &
  583 &
  0 &
  \textbf{167} &
  10110 &
  -$\infty$ &
  100\% &
  10800 &
  220 &
  583 &
  0 &
  224 \\
 &
  3 &
  11112 &
  710 &
  \textless{}$1e^{-6}$ &
  10800 &
  160 &
  710 &
  0 &
  \textbf{155} &
  9953 &
  710 &
  \textless{}$1e^{-6}$ &
  9901 &
  290 &
  710 &
  0 &
  279 \\\hline
\multirow{3}{*}{MaaS-90-6} &
  1 &
  16609 &
  -$\infty$ &
  100\% &
  10800 &
  104 &
  760 &
  0 &
  \textbf{104} &
  10031 &
  -$\infty$ &
  100\% &
  10800 &
  153 &
  760 &
  0 &
  144 \\
 &
  2 &
  18945 &
  -$\infty$ &
  100\% &
  10800 &
  126 &
  827 &
  0 &
  \textbf{166} &
  6390 &
  827 &
  \textless{}$1e^{-6}$ &
  6120 &
  234 &
  827 &
  0 &
  217 \\
 &
  3 &
  17824 &
  -$\infty$ &
  100\% &
  10800 &
  144 &
  917 &
  \textless{}$1e^{-6}$ &
  \textbf{160} &
  10756 &
  -$\infty$ &
  100\% &
  10800 &
  402 &
  917 &
  \textless{}$1e^{-6}$ &
  396 \\\hline
\multirow{3}{*}{MaaS-100-6} &
  1 &
  14435 &
  -$\infty$ &
  100\% &
  10800 &
  162 &
  758 &
  0 &
  \textbf{190} &
  9535 &
  -$\infty$  &
  100\% &
  10800 &
  585 &
  758 &
  0 &
  669 \\
 &
  2 &
  16064 &
  -$\infty$ &
  100\% &
  10800 &
  243 &
  806 &
  \textless{}$1e^{-6}$ &
  \textbf{274} &
  9422 &
  806 &
  0 &
  10719 &
  229 &
  806 &
  \textless{}$1e^{-6}$ &
  219 \\
 &
  3 &
  16412 &
  -$\infty$ &
  100\% &
  10800 &
  326 &
  882 &
  0 &
  \textbf{377} &
  8988 &
  882 &
  $1e^{-6}$  &
  10602 &
  433 &
  882 &
  0 &
  445 \\\hline
\multirow{3}{*}{MaaS-110-7} &
  1 &
  15194 &
  -$\infty$ &
  100\% &
  10800 &
  147 &
  668 &
  0 &
  \textbf{169} &
  8665 &
  -$\infty$ &
  100\% &
  10800 &
  715 &
  668 &
  0 &
  869 \\
 &
  2 &
  14945 &
  -$\infty$ &
  100\% &
  10800 &
  497 &
  1107 &
  0 &
  \textbf{694} &
  9237 &
  -$\infty$ &
  100\% &
  10800 &
  866 &
  1107 &
  0 &
  1189 \\
 &
  3 &
  14952 &
  -$\infty$ &
  100\% &
  10800 &
  128 &
  1055 &
  0 &
  \textbf{173} &
  9114 &
  -$\infty$ &
  100\% &
  10800 &
  928 &
  1055 &
  0 &
  1099 \\\hline
\multirow{3}{*}{MaaS-120-7} &
  1 &
  13923 &
  -$\infty$ &
  100\% &
  10800 &
  166 &
  1028 &
  0 &
  \textbf{246} &
  8286 &
  -$\infty$ &
  100\% &
  10800 &
  860 &
  1028 &
  0 &
  1205 \\
 &
  2 &
  14541 &
  -$\infty$ &
  100\% &
  10800 &
  109 &
  1026 &
  \textless{}$1e^{-6}$ &
  \textbf{134} &
  8450 &
  -$\infty$ &
  100\% &
  10800 &
  887 &
  1026 &
  \textless{}$1e^{-6}$ &
  1163 \\
 &
  3 &
  14169 &
  -$\infty$ &
  100\% &
  10800 &
  129 &
  918 &
  0 &
  \textbf{186} &
  9212 &
  -$\infty$ &
  100\% &
  10800 &
  347 &
  918 &
  0 &
  417 \\\hline
\multirow{3}{*}{MaaS-130-7} &
  1 &
  9219 &
  -$\infty$ &
  100\% &
  10800 &
  626 &
  945 &
  0 &
  \textbf{977} &
  5610 &
  -$\infty$ &
  100\% &
  10800 &
  1762 &
  945 &
  0 &
  2321 \\
 &
  2 &
  12685 &
  -$\infty$ &
  100\% &
  10800 &
  204 &
  949 &
  0 &
  \textbf{284} &
  8036 &
  -$\infty$ &
  100\% &
  10800 &
  1670 &
  949 &
  0 &
  2284 \\
 &
  3 &
  13396 &
  -$\infty$ &
  100\% &
  10800 &
  223 &
  1076 &
  0 &
  \textbf{336} &
  8154 &
  -$\infty$ &
  100\% &
  10800 &
  2045 &
  1076 &
  0 &
  2732 \\\hline
\multirow{3}{*}{MaaS-140-8} &
  1 &
  12631 &
  -$\infty$ &
  100\% &
  10800 &
  242 &
  1420 &
  0 &
  \textbf{449} &
  8036 &
  -$\infty$ &
  100\% &
  10800 &
  1170 &
  1420 &
  0 &
  2236 \\
 &
  2 &
  10485 &
  -$\infty$ &
  100\% &
  10800 &
  323 &
  1148 &
  0 &
  \textbf{592} &
  8154 &
  -$\infty$ &
  100\% &
  10800 &
  656 &
  1148 &
  0 &
  1033 \\
 &
  3 &
  6418 &
  -$\infty$ &
  100\% &
  10800 &
  256 &
  1162 &
  0 &
  \textbf{457} &
  3229 &
  -$\infty$ &
  100\% &
  10800 &
  2974 &
  1162 &
  0 &
  4633 \\\hline
\multirow{3}{*}{MaaS-150-8} &
  1 &
  6035 &
  -$\infty$ &
  100\% &
  10800 &
  166 &
  1351 &
  0 &
  \textbf{314} &
  3001 &
  -$\infty$ &
  100\% &
  10800 &
  757 &
  1351 &
  0 &
  1321 \\
 &
  2 &
  6145 &
  -$\infty$ &
  100\% &
  10800 &
  152 &
  1127 &
  0 &
  \textbf{300} &
  3056 &
  -$\infty$ &
  100\% &
  10800 &
  2042 &
  1127 &
  0 &
  2942 \\
 &
  3 &
  6363 &
  -$\infty$ &
  100\% &
  10800 &
  870 &
  1046 &
  0 &
  \textbf{1788} &
  3548 &
  -$\infty$ &
  100\% &
  10800 &
  2276 &
  1046 &
  0 &
  3144 \\\hline
\multirow{3}{*}{MaaS-160-8} &
  1 &
  5461 &
  -$\infty$ &
  100\% &
  10800 &
  286 &
  1282 &
  \textless{}$1e^{-6}$ &
  \textbf{564} &
  2799 &
  -$\infty$ &
  100\% &
  10800 &
  688 &
  1282 &
  \textless{}$1e^{-6}$ &
  1198 \\
 &
  2 &
  5504 &
  -$\infty$ &
  100\% &
  10800 &
  300 &
  1338 &
  \textless{}$1e^{-6}$ &
  \textbf{583} &
  3014 &
  -$\infty$ &
  100\% &
  10800 &
  1944 &
  1338 &
  \textless{}$1e^{-6}$ &
  4908 \\
 &
  3 &
  6213 &
  -$\infty$ &
  100\% &
  10800 &
  428 &
  1228 &
  \textless{}$1e^{-6}$ &
  \textbf{792} &
  3194 &
  -$\infty$ &
  100\% &
  10800 &
  757 &
  1228 &
  \textless{}$1e^{-6}$ &
  1321 \\\hline
\multirow{3}{*}{MaaS-170-9} &
  1 &
  5343 &
  -$\infty$ &
  100\% &
  10800 &
  140 &
  1643 &
  0 &
  \textbf{267} &
  2712 &
  -$\infty$ &
  100\% &
  10800 &
  661 &
  1643 &
  0 &
  1340 \\
 &
  2 &
  5325 &
  -$\infty$ &
  100\% &
  10800 &
  190 &
  1694 &
  0 &
  \textbf{323} &
  2702 &
  -$\infty$ &
  100\% &
  10800 &
 888 &
  1694 &
  0 &
  1877 \\
 &
  3 &
  5564 &
  -$\infty$ &
  100\% &
  10800 &
  200 &
  1605 &
  0 &
  \textbf{411} &
  2709 &
  -$\infty$ &
  100\% &
  10800 &
  749 &
  1605 &
  0 &
  1481 \\\hline
\multirow{3}{*}{MaaS-180-9} &
  1 &
  4860 &
  -$\infty$ &
  100\% &
  10800 &
  149 &
  1893 &
  0 &
  \textbf{303} &
  2305 &
  -$\infty$ &
  100\% &
  10800 &
  733 &
  1893 &
  0 &
  1546 \\
 &
  2 &
  5213 &
  -$\infty$ &
  100\% &
  10800 &
  282 &
  1663 &
  0 &
  \textbf{596} &
  2504 &
  -$\infty$ &
  100\% &
  10800 &
  2084 &
  1663 &
  0 &
  5144 \\
 &
  3 &
  5126 &
  -$\infty$ &
  100\% &
  10800 &
  419 &
  1554 &
  \textless{}$1e^{-6}$ &
  \textbf{772} &
  2725 &
  -$\infty$ &
  100\% &
  10800 &
  1612 &
  1554 &
  \textless{}$1e^{-6}$ &
  3159 \\\hline
\multirow{3}{*}{MaaS-190-10} &
  1 &
  4728 &
  -$\infty$ &
  100\% &
  10800 &
  364 &
  1724 &
  \textless{}$1e^{-6}$ &
  \textbf{855} &
  2093 &
  -$\infty$ &
  100\% &
  10800 &
  1036 &
  1724 &
  \textless{}$1e^{-6}$ &
  2363 \\
 &
  2 &
  5104 &
  -$\infty$ &
  100\% &
  10800 &
  231 &
  1673 &
  \textless{}$1e^{-6}$ &
  \textbf{491} &
  2143 &
  -$\infty$ &
  100\% &
  10800 &
  1374 &
  1673 &
  \textless{}$1e^{-6}$ &
  2965 \\
 &
  3 &
  4889 &
  -$\infty$ &
  100\% &
  10800 &
  270 &
  1730 &
  \textless{}$1e^{-6}$ &
  \textbf{626} &
  2246 &
  -$\infty$ &
  100\% &
  10800 &
  1297 &
  1730 &
  \textless{}$1e^{-6}$ &
  2736 \\\hline
\multirow{3}{*}{MaaS-200-10} &
  1 &
  4741 &
  -$\infty$ &
  100\% &
  10800 &
  565 &
  1621 &
  \textless{}$1e^{-6}$ &
  \textbf{1230} &
  2202 &
  -$\infty$ &
  100\% &
  10800 &
  2487 &
  1621 &
  \textless{}$1e^{-6}$ &
  5196\\
 &
  2 &
  4607 &
  -$\infty$ &
  100\% &
  10800 &
  256 &
  1584 &
  \textless{}$1e^{-6}$ &
  \textbf{561} &
  1986 &
  -$\infty$ &
  100\% &
  10800 &
  2474 &
  1584 &
  \textless{}$1e^{-6}$ &
 4783 \\
 &
  3 &
  4880 &
  -$\infty$ &
  100\% &
  10800 &
  509 &
  1814 &
  \textless{}$1e^{-6}$ &
  \textbf{1194} &
  1974 &
  -$\infty$ &
  100\% &
  10800 &
  2480 &
  1814 &
  \textless{}$1e^{-6}$ &
  4962 \\ \bottomrule
\end{tabular}}}	
{Numbers in bold denotes the smallest CPU runtime for each instance}
\end{table}
\end{document}